\documentclass[reqno]{amsart}
\usepackage[margin = 1in]{geometry}
\usepackage{amsmath, amssymb, amsthm, fancyhdr, verbatim, graphicx}
\usepackage{enumerate}
\usepackage[all]{xy}
\usepackage[usenames,dvipsnames]{xcolor}
\usepackage{mathrsfs}
\usepackage{tikz-cd}
\usepackage{framed, hyperref}
\usepackage[titletoc]{appendix}
\usepackage{bbm}
\usepackage{lipsum}
\usepackage{adjustbox}
\usepackage{bbm}

\usepackage{stmaryrd}

\numberwithin{equation}{section}

\newtheorem{theorem}{Theorem}
\newtheorem{lemma}[theorem]{Lemma}

\newtheorem{proposition}[theorem]{Proposition}
\newtheorem{corollary}[theorem]{Corollary}
\newtheorem{conjecture}[theorem]{Conjecture}

\theoremstyle{remark}
\newtheorem{definition}[theorem]{Definition}
\newtheorem{examplen}[theorem]{Example}
\newtheorem*{example}{Example}
\newtheorem{remark}[theorem]{Remark}

\numberwithin{theorem}{section}
\numberwithin{equation}{section}

\DeclareMathOperator{\Ad}{Ad}
\DeclareMathOperator{\Bun}{Bun}
\DeclareMathOperator{\Aut}{Aut}

\DeclareMathOperator{\GL}{GL}

\DeclareMathOperator{\PGL}{PGL}

\DeclareMathOperator{\depth}{depth}

\DeclareMathOperator{\Def}{Def}

\DeclareMathOperator{\Frob}{Frob}

\DeclareMathOperator{\Hom}{Hom}

\DeclareMathOperator{\rank}{rank}

\DeclareMathOperator{\Sets}{Sets}

\DeclareMathOperator{\Gal}{Gal}

\DeclareMathOperator{\ad}{ad}
\DeclareMathOperator{\diag}{diag}

\DeclareMathOperator{\End}{End}
\DeclareMathOperator{\tr}{tr}
\DeclareMathOperator{\pr}{pr}
\DeclareMathOperator{\Lie}{Lie}

\DeclareMathOperator{\Res}{Res}

\DeclareMathOperator{\Ind}{Ind}

\DeclareMathOperator{\res}{res}
\DeclareMathOperator{\Fr}{Fr}

\DeclareMathOperator{\rec}{rec}

\DeclareMathOperator{\Spec}{Spec}

\DeclareMathOperator{\IC}{IC}
\DeclareMathOperator{\CT}{CT}

\newcommand{\GG}{\mathbb{G}}
\newcommand{\TT}{\mathbb{T}}

\newcommand{\UU}{\mathbb{U}}
\newcommand{\rT}{\mathrm{T}}
\newcommand{\F}{\mathbb{F}}
\newcommand{\co}{\colon}
\newcommand{\ul}[1]{\underline{#1}}
\newcommand{\CC}{\mathbf{C}}

\newcommand{\PP}{\mathbf{P}}

\newcommand{\A}{\mathbb{A}}
\def\acts{\curvearrowright}
\newcommand{\surj}{\twoheadrightarrow}
\newcommand{\msf}[1]{\mathsf{#1}}

\DeclareMathOperator{\Nm}{Nm}

\newcommand{\cB}{{\mathcal B}}
\newcommand{\cC}{{\mathcal C}}
\newcommand{\cD}{{\mathcal D}}
\newcommand{\cE}{{\mathcal E}}
\newcommand{\cF}{{\mathcal F}}
\newcommand{\cG}{{\mathcal G}}
\newcommand{\cH}{{\mathcal H}}

\newcommand{\cJ}{{\mathcal J}}

\newcommand{\cL}{{\mathcal L}}

\newcommand{\cO}{{\mathcal O}}

\newcommand{\cW}{{\mathcal W}}

\newcommand{\cZ}{{\mathcal Z}}
\newcommand{\fra}{{\mathfrak a}}
\newcommand{\frb}{{\mathfrak b}}

\newcommand{\ra}{{\rightarrow}}

\newcommand{\frg}{{\mathfrak g}}
\newcommand{\frh}{{\mathfrak h}}

\newcommand{\ffrm}{{\mathfrak m}}
\newcommand{\frn}{{\mathfrak n}}

\newcommand{\bbA}{{\mathbb A}}

\newcommand{\bbF}{{\mathbb F}}

\newcommand{\bbQ}{{\mathbb Q}}

\newcommand{\bbZ}{{\mathbb Z}}

\newcommand{\hG}{{\widehat{G}}}
\newcommand{\G}{{\Gamma}}
\newcommand{\Gv}{{\Gamma_{K_v}}}

\newcommand{\hT}{{\widehat{T}}}
\newcommand{\bs}{{\overline{\sigma}}}

\newcommand{\wh}[1]{\widehat{#1}}
\newcommand{\mf}[1]{\mathfrak{#1}}

\DeclareMathOperator{\cInd}{c-Ind}
\DeclareMathOperator{\pt}{pt}

\newcommand{\nc}{\newcommand}
\nc{\renc}{\renewcommand}
\nc{\ssec}{\subsection}
\nc{\sssec}{\subsubsection}
\nc{\on}{\operatorname}

\nc\ol{\overline}
\nc\wt{\widetilde}
\nc\tboxtimes{\wt{\boxtimes}}
\nc{\alp}{\alpha}

\nc{\BunBb}{\overline{\Bun}_B}

\newcommand{\Z}{\mathbb{Z}}
\newcommand{\Q}{\mathbb{Q}}
\newcommand{\inj}{\hookrightarrow}
\newcommand{\pc}{{}^{\mf{p}}}
\newcommand{\pch}{{}^{\mf{p}} \cH}
\newcommand{\ld}{{}^L}

\newcommand{\mrm}[1]{\mathrm{#1}}

\DeclareMathOperator{\Sht}{Sht}
\DeclareMathOperator{\gen}{gen}
\DeclareMathOperator{\cusp}{cusp}

\DeclareMathOperator{\Tr}{Tr}
\DeclareMathOperator{\Ima}{Im}
\DeclareMathOperator{\Jord}{Jord}
\DeclareMathOperator{\Exc}{Exc}

\DeclareMathOperator{\tors}{tors}

\title[Cyclic base change over function fields]{Cyclic base change of cuspidal automorphic representations over function fields}
\author[G. B\"ockle, T. Feng, M. Harris, C. Khare, and J. A. Thorne]{Gebhard B\"ockle, Tony Feng, Michael Harris, \\Chandrashekhar Khare, and Jack A. Thorne}

\begin{document}

\begin{abstract}
Let $G$ be a split semi-simple group over a global function field $K$. Given a cuspidal automorphic representation $\Pi$ of $G$ satisfying a technical hypothesis, we prove that for almost all primes $\ell$, there is a cyclic base change lifting of $\Pi$ along any $\Z/\ell\Z$-extension of $K$. Our proof does not rely on any trace formulas; instead it is based on modularity lifting theorems, together with a Smith theory argument to obtain base change for residual representations. As an application, we also prove that for any split semisimple group $G$ over a local function field $F$, and almost all primes $\ell$, any irreducible admissible representation of $G(F)$ admits a base change along any $\Z/\ell \Z$-extension of $F$. Finally, we characterize local base change more explicitly for a class of representations called toral supercuspidal representations. \end{abstract}

\maketitle

\setcounter{tocdepth}{1}
\tableofcontents

\section{Introduction}






A particular case of Langlands' principle of functoriality is \emph{cyclic base change}, which has been established in significant generality over global fields of characteristic $0$, but \emph{not} global fields of characteristic $p$. The results towards cyclic base change functoriality in characteristic 0 were worked out in a long series of papers by many authors including Saito, Shintani, Langlands, Arthur, Clozel, Labesse, and Harris (see the introduction to \cite{F20} for references), all following a strategy proposed by Langlands that is based on studying the twisted trace formula. But this strategy encounters complications in characteristic $p$, for example because of the contribution of inseparable elements to the trace formula. 

In this paper we establish (under technical assumptions) cyclic base change liftings over characteristic $p$ global fields, using a completely different strategy that is instead predicated upon the theory of automorphy lifting. Furthermore, we prove novel types of results regarding cyclic base change in the Genestier-Lafforgue correspondence. The results and proofs will be described in more detail in the rest of the introduction. 	

\subsection{Global base change}

Let $X$ be a smooth, projective, geometrically connected curve over the finite field $\bbF_q$ with function field $K := \bbF_q(X)$, and let $G$ be a split semisimple algebraic group over $\bbF_q$. Let $\Pi$ be a cuspidal automorphic representation of $G(\A_K)$ defined over a number field\footnote{Every cuspidal automorphic representation over a field of characteristic zero admits a model over a number field, so this assumption incurs no loss of generality for the questions we will study.}. For a field extension $K'/K$, a \emph{(weak) base change lifting} of $\Pi$ to $K'$ is a cuspidal automorphic representation $\Pi'$ of $G(\A_{K'})$ such that for almost all places $v$ of $K$ where $\Pi$ is unramified, $\Pi'$ is also unramified at all places of $K'$ above $v$, and the Satake parameters of $\Pi$ and of $\Pi'$ over $v$ are related by the \emph{base change homomorphism} for spherical Hecke algebras. (Later in the paper we will define and prove stronger notions of base change lifting.)

Over number fields, the existence of such a base change lifting is proved in complete generality for cuspidal automorphic representations of $\GL_n$ \cite{Art89}, and it is proved under some technical assumptions by Labesse \cite{Lab00} for general groups, still over number fields, by a comparison of the stable trace formula for $G(\bbA_{K})$ with the  stable trace formula for $G(\bbA_{K'})$ twisted with respect to a generator of $\Gal(K'/K)$.  It is expected that the stable trace formula will eventually provide weak base change over function fields as well, but the proof of such a formula faces a number of obstacles and even the first steps of the proof of the general trace formula are only now becoming available \cite{LL21}.  

In this paper we take a different approach to establish the existence of weak base change liftings (under technical hypotheses), based on automorphy lifting arguments instead of the trace formula. In particular, we use the results of Vincent Lafforgue on the global Langlands correspondence \cite{Laf18}. Let $\hG$ denote the dual group of $G$, considered as a split semisimple group scheme over $\bbZ$. For each prime $\ell$ we fix an algebraic closure $\overline{\bbQ}_\ell$ of $\bbQ_\ell$. As a consequence of V. Lafforgue's work in \cite{Laf18} is the construction, for each cuspidal automorphic representation $\Pi$ of $G(\A_K)$ valued in $\ol{\Q}_{\ell}$, of a finite set of $\hG(\ol{\Q}_{\ell})$-conjugacy classes of continuous Galois representations $\sigma : \Gal(K^s/K) \to \hG(\ol{\Q}_{\ell})$, each of which has the property that its Frobenius eigenvalues at unramified places match the Satake parameters of $\Pi$ (see \S \ref{sec: automorphy lifting} for the precise formulation). 

Suppose $\Pi$ is a cuspidal automorphic representation defined over $\ol{\Q}$. Then for each $\ell \neq p$, we may choose an embedding $\ol{\Q} \inj \ol{\Q}_{\ell}$, and then apply Lafforgue's theory, obtaining as above a finite set of Galois representations. These will be referred to as the Galois representations ``attached to $\Pi$ by Lafforgue's correspondence''. A consequence of our main result is the following:  

\begin{theorem}\label{thm: intro BC}
Let $\Pi$ be a cuspidal automorphic representation of $G(\A_K)$ over $\ol{\Q}$. Suppose that for some (equivalently, any) prime $\ell' \neq p$, some (equivalently, each) of the Galois representations attached to $\Pi$ by Lafforgue's correspondence has Zariski-dense image in $\wh{G}$. Then there exists a constant $c(\Pi)$ such that for all primes $\ell > c(\Pi)$, and all cyclic $\ell$-extensions $K'/K$, there is a weak base change lifting of $\Pi$ to a cuspidal automorphic representation $\Pi'$ of $G(\A_{K'})$. 
\end{theorem}

This follows from a more precise and refined result, Theorem \ref{thm: global existence}, which guarantees that $\Pi'$ may be chosen to be a base change lifting in the strongest possible sense (including a compatibility at ramified places). However, the latter requires more discussion in order to formulate, so we postpone it for now.

\begin{remark}
In \cite{ST21}, the Ramanujan Conjecture is proved for cuspidal automorphic representations satisfying certain types of local conditions, \emph{conditionally} upon the existence of base change liftings (with compatibility at ramified places as well) for constant extensions of large enough degree. Our Theorem \ref{thm: global existence} provides this type of base change lifting. However, it is not hard to deduce the Ramanujan Conjecture directly when our hypothesis is satisfied. 
\end{remark}

We sketch the approach of the proof. It is based on the automorphy lifting techniques pioneered by Taylor-Wiles. The first ingredient in this program is residual automorphy, and this is provided by \cite{F20}, which established the existence of base change for automorphic forms over $\ol{\F}_{\ell}$ in the special case where the extension $K'/K$ is cyclic of order $\ell$. (It was essential for the argument of \cite{F20} that the characteristic $\ell$ of the automorphic forms coincides with the order $\ell$ of the extension.) To prove Theorem \ref{thm: intro BC}, we use \emph{compatible systems} of Galois representations (whose existence is guaranteed by the Zariski density assumption) to bootstrap from this case and lift the automorphy to characteristic zero. In particular, we prove an automorphy lifting theorem for representations with arbitrary ramification, going beyond the everywhere unramified case which was established by four of the authors in \cite{BHKT}. 

We briefly discuss some of the new issues that go into proving this ramified automorphy lifting result. First, we need to arrange local deformation rings that are well-behaved, which is one reason for requiring $\ell > c(\Pi)$. We show that such a condition guarantees that the unrestricted local deformation condition is smooth (see \S \ref{mwf} for discussion of this). A second issue is that we need to know that the image of the associated $\ell$-adic Galois representation is large enough to provide a sufficient supply of Taylor-Wiles primes, for all $\ell > c(\Pi)$.  The condition used in \cite{BHKT} is that the image is {\it $\hG$-abundant}. For the ``potential automorphy'' theorem of \cite{BHKT}, it was enough to know that the image mod $\ell$ is abundant outside a set of $\ell$ having Dirichlet density 0.  Since the $\ell$-adic image is Zariski dense, this follows from a result of Larsen \cite{Lar95}. In order to obtain base change liftings for $\ell > c(\Pi)$, we need the stronger result that the image mod $\ell$ is abundant for all but finitely many $\ell$.   For this we use a recent theorem proved by one of us with Gajda and Petersen \cite{BGP}; see also \cite[E.10]{Dri16}.

\subsection{Local base change} From Theorem \ref{thm: intro BC} we draw some applications to the Local Langlands Correspondence. To state them, let $F$ be a local function field of characteristic $p$, $W_F$ the Weil group of $F$, and $G$ a split semisimple group over $F$. Genestier-Lafforgue \cite{Gen17} have constructed a semi-simplified local Langlands correspondence
\begin{align}\label{eq: LLC}
\left\{ \begin{array}{@{}c@{}}  \text{irreducible admissible representations} \\  \text{$\pi$ of $G(F)$ over $\ol{\Q}_{\ell}$}\end{array} \right\}/\sim  & \longrightarrow \left\{ \begin{array}{@{}c@{}}  \text{semi-simple $L$-parameters} \\ 
\sigma_{\pi} \co W_F \rightarrow G(\ol{\Q}_{\ell})  \end{array} \right\}/\sim.
\end{align}
If $F'/F$ is a field extension, we say that an irreducible admissible representation $\pi'$ of $G(F')$ is a \emph{base change lifting} of an irreducible admissible representation $\pi$ of $G(F)$ if $\sigma_{\pi'} \cong \sigma_{\pi}|_{W_{F'}}$. 

Using our global results, we deduce in \S \ref{sec: local base change} the following Theorem.

\begin{theorem}
Let $\pi$ be an irreducible admissible representation of $G(F)$ over $\ol{\Q}$. There exists a constant $c(\pi)$ such that for all primes $\ell > c(\pi)$, for any $\Z/\ell \Z$-extension $F'/F$ there exists a base change lifting of $\pi \otimes_{\ol{\Q}} \ol{\Q}_{\ell}$ to $G(F')$.
\end{theorem}

An analogous theorem was proved for mod $\ell$ representations, when $F'/F$ is cyclic of degree $\ell$, in \cite{F20}. The strategy here is completely different. The idea is to globalize $\pi$ to a cuspidal automorphic representation $\Pi$, and the extension $F'/F$ to a global extension $K'/K$, to which we can apply Theorem \ref{thm: intro BC}. If we can do this, then we can extract the local component of $\Pi'$ to obtain a local base change $\pi'$. However, we must take care to construct a globalization $\Pi$ satisfying the hypotheses of Theorem \ref{thm: intro BC}. This is accomplished by setting up $\Pi$ with specified local conditions at a finite number of auxiliary places, so that the corresponding Galois representation $\sigma_{\Pi}$ has Zariski-dense image. For example, we put a supercuspidal local component ``$V_{\phi}$'' at one auxiliary place whose Genestier-Lafforgue parameter $\sigma_{V_{\phi}}$ is already absolutely irreducible, in order to guarantee that $\sigma_{\Pi}$ is absolutely irreducible. For this, we require knowledge of $\sigma_{V_{\phi}}$. This is obtained by finding a \emph{different} globalization $\Pi^{\mrm{aux}}$ of $V_{\phi}$ whose corresponding global Galois representation $\sigma_{\Pi^{\mrm{aux}}}$ can be calculated explicitly, from which we extract the Genestier-Lafforgue paramater of $V_{\phi}$ using local-global compatibility; the work of Heinloth-Ng\^{o}-Yun \cite{HNY} provides a convenient such ``auxiliary globalization''. 

In the final section \S \ref{sec: toral supercuspidal}, we study local base change more explicitly for a class of supercuspidal representations singled out in \cite{Kal19} under the name ``toral supercuspidal representations''. This is a fairly broad class of supercuspidal representations, encompassing examples of arbitrary depth, for which Kaletha has constructed in \cite[\S 6]{Kal19} an explicit parametrization by $L$-parameters, satisfying good properties of the expected Local Langlands correspondence. For such representations with unramified input data, we are able to explicitly characterize unramified base change along the Genestier-Lafforgue correspondence under some technical assumptions. A motivation behind the results of \S \ref{sec: toral supercuspidal} is consistency between Kaletha's parametrization and the Genestier-Lafforgue correspondence. Roughly speaking, we prove, under genericity conditions, that the base change of mod $\ell$ toral supercuspidal representations along unramified $\Z/\ell \Z$-extensions, under the Genestier-Lafforgue correspondence, behaves as would be predicted by Kaletha's parametrization. See \S \ref{ssec: toral supercuspidal base change} for the precise formulations. 

Next we comment on the proof. The mechanism for understanding base change of mod $\ell$ representations is a Conjecture of Treumann-Venkatesh \cite[Conjecture 6.6]{TV16}, which predicts that Tate cohomology should realize functoriality in the Local Langlands correspondence for mod $\ell$ representations. This was proved in \cite[Theorem 1.2]{F20} for the Genestier-Lafforgue correspondence, thanks to which our task amounts to computing the Tate cohomology of toral supercuspidal representations. However, Tate cohomology is tricky to calculate in general; for a general supercuspidal representation presented via a Yu datum, as considered in \cite{Kal19}, it would be a challenge even to determine whether the Tate cohomology is non-zero. Crucial traction for this problem is provided by recent work of Chan-Ivanov and Chan-Oi \cite{CI19, CO21} on geometric models for toral supercuspidal representations in terms of ``deep level Deligne-Lusztig induction''; we ultimately compute the Tate cohomology by applying equivariant localization tools to their deep level Deligne-Lusztig varieties.

\subsection{Acknowledgments}

We thank Anna Cadoret, Charlotte Chan, H\'{e}l\`{e}ne Esnault, Kazuhiro Ito, David Hansen, Mikko Korhonen, Marie-France Vign\'eras, and Weizhe Zheng for helpful conversations. 

G.B.  was supported by the DFG grant CRC/TRR 326, project number 444845124. T.F. was partially supported by by an NSF Postdoctoral Fellowship under grant No. 1902927, as well as the Friends of the Institute for Advanced Study.    M.H. was partially supported by NSF Grant DMS-2001369 and by a Simons Foundation Fellowship, Award number 663678. C.K.  was partially supported by NSF Grant DMS-2200390. J.T.’s work received funding from the European Research Council (ERC) under the European Union’s Horizon 2020 research and innovation programme (grant agreement No 714405). 

\section{Notation and terminology}\label{sec_notation}

We fix a finite field $\bbF_q$ of characteristic $p$. Let $X$ be a smooth, geometrically connected, and projective curve over $\bbF_q$ and $K = \F_q(X)$, and let $G$ be a  split reductive group over $K$. The notation $\ell$ always denotes a prime not equal to $p$.

\subsection{Notation related to global and local fields} We write $K^s$ for a fixed choice of separable closure and $\Gamma_K := \Gal(K^s / K)$ for the corresponding Galois group. For $S$ a finite set of places of $K$, we write $K_S$ for the maximal subextension of $K^s$ unramified outside $S$, and $\Gamma_{K, S} := \Gal(K_S / K)$. If $v$ is a place of $K$, then $\Gamma_{K_v} = \Gal(K_v^s / K_v)$ will denote the decomposition group, and $\Gamma_{K_v} \to \Gamma_K$ the homomorphism corresponding to a fixed choice of $K$-embedding $K^s \hookrightarrow K_v^s$. If $v \not\in S$, then $\Frob_v \in \Gamma_{K, S}$ denotes a choice of geometric Frobenius element at the place $v$. We will identify the set of places of $K$ with the set of closed points of $X$. For a place $v \in |X|$ we write $q_v = \# k(v) = \# (\cO_{K_v} / \varpi_v \cO_{K_v})$ for the size of the residue field at $v$. We write $| \cdot |_v$ for the norm on $K_v$, normalized so that $|\varpi_v|_v = q_v^{-1}$; then the product formula holds. We write $\widehat{\cO}_{K} = \prod_{v \in X} \cO_{K_v}$. We will write $W_{K_v}$ for the Weil group of the local field $K_v$. 


\subsection{Notation related to group schemes}
If $G, H, \dots$ are group schemes over a base $S$, then we use Gothic letters $\frg, \frh, \dots$ to denote their Lie algebras, and $G_T, \frg_T, \dots$ to denote the base changes of these objects relative to a scheme $T \to S$. If $G$ acts on an $S$-scheme $X$ and $x \in X(T)$, then we write $Z_G(x)$ or $Z_{G_T}(x)$ for the scheme-theoretic stabilizer of $x$; it is a group scheme over $T$. We denote the centre of $G$ by $Z_G$. We say that a group scheme $G$ over $S$ is \emph{reductive} if $G$ is smooth and affine with reductive (and therefore connected) geometric fibres.

\subsection{Notation related to coefficient rings} When doing deformation theory, we will generally fix a prime $\ell$ and an algebraic closure $\overline{\bbQ}_\ell$ of $\bbQ_\ell$. A finite extension $E / \bbQ_\ell$ inside $\overline{\bbQ}_\ell$ will be called a coefficient field; when such a field $E$ has been fixed, we will write $\cO$ or $\cO_E$ for its ring of integers, $k$ or $k_E$ for its residue field, and $\varpi$ or $\varpi_E$ for a choice of uniformizer of $\cO_E$. We write $\cC_\cO$ for the category of Artinian local $\cO$-algebras with residue field $k$; if $A \in \cC_\cO$, then we write $\ffrm_A$ for its maximal ideal. Then $A$ comes with the data of an isomorphism $k \cong A / \ffrm_A$. 

\subsection{Dual groups} In this paper, we will view the dual group $\wh{G}$ of $G$ as a split reductive group over $\Z$. Our definition of $\wh{G}$ follows \cite[\S 2.1]{BHKT}. A prime $\ell$ is called a \emph{very good characteristic} for $\hG$ if it satisfies the conditions in the table below for all the simple factors of $\wh{G}$ (referring to the absolute root system types): 
\begin{center}
\begin{tabular}{|l|l|}
\hline
Condition     & Types              \\ \hline
$\ell \nmid n+1$ & $A_n$              \\ \hline
$\ell \neq 2$    & $B, C, D, E, F, G$ \\ \hline
$\ell \neq 3$    & $E, F, G$          \\ \hline
$\ell \neq 5$    & $E_8$              \\ \hline
\end{tabular}
\end{center}

\section{Deformation problems}\label{sec: deformation}

In this section we set up the results on Galois deformation theory that will be used later for automorphy lifting. 

\subsection{Setup for the general theory}

Let $\hG$ be a split semisimple group over $\bbZ$. We fix a prime $\ell$ which is a very good characteristic for $\hG$, as well as a coefficient field $E \subset \overline{\bbQ}_\ell$. We also fix an absolutely irreducible representation $\overline{\sigma}$ of $\Gamma_{K,S}$. We recall the results of \cite{BHKT} on the deformation theory of representations of $\Gamma_{K,S}$ to $\hG$ with $\ell$-adic coefficients.

\begin{lemma}\cite[Lemma 5.1]{BHKT}\label{lem_triviality_of_scheme_theoretic_centralizer_in_very_good_characteristic}
Let $\bs:  \Gamma_{K,S} \ra \hG(k)$ be an absolutely $\hG$-irreducible homomorphism.
The scheme-theoretic centralizer of $\bs(\Gamma_{K,S})$ in $\hG_k^{\ad}$ is \'etale over $k$, and $H^0(\Gamma_{K,S}, \ad \bs) = 0$.
\end{lemma}
Let $A \in \cC_\cO$. We define {\it lifings} and {\it deformations} of $\bs$ over $A$ as in \cite[\S 5.1]{BHKT}, and let  
$\Def_{\bs} : \cC_\cO \to \Sets$ be the functor that associates to $A \in \cC_\cO$ the set of deformations of $\bs$ over $A$.  Then


\begin{proposition}\cite[Propositions 5.10]{BHKT}\label{prorep} The functor $\Def_{\bs}$ is pro-represented by a complete Noetherian local $\cO$-algebra $R_{\bs} = R_{\bs,S}$.
(We will include the subscript $S$ in the notation when we want to let $S$ vary.)
\end{proposition}




\subsection{Local deformation conditions}
So far we have considered deformations of $\bs$ with no restriction on the points in $S$.  We now introduce  local deformation conditions, following the discussion
in \cite[\S 2.2]{CHT08}, and for general groups in \cite[\S 3.2]{Pat16}, with notation as in \cite[\S 4]{Kha15}.   For each $v \in S$ let $\Gv \subset \G_{K,S}$ be a decomposition group at $v$, and let $\bs_v$ denote the restriction of $\bs$ to $\Gv$.  Let $\cD^\square_{\bs,v}$ denote the functor on $\cC_\cO$ which to $A \in \cC_\cO$ associates
the set of liftings of $\bs$ to $A$; then $\cD^\square_{\bs,v}$ is represented by a ring $R^\square_{\bs,v}$.  As in \cite[Definition 4.1]{Kha15}, we define a {\it local deformation problem for $\bs_v$} to be a subfunctor $\cD_v \subset \cD^\square_{\bs,v}$ such that
\begin{itemize}
\item $\cD_v$ is represented by a quotient of  $R^\square_{\bs,v}$; and
\item For any $A \in \cO$, if $\sigma \in \cD_v(A)$ and $\alpha \in \ker [\hG(A) \ra \hG(k)]$ then $\ad(\alpha)(\sigma) \in \cD_v(A)$.
\end{itemize}

\begin{definition}\label{globaldef}  A \emph{global deformation problem} for $\bs$ is a collection $\cD := (S, 
\{\cD_v\}_{v \in S})$ where for each $v \in S$, $\cD_v$ is a local deformation problem for $\bs_v$.

For $A \in \cC_\cO$, a deformation $\sigma:  \G_{K,S} \ra \hG(A)$ is of \emph{type $\cD$} if for every $v \in S$, the restriction of $\sigma$ to $\Gv$ belongs to $\cD_v$.

The functor of deformations of $\bs$ of type $\cD$ is denoted $\Def_{\bs,\cD} = \Def_{\bs,S,\cD}$. Again, we will emphasize the subscript $S$ when we want to contemplate varying $S$.
\end{definition}

To each local deformation problem $\cD_v$ we attach a $k$-subspace $L_v \subset H^1(G_{K_v}, \ad(\bs))$ as in \cite[\S 3.2]{Pat16}.  
\begin{definition}\label{minimal_deformation_problems}The deformation problem $\cD_v$ is {\em balanced} $\dim L_v = \dim H^0(G_{K_v}, \ad(\bs))$. 
\end{definition}

\begin{example}
The minimal deformation problems considered in \cite[\S 2.4.4]{CHT08} are balanced. 
\end{example}

\begin{examplen}\label{ex: h^2=0} The deformation problem $\cD_v$ is {\em unrestricted} if $\cD_v = \cD^\square_{\bs,v}$. In this case we have $L_v = H^1(\Gv, \ad(\bs))$.  

Recall that for all $v \in S$, Tate's local Euler characteristic formula for any $k[\Gv]$-module $M$ reduces to
$$h^0(\Gv,M) - h^1(\Gv,M) + h^2(\Gv,M) = 0$$
because $\ell \neq p$. Hence, if $H^2(\Gv, \ad(\bs)) = 0 $ then the unrestricted deformation problem $\cD_v$ is balanced. 

Note that local duality implies that
\begin{equation}\label{h0h2}
H^2(\Gv, \ad(\bs)) = 0 \Leftrightarrow H^0(\Gv, \ad(\bs)^\vee(1)) = 0.
\end{equation}
\end{examplen}

\subsection{Selmer groups}

Given a global deformation problem $\cD$, we define 
\begin{equation}\label{h1cD}
H^1_{\cD}(\Gamma_{K,S},\ad(\bs)) := \ker  \left[H^1(\Gamma_{K,S},\ad(\bs)) \ra \bigoplus_{v \in S} H^1(\Gv, \ad(\bs))/L_v \right],
\end{equation}
and let $h^1_{\cD}(\Gamma_{K,S},\ad(\bs) ) := \dim_k H^1_{\cD}(\Gamma_{K,S},\ad(\bs)).$

With $L_v$ as above, let $L_v^\perp \subset H^1(\Gv, \ad(\bs)^{\vee}(1))$ denote its orthogonal complement under the local duality pairing. We define 
\[
H^1_{\cD^{\perp}}(\Gamma_{K,S},\ad(\bs)^{\vee}(1)) := \ker  \left[H^1(\Gamma_{K,S},\ad(\bs)^{\vee}(1)) \ra \bigoplus_{v \in S} H^1(\Gv, \ad(\bs)^{\vee}(1))/L_v^{\perp} \right],
\]
and let $h^1_{\cD^{\perp}}(\Gamma_{K,S},\ad(\bs)^{\vee} (1)) := \dim_k H^1_{\cD^{\perp}}(\Gamma_{K,S},\ad(\bs)^{\vee}(1))$.

\begin{proposition}\label{repcD} (i)  The functor $\Def_{\bs,\cD}$ is pro-represented by a complete Noetherian local $\cO$-algebra $R_{\bs,\cD} = R_{\bs,S,\cD}$.

(ii) There is a surjection $ \cO \llbracket  X_1, \dots X_g \rrbracket  \surj R_{\bs, S, \cD }$, where $g = \dim_k H^1_{\cD}(\Gamma_{K, S}, \ad \bs)$.
\end{proposition} 
\begin{proof} For (i) the proof is the same as for \cite[Proposition 2.2.9]{CHT08}, and for (ii) the proof is the same as for \cite[Corollary 2.2.12]{CHT08}.
\end{proof}

\subsection{Abundant subgroups}

We recall two definitions from \cite{BHKT}.

\begin{definition}\label{conditions}  Let $\Gamma$ be an abstract group, and let $\hG$ be a reductive group over $k$.
\begin{enumerate} \item A homomorphism $\sigma : \Gamma \to \hG(k)$ is said to be \emph{absolutely $\hG$-completely reducible} (resp. \emph{absolutely $\hG$-irreducible}) if it is $\hG$-completely reducible (resp. $\hG$-irreducible) after extension of scalars to an algebraic closure of $k$.
\end{enumerate} 
Suppose now that $\hG$ is semisimple, and let $\overline k$ be an algebraic closure of~$k$.
\begin{enumerate} \item[(iii)] 
A finite subgroup  $H \subset \hG(\overline k)$ is said to be \emph{$\hG$-abundant} over $k$, if it satisfies the following conditions:
\begin{enumerate}
\item[(a)] The group $H$ is contained in $\hG(k)$.
\item[(b)] The groups $H^0(H, \widehat{\frg}_{k})$, $H^0(H, \widehat{\frg}^\vee_{k})$, $H^1(H, \widehat{\frg}^\vee_{k})$ and $H^1(H, k)$ all vanish. 
\item[(c)]  For each regular semisimple element $h \in H$, the torus $Z_{\hG}(h)^\circ \subset \hG$ is split over~$k$.
\item[(d)] For every simple $k[H]$-submodule $W \subset \widehat{\frg}^\vee$, there exists a regular semisimple element $h \in H$ such that $W^h \neq 0$ and $Z_{\hG}(h)$ is connected. (We recall that $Z_{\hG}(h)$ is always connected if $\hG$ is simply~connected.)
\end{enumerate}
\end{enumerate}
\end{definition}

\begin{lemma}\label{abundantimage}  Suppose $\hG$ is simple and split over $k$. Let $\hG^{sc}$ denote the simply-connected cover of $\hG$, and let $\hG(k)^+$ denote the image of $\hG^{sc}(k)$ in $\hG(k)$. Let $H$ be a group with  $\hG^+(k)\subset H\subset \hG(k)$ and suppose that  $(k,\mathrm{type\, of\, }\hG)$ is not in the following list:
\[
\{(\mathbb F_3,A_1), (\mathbb F_5,A_1)
 \}\cup \{(\mathbb F_q,C_n)\mid q\in\{3,5,9\},n\ge2\}
\]
Then there exists a finite field $k_1\supset k$ such that for all finite fields $k'\supset  k_1$ the group $H$ is $\hG$-abundant over~$k'$.
\end{lemma}

\begin{proof}      
We will prove the following claims for $(k,\mathrm{type\, of\, }\hG)$ not in the above list:
\begin{enumerate}
\item[(i)]  The adjoint module $\widehat{\frg}_{k}$ is self-dual as a $k[H]$-representation,
\item[(ii)] The groups $H^0(H, \widehat{\frg}_{k})$, $H^1(H, \widehat{\frg}_{k})$ and $H^1(H, k)$ all vanish.
\item[(iii)]  The action of $H$ on $\widehat{\frg}_{k}$ is absolutely irreducible.
\item[(iv)]  The group $H$ contains strongly regular semisimple elements~$h$, i.e., $h$ such that $Z_{\hG}(h)$ is connected. 
\end{enumerate}
Let us first show that the claims suffices to prove the lemma.  Indeed, condition (a) is obvious, condition (b) is the same as (ii) (given (i)), while condition (c) is automatic if $k_1$ is chosen to split all of the finitely many $k$-rational tori $Z_{\hG}(h)^\circ$ in $\hG$. Finally, given (iii), condition (d) comes down to the existence of strongly regular semisimple $h$, as in (iv), because $(\widehat{\frg}^\vee)^h$ always contains the Lie algebra of the centralizer of $h$.

It remains to establish the claim. Parts (i) to (iii) with many details can be found in \cite[\S~3]{AB20}, which however largely builds on further references, which we now briefly recall. We first observe that $\widehat{\frg}_{k}$ is an absolutely irreducible $k[H]$-representation, essentially by results of Hiss \cite{Hiss} and Hogeweij \cite{Hogeweij}; see \cite[Proposition~3.22]{AB20}. Now because $\widehat{\frg}_{k}$ is absolutely irreducible, to prove self-duality, it suffices to see that $\widehat{\frg}_{k}$ and $\widehat{\frg}_{k}^\vee$ have the same weights. But this is clear since with any root $\alpha$ of the Lie algebra $\widehat{\frg}_{k}$ also $-\alpha$ is a root, and both have multiplicity one. This settles (i) and~(iii), and the first assertion of~(ii). 

Next, a largely classical result, completed by Tits, asserts that  $\hG(k)^+$ is perfect, unless $(k,\mathrm{type\, of\, }\hG)$ is $(\mathbb F_3,A_1)$, see \cite[Thm.~24.17]{MalleTesterman}; because $\hG$ is split and $l$ is of very good characteristic for $\hG$, there is only a single exception. The latter also implies that $\hG(k)/\hG(k)^+$ is of order prime to the characteristic  $l$ of $k$. Hence $H^1(H, k)=\Hom(H, k)$ vanishes, and this shows the third part of (ii); see \cite[Corollary~3.12]{AB20}.

We now turn to the remaining condition $H^1(\hG(k), \widehat{\frg}_{k}) = 0$ from (ii). The most complete vanishing results for $H^1(H, \widehat{\frg}_{k})$ are due to V\"olklein with earlier work by Cline, Parshall, Scott in the A-D-E cases and by Hertzig (unpublished). We observe that  $\Lie\hG$ has trivial center because $l$ is very good for $\hG$. The latter also excludes the case  $(k,\mathrm{type\, of\, }\hG)= (\mathbb F_2,A_1)$. The cases $ (\mathbb F_5,A_1)$ and $(\mathbb F_q,C_n)$, $q\in\{3,5,9\}$, $n\ge2$, for  $(k,\mathrm{type\, of\, }\hG)$ are ruled out by our list of exceptions. It follows from \cite[Theorem and Remarks]{Vol89} that $H^1(\hG(k), \widehat{\frg}_{k}) = 0$.

\smallskip

We finally prove (iv). If $\hG$ is simply connected, then by \cite[II.3.9]{Spr70} the group $Z_{\hG}(h)$ is connected for all regular semisimple $h\in H$, and we are done. We now focus on the case where $\hG$ of adjoint type, where we find a suitable $h$ using \cite{FS88}. The last paragraph explains how this gives (iv) for $\hG$ of any other type.

Suppose that $\hG$ is of adjoint type. Then by \cite[Theorem~3.1]{FS88} there exists a maximal torus $T\subset \hG$ defined over $k$ such that the following holds: 
\begin{enumerate}
\item[(1)] The group $T^+(k):=\hG^+(k)\cap T(k)$ is cyclic. Let $h$ be a generator.
\item[(2)] One has $Z_{\hG(k)}(h)=T(k)$, and the order of $h$ is given in  \cite[Table III]{FS88}.
\item[(3)] If $\hG$ is of type $A_n$, then $T$ is quotient of a totally anisotropic maximal torus of $\GL_{n+1}$ (over $k$) modulo its center, i.e. $T\cong \Res_{k^{(n+1)}}^k \mathbb G_m/\mathbb G_m$, 
the quotient of the Weil restriction of $\mathbb G_m$ from the unique degree $n+1$-extension $k^{(n+1)}$ of $k$ to $k$ modulo $\mathbb G_m$.
\end{enumerate}

Let $C:=Z_{\hG}(h)$, and let $\pi_0(C)=C/C^0$ be its component group scheme. It follows from \cite[2.9 and its proof]{Sei83} that the identity component $C^0$ is a (commuting) product of Chevalley groups over $k$ and of the center of $C^0$, and so from (1) above the natural map $T\to C^0$ must be an isomorphism because any Chevalley group has $k$-points of order $l$. 
 In particular, $h$ is regular semisimple.

Next we gather results on $\pi_0( C)$ that will allow us to identify it with the trivial group.
By \cite[Corollary~II.4.4]{Spr70} there is an injective morphism of group schemes $\pi_0(C)\to Z(\hG^{sc})$ to the center $Z(\hG^{sc})$ of $\hG^{sc}$. Because $l$ is good for $\hG$, the finite commutative group scheme $Z(\hG^{sc})$ is \'etale over $k$. The centers are completely known, and a list can be found for instance in \cite[Table 1]{AB20}. Using the Bruhat-decomposition of $\hG$ over $\overline k$, and that $h$ is regular semisimple, it is not difficult to show that $C$ is a subgroup of the normalizer of $T$ in $\hG$ (over $\overline k$), and hence $\pi_0(C)$ is a subgroup of the Weyl group of $G$. Moreover by \cite[I.2.11]{Spr70} one has a short exact sequence 
\[ 0\to T(k) \to C(k) \to   \pi_0(C)(k)\to 0.\]
From (2) above we deduce that $T(k)\to C(k)$ is an isomorphism, and hence $\pi_0( C)(k)$ is trivial.

Now suppose the type is $B_n$, $C_n$, $D_n$ or $E_7$. Then $l$ cannot be $2$, and so $2$ divides the order of $k^\times$. It follows from the classification of the centers $Z(\hG^{sc})$, that any non-trivial subgroup scheme will contain an element of exact order $2$. But then $\pi_0( C)$ must be trivial, because we know already that $\pi_0( C) (k)$ is trivial. 

The types $E_8$, $F_4$ and $G_2$ need not be considered, since here $Z(\hG^{sc})$ is trivial and hence $\hG$ itself is simply connected. Next we consider type $E_6$, so that $Z(\hG^{sc})\cong \mu_3$. Let $q=\#k$. If $q \equiv 1 \pmod 3$, we can argue as in the previous paragraph to deduce that $\pi_0( C)$ is trivial. If on the other hand $q \equiv -1 \pmod 3$, then from \cite[Table III]{FS88} we find that
\[ \#T(k) = (q^2+q+1)(q^4-q^2+1)\equiv -1\pmod 3,\]
so that the order of $h$ is not divisible by $3$. It follows from  \cite[Corollary~4.6]{Spr70} that $\pi_0( C)$ is connected.

It remains to discuss the type $A_n$. For $n=1$, by (3) the element $h$ is given as a diagonal matrix $\diag(x,x^q)$ in $\GL_2/\mathbb G_m$ for some $x\in (k^{(2)} )^\times$ of order $q^2-1$, so that $h$ has order $q+1$. Because $C$ is spanned by $T$ and at most some Weyl group elements, we need to understand whether the unique non-trivial Weyl group element $w$ of $\hG=\PGL_2$ lies in $C$. The element $w$ exchanges the entries of $h$. If it would fix $h$, then $\diag(x,x^q)\equiv \diag(x^q,x)\pmod{\overline k^\times}$, so that there exists $a\in \overline k^\times$ with $ax=x^q$ and $ax^q=x$. This firstly implies $a\in\{\pm1\}$, and then that the order of $x$ divides $2(q-1)$. Because $q+1$ divides the order of $x$, this could only happen  for $q=2$, which is forbidden, since $l$ is assumed to be very good for $\hG=\PGL_2$.

Now suppose $n\ge2$. Then $h$ is given by $\diag(x,x^q,\ldots,x^{q^n})\in\GL_{n+1}/\mathbb G_m(\overline k)$ for some $x\in (k^{(n+1)})^\times$ of order $q^{n+1}-1$. Let $w$ be in the Weyl group of $\PGL_n$, which we identify with $S_{n+1}$. Suppose that $h=w\circ h$, i.e., that $w$ is the image of an element of $C$. Because $\pi_0( C)\subset Z(\PGL_{n+1})\cong \mu_{n+1}$ we may assume that $w$ is cyclic of order dividing $n+1$. Let $r_1\ge\ldots r_s\ge2$ be the lengths in the cycle decomposition of $w$, and let $r=\gcd(r_1,\ldots r_s)$, which is now a divisor of $n+1$. As in the case $n=1$, we find $a\in \overline k^\times$ of order dividing $r$ and such that $x^{q^{w(i)}}=ax^{q^i}$, or equivalently $x^{q^{w(i)}-q^i}=a$, for $i=0,\ldots,n$.

To get further, let $t$ be a large Zsygmondy prime for $(q,n+1)$ in the sense of  \cite[Theorem~2.1]{FS88}, i.e., $t$ is a prime that divides $q^{n+1}-1$, that does not divide $q^m-1$ for $1\le m \le n$ and with either $t>n+1$ or $t^2|(q^{n+1}-1)$. We deduce that $t$ divides the order of $x^{q^{w(i)}-q^i}$ for any $i$ with $w(i)\neq i$, and it follows that $t$ divides $r$ and in turn $n+1$, unless $w$ is trivial. But then modulo the prime $t$ and using Fermat's Little Theorem we have
\[ q^{n+1}-1 = (q^t)^{(n+1)/t} -1\equiv q^{(n+1)/t} -1 \pmod t\]
This contradicts that $t$ is a Zsygmondy prime for $(q,n+1)$. Hence $w$ must be trivial and $T=C$.

Finally, for a general $\hG$ consider the central isogeny $\hG\to \hG^{ad}$ to the adjoint quotient $\hG^{ad}$ of $\hG$. Let $H^{ad}\subset \hG^{ad}(k)$ be the image of $H$. Now the element $h$ in $H^{ad,+}(k)$ from (1) above clearly lifts to some $h'\in \hG^{+}(k)$ with $h'$ regular semisimple. Moreover we have a left exact sequence $0\to Z(\hG)\to Z_{\hG}(h')\to Z_{\hG^{ad}}(h)$. In it, $Z_{\hG^{ad}}(h)$ is a maximal torus, $Z_{\hG}(h')$ contains a maximal torus, and the center $Z(\hG)$ lies in any maximal torus of $\hG$. Hence $Z_{\hG}(h')$ is a (maximal) torus, and thus connected.
\end{proof}

\subsection{Taylor-Wiles primes}\label{ssec: T-W}

Folllowing \cite[Definition 5.16]{BHKT}, we define a \emph{Taylor-Wiles datum for $\overline{\sigma}$} to be a pair $(Q, \{ \varphi_v \}_{v \in Q})$, where:
\begin{itemize}
\item $Q$ is a finite set of places $v$ of $K$ such that $\overline{\sigma}(\Frob_v)$ is regular semisimple and $q_v \equiv 1 \text{ mod }\ell $. 
\item For each $v \in Q$, $\varphi_v : \hT_k \cong Z_\hG(\overline{\sigma}(\Frob_v))$ is a choice of inner isomorphism. In particular, this forces $Z_\hG(\overline{\sigma}(\Frob_v))$ to be connected. 
\end{itemize}

For a global deformation problem $\cD = (S, \{\cD_v\}_{v \in S})$ and a Taylor-Wiles datum $(Q, \{ \varphi_v \}_{v \in Q})$ with $Q \cap S = \emptyset$, we define a new global deformation problem $\cD_Q = (S \cup Q, (\cD_Q)_v )$ with
\[
(\cD_Q)_v = \begin{cases} \cD_v & v \in S, \\ \cD_{\bs, v}^{\square} & v \in Q. \end{cases}
\]

The Proposition below guarantees the existence of sets of Taylor-Wiles primes suitable for patching, generalizing \cite[Proposition 5.19]{BHKT}. Actually there are some typos in the statement and proof of \emph{loc. cit.}; the argument below can be used to rectify them. 

\begin{proposition}\label{minimal_taylor_wiles_data}
Suppose that the group $\overline{\rho}(\Gamma_{K(\zeta_\ell)}) \subset \hG(k)$ is $\hG$-abundant. Let $\cD := (S, \{\cD_v\}_{v \in S})$ be a global deformation problem such that for each $v \in S$ the local deformation problem $\cD_v$is balanced, with corresponding subspace $L_v \subset H^1(\Gv, \ad(\bs))$.  Then for each $j \geq 1$, there is a Taylor--Wiles datum $(Q, \{ \varphi_v \}_{v \in Q})$, satisfying the following conditions:
\begin{enumerate}
\item $Q \cap S = \emptyset$, and for each $v \in Q$, we have $q_v \equiv 1 \text{ mod }\ell^j$, and 
\[
\# Q = h^1_\cD(\Gamma_{K, S}, \ad(\bs)) = h^1_{\cD^{\perp}}(\Gamma_{K, S}, \ad(\bs)^\vee(1)).
\]
\item We have $h^1_{\cD_Q}(\Gamma_{K, S \cup Q}, \ad(\bs)^\vee(1))  = 0$. 
\item There exists a surjection $\cO \llbracket X_1, \dots, X_g \rrbracket \surj R_{\bs, S \cup Q, \cD}$ with $g = h^1_{\cD}(\Gamma_{K, S}, \ad(\bs)) + (r-1)\# Q$, where $r = \rank \hG$.
\end{enumerate}
\end{proposition}

\begin{proof} 
If $S \cup Q = \emptyset$ then (i) and (ii) are vacuous and (iii) is an immediate consequence of Proposition \ref{repcD} (ii). For the rest of the proof, we assume that $S \cup Q$ is non-empty. 

We apply \cite[Theorem 6.2]{Ces15}, using that $H^0(\Gamma_{K, S \cup Q}, \ad (\bs)(1)) = H^0(\Gamma_{K, S \cup Q}, \ad (\bs)) = 0$ by the abundance assumption, and the local and global Euler characteristic formulas (note our assumption that $S \cup Q$ is non-empty) to deduce that
\begin{equation}\label{eq: dimension equation}
h^1_{\cD_Q}(\Gamma_{K, S \cup Q} ,\ad(\bs)) - h^1_{\cD_Q^{\perp}}(\Gamma_{K, S \cup Q} ,  \ad (\bs)^{\vee}(1))  =  \sum_{v \in S}  \left(\dim \cL_v - h^0(\Gamma_v, \ad (\bs)) \right)  + \sum_{v \in Q} r .
\end{equation}
The assumption that $\cD_v$ is balanced for all $v \in S$ implies that $\dim \cL_v - h^0(G_{F_v}, \ad (\bs)) = 0$ for all $v \in S$. We claim that once we can arrange (i) and (ii), then (iii) follows automatically. Indeed, if $h^1_{\cD_Q^{\perp}}(\Gamma_{K, S \cup Q} , \ad (\bs)^{\vee}(1))$ vanishes then \eqref{eq: dimension equation} implies that 
\[
h^1_{\cD_Q}(\Gamma_{K, S \cup Q} , \ad (\bs))  = r \# Q = h^1_{\cD}(\Gamma_{K, S}, \ad(\bs)) + (r-1)\# Q,
\]
from which (iii) follows by invoking Proposition \ref{repcD}.

It therefore remains to show that $Q$ can be chosen so that 
\[
H^1_{\cD_Q^{\perp}}(\Gamma_{K, S \cup Q} , \ad(\bs)^{\vee}(1)) := \ker \left(  H^1(\Gamma_{K, S \cup Q} , \ad(\bs)^{\vee}(1)) \rightarrow \left( \bigoplus_{v \in S} \frac{H^1(K_v, \ad(\bs)^{\vee}(1)) }{ L_v}  \oplus \bigoplus_{v \in Q} H^1(K_v, \ad(\bs)^{\vee}(1))  \right) \right) 
\]
vanishes. By a comparison of inflation-restriction exact sequences, the inflation map $H^1(\Gamma_{K, S }, \ad(\bs)^{\vee}(1)) \rightarrow H^1(\Gamma_{K, S \cup Q}, \ad(\bs)^{\vee}(1))$ takes 
\begin{align*}
&\ker \left(  H^1(\Gamma_{K, S} , \ad(\bs)^{\vee}(1)) \rightarrow  \bigoplus_{v \in Q} H^1(k(v), \ad(\bs)^{\vee}(1))  \right)  \\ 
&\hspace{2cm} \xrightarrow{\sim}\ker \left(  H^1(\Gamma_{K, S \cup Q} , \ad(\bs)^{\vee}(1)) \rightarrow  \bigoplus_{v \in Q} H^1(K_v, \ad(\bs)^{\vee}(1))  \right) .
\end{align*}
Therefore it suffices by induction to show that for any $j$ and any non-zero $[\psi] \in H^1(\Gamma_{K, S}, \wh{\mf{g}}_k^{\vee}(1)) )$, we can find infinitely many places $v \not\in S$ of $K$ such that $q_v \equiv 1 \text{ mod }\ell^j$, $\overline{\rho}(\Frob_v)$ is regular semisimple with connected centralizer in $\hG_k$, and $\res_{K_v} [\psi] \neq 0$. The rest of the argument concludes exactly as in \cite[Proof of Proposition 5.19]{BHKT}.

\end{proof}

\section{Automorphy lifting}\label{sec: automorphy lifting}


In this section we prove the automorphy lifting theorem that will be used to deduce cyclic base change. Such automorphy lifting theorems have been the subject of much study over number fields, and our proof employs similar techniques, but we are able to obtain results for more general groups thanks to the different numerical behavior of the Euler characteristics of function fields.  

To put our results in context, we compare them to \cite[Theorem 8.20]{BHKT}, which handled everywhere unramified Galois representations. The main novelty of the present situation is that we need to implement automorphy lifting for arbitrarily ramified residual representations. In this section we use the Taylor-Wiles method to prove an automorphy lifting theorem (Theorem \ref{thm_R_equals_B}) whenever one knows that all the local deformation problems are well-behaved, e.g., formally smooth. In the unramified case it was easy to see that the unramified local deformation problem is formally smooth; in the next section we give a suitable extension to unrestricted deformations for ramified representations, when $\ell$ is large enough.

\subsection{Spaces of automorphic forms}\label{ssec: automorphic forms}
We define integral spaces of automorphic forms as in \cite[\S 8.1]{BHKT}. For any open subgroup $U \subset G(\wh{\cO})$ and any $\Z[1/p]$-algebra $R$, we define:
\begin{itemize}
\item $C(U, R)$ to be the $R$-module of functions $f : G(K) \backslash G(\bbA_K) / U \to R$;
\item $C_c(U, R) \subset C(U, R)$ to be the $R$-submodule of functions $f$ which have finite support;
\item and $C_{\text{cusp}}(U, R)$ to be the $R$-submodule of functions $f$ which are cuspidal, in the sense that for all proper parabolic subgroups $P \subset G$ and for all $g \in G(\bbA_K)$, the integral
\[ \int_{n \in N(K) \backslash N(\bbA_K)} f(ng)\, dn \]
vanishes, where $N$ is the unipotent radical of $P$.  
\end{itemize}
This last integral is normalized by endowing $N(K) \backslash N(\bbA_K)$ with its probability Haar measure (which makes sense because we are assuming that $p$ is a unit in $R$). We have $C_{\cusp}(U,R) \subset C_{c}(U, R)$ by\cite[Proposition 8.2]{BHKT}. 

We define $C_{\cusp}(G,R) := \varinjlim_U C_{\text{cusp}}(U, R)$. 

Let  $N = \sum_v n_v \cdot v \subset X$ be an effective divisor and let $U(N) := \ker( \prod_v G(\cO_{K_v}) \to G(\cO_N) )$. The underlying set of places $|N| = \bigcup_{n_v>0} \{v\}$ will play the role of the set $S$ from \S \ref{sec: deformation}. 


\subsection{The excursion algebra and Lafforgue's parametrization}\label{sec: excursion algebra} For the remainder of this section we choose a prime $\ell \neq p$ and coefficient field $E \subset \overline{\bbQ}_\ell$. For $R \in \{k, \cO,E\}$ we denote by $\cB(U(N), R)$ the $R$-subalgebra of $\End_{R}(C_{\text{cusp}}(U(N), R))$ generated by V. Lafforgue's {\it excursion operators}, as in \cite[\S 8.4]{BHKT}. For any fixed $N$, this is a finite $R$-algebra. The Hecke operators $T_{V,v}$ as lying in $\cB(U(N), R)$, for
$v \notin |N|$ and $V$ a representation of $\hG$, are defined as in \cite{BHKT}.

The points of $\Spec \cB(U(N),R)$ are naturally identified with semi-simple $L$-parameters, in a manner that we will presently review. Combining \cite[Corollary 8.6, Corollary 8.11]{BHKT}, we have the following (with notation as in \emph{loc. cit.}):

\begin{theorem}\label{BHKT}
Let $R = E$ or $k$ and $\mf{m} \subset \cB(U(N),R)$ be a maximal ideal. Fix an embedding $\cB(U(N),R)/\mf{m} \inj \ol{R}$. Then:
\begin{enumerate} \item There exists a continuous, absolutely $\hG$-completely reducible representation $\sigma_{\mf{m}} : \Gamma_K \to \hG( \overline{R} )$ satisfying the following condition: for all excursion operators $S_{I, (\gamma_i)_{i \in I}, f}$, we have
\begin{equation}\label{eqn_compatibility_with_excursion_operators} f( (\sigma_{\mf{m}}(\gamma_i))_{i \in I}) = S_{I, (\gamma_i)_{i \in I}, f} \text{ mod }{\mf{m}}. 
\end{equation}
\item The representation $\sigma_{\mf{m}}$ is uniquely determined up to $\hG(\ol{R})$-conjugacy by (\ref{eqn_compatibility_with_excursion_operators}).
\item The representation $\sigma_{\mf{m}}$ is unramified outside $|N|$. If $v \not\in |N|$, then it satisfies the expected local-global compatibility relation at $v$: for all irreducible representations $V$ of $\hG_{\overline{\bbQ}_\ell}$, we have $T_{V, v} \in \cB(U(N), R)$ and
\[ \chi_V( \sigma_{\mf{m}}(\Frob_v) ) = T_{V, v}  \in \cB(U(N), R)/\mf{m}. \]
\end{enumerate}


\end{theorem}

For $R \in \{k, E \}$ and any $N$, V. Lafforgue constructs a decomposition of $C_{\cusp}(U(N), \ol{R})$ into summands indexed by $L$-parameters $\sigma \in H^1(\Gamma_{K,N} , \wh{G}(\ol{R}))$. More specifically, this decomposition comes from the generalized eigenspace decomposition for the action of $\cB(U(N), R)$ on $C_{\cusp}(U(N), R)$.  If $\Pi \subset C_{\cusp}(G, \ol{R})$ is a subspace such that $\Pi^{U(N)}$ is stable under the $\cB(U(N), \ol{R})$-action and supported over a unique maximal ideal, we will denote this maximal ideal by $\ffrm = \ffrm_\Pi \subset \cB(U(N), \ol{R})$. The corresponding $L$-parameter is then denoted $\sigma_{\Pi} := \sigma_{\ffrm}$, and we say that $\sigma_{\Pi} $ is the $L$-parameter \emph{attached} to $\Pi$. This assignment $\Pi \mapsto \sigma_{\Pi}$ is independent of the choice of $N$ such that $\Pi^{U(N)} \neq 0$. This correspondence has the property that $\Pi$ and $\sigma_{\Pi}$ match under the Local Langlands correspondence at all places where $\Pi$ is unramified, although we caution that this property does not characterize it for general groups.

\begin{remark}[A one-to-finite correspondence]\label{rem: weak correspondence}
We continue to assume $R \in \{k, E\}$. Not every irreducible $G(\A)$-subspace $\Pi$ of $C_{\cusp}(G,R)$ is attached to an $L$-parameter $\sigma_{\Pi}$ in the sense of the previous paragraph. The notation $\sigma_{\Pi}$ will only be used if $\Pi$ is attached to $\sigma_{\Pi}$ in the manner of the above paragraph. However, we do have the following construction which produces a finite set of $L$-parameters, each of which matches $\Pi$ locally at all unramified places. 

Let $C_{\cusp}(G,R)[\Pi]$ be the $\Pi$-isotypic subspace of $C_{\cusp}(G,R)$ for the $G(\A_K)$-action. Then $C_{\cusp}(G,R)[\Pi]^{U(N)}$ is stable under the action of $\cB(U(N), R)$, and finite-dimensional for every $N$. We define $\ul{\rho}_{\Pi}$ to be the \emph{finite set} of $L$-parameters corresponding to the maximal ideals of $\cB(U(N), R)$ in the support of $C_{\cusp}(G,R)[\Pi]^{U(N)}$. Since each $\sigma \in \ul{\rho}_{\Pi}$ is associated to an automorphic representation isomorphic to $\Pi$ at all local places, in particular the unramified ones, we have that each $\rho$ matches $\Pi$ at every unramified place under the (unramified) Local Langlands correspondence.

We note that by \cite[Proposition 6.4]{BHKT}, if some $\sigma \in \ul{\rho}_{\Pi}$ has Zariski-dense image in $\wh{G}$, then in fact $\ul{\rho}_{\Pi}$ is a singleton, so $\ul{\rho}_{\Pi} = \{ \sigma_{\Pi}\}$.

In \S \ref{mwf} we will need to work with collections of representations with $\ol{\Q}_{\ell}$ or $\overline{\bbF}_\ell$ coefficients, as $\ell$ varies.  We may therefore 
write $\sigma_{\Pi, \ell}$ for $\sigma_{\Pi}$ and $\bs_{\Pi,\ell}$ for $\bs_\Pi$ when it is necessary to specify the coefficients.  
\end{remark}

\begin{definition}\label{MTgroup}  Let $\Pi$ be a cuspidal automorphic representation of $G(\bbA_K)$ defined over a number field $L$.  We say the \emph{Mumford-Tate group of $\Pi$ is $\hG$} if for some (equivalently, every) prime $\ell \neq p$ and embedding $L \to \overline{\bbQ}_\ell$, the image of $\sigma_{\Pi, \ell}$ is Zariski dense in $\hG$. As explained in Remark \ref{rem: weak correspondence}, this notion does not depend on the choice of $\sigma_{\Pi, \ell} \in  \ul{\rho}_{\Pi \otimes_L \overline{\bbQ}_\ell}$.
\end{definition}

\subsection{Automorphy lifting theorems}

Let $\chi \co \cB(U(N), \cO)  \rightarrow \overline{\Z}_{\ell}$ be a homomorphism. After possibly enlarging $E$, we can assume that  $\chi$ takes values in $\cO$. Let $\mf{m} \subset \cB(U(N), \cO) $ be the maximal ideal which is the kernel of the composition of $\chi$ with  $\cO \rightarrow k$, and $\bs_{\mf{m}}$ the Galois representation corresponding to $\mf{m}$ under Theorem \ref{BHKT}. After possibly further enlarging $E$, we can assume that $\bs_{\mf{m}}$ takes values in $\widehat{G}(k)$.

With Definition \ref{conditions} in mind, we make the following assumptions on $\bs_{\mf{m}}$:
\begin{enumerate}
\item $\ell \nmid \# W$. This implies in particular that $\ell$ is a very good characteristic for $\hG$.
\item The subgroup $Z_{\hG^\text{ad}}(\overline{\sigma}_\ffrm(\Gamma_{K,|N|}))$ of $\hG_k$ is scheme-theoretically trivial.
\item  The representation $\overline{\sigma}_\ffrm$ is absolutely $\hG$-irreducible. 
\item The subgroup $\overline{\sigma}_\ffrm(\Gamma_{K(\zeta_\ell), |N|})$ of $\hG(k)$ is $\hG$-abundant. 
\end{enumerate}

In addition, we choose a global deformation problem $\cD := (|N|, \{\cD_v\}_{v \in |N|})$ for $\bs_{\mf{m}}$ such that all places of ramification for $\bs_{\mf{m}}$ are contained in $|N|$, and: 
\begin{enumerate}
\item[(v)] For each $v \in |N|$, the local deformation problem $\cD_v$ is \emph{balanced} in the sense of Definition \ref{minimal_deformation_problems} and unrestricted in the sense of Example \ref{ex: h^2=0}. (So we require that $H^2(\Gamma_{K_v}, \ad \ol{\sigma}) = 0$ for $v \in |N|$.) 
\end{enumerate}

\begin{remark}Points (i)-(iv) are copied from \cite[\S 8.4]{BHKT}, where their roles are explained. Point (v) allows us to invoke Proposition \ref{minimal_taylor_wiles_data} to produce Taylor-Wiles primes for patching. 
\end{remark}

Let $R_{\overline{\sigma}_\ffrm, |N|, \cD}$ be the corresponding global deformation ring (Proposition \ref{repcD}). From \cite[Theorem 8.5]{BHKT} we have a pseudocharacter $\Theta_{U(N)} = (\Theta_{U(N),n})_{n \geq 1}$ valued in $\cB(U(N), \cO)$, which factors through the quotient $\Gamma_K \to \Gamma_{K, |N|}$. Write $\Theta_{U(N), \mf{m}}$ for the projection of the pseudocharacter $\Theta_{U(N)}$ to $\cB(U(N), \cO)_{\mf{m}}$, and $\sigma^{\mrm{univ}} \co \Gamma_K \rightarrow \wh{G}(R_{\overline{\sigma}_\ffrm, |N|, \cD})$ for a representative of the universal deformation. 

The following Lemma is a variant of \cite[Lemma 8.19]{BHKT}.
\begin{lemma}\label{lem_existence_of_map_R_to_B}
There is a unique morphism $f_\ffrm : R_{\overline{\sigma}_\ffrm, |N|, \cD} \to \cB(U(N), \cO)_\ffrm$ of $\cO$-algebras such that $f_{\ffrm, \ast} \tr \sigma^{\mrm{univ}} = \Theta_{U(N), \ffrm}$. It is surjective.
\end{lemma}
\begin{proof} The existence and uniqueness of $f_\ffrm$ follows  from \cite[Theorem 4.10]{BHKT}. It is surjective because the ring $\cB(U(N), \cO)_\ffrm$ as being generated by the excursion operators $S_{I, (\gamma_i)_{i \in I}, f}$, each of which is explicitly realized as the image of the element  $f(\sigma^{\mrm{univ}}(\gamma_i)_{i \in I}) \in R_{\overline{\sigma}_\ffrm, |N|, \cD}$.
\end{proof}

We can now  prove the following generalization of \cite[Theorem 8.20]{BHKT} that allows for ramification: 

\begin{theorem}\label{thm_R_equals_B}  
Suppose $\mf{m}$ satisfies assumptions (i)-(v) as above. Then $C_{\text{cusp}}(U(N), \cO)_\ffrm$ is a free $R_{\overline{\sigma}_\ffrm, |N|, \cD}$-module, $f_\ffrm$ is an isomorphism, and $R_{\overline{\sigma}_\ffrm, |N|, \cD}$ is a complete intersection $\cO$-algebra. 
\end{theorem}

\begin{proof}[Proof of Theorem \ref{thm_R_equals_B}]
The proof is an implementation of the Taylor--Wiles method. The argument runs similarly to \cite[Theorem 8.20]{BHKT}, but we give a summarized version with the necessary changes for our situation, and only sketching the parts which are the same. 

In fact, since the unramified case is established in \emph{loc.~cit.}, it will be convenient for us to focus our attention on the case of non-trivial level structure -- this will allow to simplify certain aspects of the argument. 
Therefore, we may and do assume for the rest of the proof that $N$ is non-empty. For example, \cite[Theorem 8.17]{BHKT} simplifies in this situation to the following statement. 

\begin{lemma}\label{lem: freeness_over_rings_of_diamond_operators}
Let $U = \prod_v U_v$ be an open compact subgroup of $G(\widehat{\cO}_K)$ contained in $U(N)$ for some non-empty $N$. Let $V \subset U$ be an open normal subgroup such that $U/V$ is abelian of $\ell$-power order. Let $\ffrm \subset \cB(W, \cO)$ be a maximal ideal such that $\overline{\sigma}_\ffrm$ is $\hG$-irreducible. Then $C_{\text{cusp}}(W, \cO)_\ffrm$ is a finite free $\cO[U/W]$-module.
\end{lemma}

\begin{proof}
The finiteness is proved exactly as in \cite[Theorem 8.17]{BHKT}. For the freeness, it suffices to show that $U/W$ acts freely on $G(K) \backslash G(\A_K) / W$. An equivalent formulation is that if we regard $[G(K) \backslash G(\A_K) / U]$ as a groupoid, then all of its stabilizer groups map trivially to $U/W$. By Weil's uniformization theorem (explained for example in \cite[\S 1.2.4]{GL19}), we may interpret the groupoid $[G(K) \backslash G(\A_K) / U]$ as a subgroupoid of the $\F_q$-points of $\Bun_G^U$, the moduli stack of $G$-bundles on $X$ with level structure corresponding to $U$. For any $v \in X$ and any $G$-bundle $\cE$ over $X$, the restriction map $\Aut(\cE) \rightarrow \Aut(\cE|_{X_v^{\wedge}})$ is injective, so the kernel of $\Aut(\cE) \rightarrow \Aut(\cE_v)$ is pro-unipotent, hence its $\F_q$-points form a pro-$p$ group. Therefore, the automorphism groups of the points of the groupoid $[G(K) \backslash G(\A_K) / U]$ are pro-$p$, so they can only map trivially to $U/W$. 
\end{proof}


We begin by preparing the usual setup for patching. Recall that we have fixed a choice $T \subset B \subset G$ of split maximal torus and Borel subgroup of $G$. For a Taylor-Wiles datum $(Q, \{\varphi_v\}_{v \in Q})$, as in \S \ref{ssec: T-W}, with $Q$ disjoint from $|N|$, we introduce the following notation. Extend the global deformation problem $\cD = (|N|, \{\cD_v\}_{v \in |N|})$ to $\cD_Q := (|N| \cup Q, \{\cD_{Q,v}\}_{v\in |N| \cup Q})$ by setting $\cD_{Q, v} = \cD_v$ if $v \in |N|$ and $\cD_{Q,v} = \cD_{\ol{\sigma},v}^{\square}$ (the unrestricted deformation problem) if $v \in Q$. Define $\Delta_Q$ to be the maximal $\ell$-power order quotient of the group $\prod_{v \in Q} T(k(v))$. Using local class field theory, the action of the universal deformation ring for the restriction of $\overline{\sigma}_\ffrm$ to the tame inertia groups at places in $Q$ equips $R_{\overline{\sigma}_\ffrm, |N| \cup Q, \cD_Q}$ with a $\cO[\Delta_Q]$-algebra structure. Writing $\fra_Q \subset \cO[\Delta_Q]$ for the augmentation ideal, we have a canonical isomorphism $R_{\overline{\sigma}_\ffrm, |N| \cup Q, \cD_Q} / \fra_Q \cong R_{\overline{\sigma}_\ffrm, |N|, \cD}$.

We now define the relevant level structures for patching. We define open compact subgroups $U_1(Q) \subset U_0(Q) \subset U(N)$ as follows:
\begin{itemize}
\item $U_0(Q) = \prod_v U_0(Q)_v$, where $U_0(Q)_v = U(N)_v = \ker ( G(\cO_{K_v}) \rightarrow G(\cO_{n_v \cdot v}))$ if $v \not\in Q$, and $U_0(Q)_v$ is an Iwahori group if $v \in Q$. 
\item $U_1(Q) = \prod_v U_1(Q)_v$, where $U_1(Q)_v = U(N)_v$ if $v \not\in Q$, and $U_1(Q)_v$ is the maximal pro-prime-to-$\ell$ subgroup of $U_0(Q)_v$ if $v \in Q$.
\end{itemize}
Thus $U_1(Q) \triangleleft U_0(Q)$ is a normal subgroup, and there is a canonical isomorphism $U_0 / U_1 \cong \Delta_Q$. 

We now need to define auxiliary spaces of modular forms. We define $H'_0 = C_{\text{cusp}}(U(N), \cO)_\ffrm$. There are surjective maps $\cB(U_1(Q), \cO) \surj \cB(U_0(Q), \cO) \surj \cB(U(N), \cO)$, and we write $\ffrm$ as well for the pullback of $\ffrm \subset \cB(U(N), \cO)$ to these two algebras. As in Lemma \ref{lem_existence_of_map_R_to_B}, we have surjective morphisms
\[ R_{\overline{\sigma}_\ffrm, |N| \cup Q, \cD_Q} \surj \cB(U_1(Q), \cO)_\ffrm \surj \cB(U_0(Q), \cO)_\ffrm. \]
We define $H'_{Q, 1} := C_{\text{cusp}}(U_1(Q), \cO)_\ffrm$ and $H'_{Q, 0} := C_{\text{cusp}}(U_0(Q), \cO)_\ffrm$. 

We discuss how to cut $H'_0$ out of $H'_{Q,0}$ as a direct summand, and an analogous construction for $H'_{Q,1}$. This will make use of the Hecke algebras $\cH_{U_0(Q),v}$ at a place $v \in Q$ (where it is the Iwahori-Hecke algebra by our choice of level structure) and $\cH_{U_1(Q),v}$. The ``translation part'' of Bernstein's presentation of the Iwahori-Hecke algebra is an embedding $\cO[X_\ast(T)] \to \cH_{U_0(Q)_v}$, whose action on $H'_{Q,0}$ for $v \in Q$ induces an $R_{\overline{\sigma}_\ffrm, |N| \cup Q, \cD_Q}[ \prod_{v \in Q} X_\ast(T) ]$-module structure on $H'_{Q, 0}$. We write $\frn_{Q, 0} \subset \cO[ \prod_{v \in Q} X_\ast(T) ]$ for the maximal ideal which is associated to the tuple of characters ($v \in Q$):
\begin{equation}\label{eqn_tuple_of_characters_associated_to_TW_datum} \varphi_v^{-1} \circ \overline{\sigma}_\ffrm|_{W_{K_v}} : W_{K_v} \to \hT(k), 
\end{equation}
Then $H'_{Q, 0, \frn_{Q, 0}}$ is a direct factor $R_{\overline{\sigma}_\ffrm, |N| \cup Q, \cD_Q}$-module of $H'_{Q, 0}$, and there is a canonical isomorphism $H'_{Q, 0, \frn_{Q, 0}} \cong H'_0$ of $R_{\overline{\sigma}_\ffrm, |N| \cup Q, \cD_Q}$-modules.\footnote{This uses that for $v \in Q$, the stabilizer in the Weyl group of the regular semisimple element $\overline{\rho}_\ffrm(\Frob_v) \in \hG(k)$ is trivial, which is equivalent to the condition that the centralizer in $\hG_k$ of $\overline{\rho}_\ffrm(\Frob_v)$ is connected, which is part of the definition of a Taylor--Wiles datum.}

Similarly, for $v \in Q$ we write $T(K_v)_\ell$ for the quotient of $T(K_v)$ by its maximal pro-prime-to-$\ell$ subgroup. Then there is a structure of $R_{\overline{\sigma}_\ffrm, |N| \cup Q, \cD_Q}[ \prod_{v \in Q} T(K_v)_\ell]$-module on $H'_{Q, 1}$, where the copy of $T(K_v)_\ell$ corresponding to $v \in Q$ acts via the analogous embedding $\cO[T(K_v)_\ell] \to \cH_{U_1(Q)_v}$. We write $\frn_{Q, 1} \subset \cO[ \prod_{v \in Q} T(K_v)_\ell ]$ for the maximal ideal which is associated to the tuple of characters (\ref{eqn_tuple_of_characters_associated_to_TW_datum}). Then $H'_{Q, 1, \frn_{Q, 1}}$ is a direct factor $R_{\overline{\sigma}_\ffrm, |N| \cup  Q, \cD_Q}$-module of $H'_{Q, 1}$, and the two structures of $\cO[\Delta_Q]$-module on $H'_{Q, 1, \frn_{Q, 1}}$, one arising from the homomorphism $\cO[\Delta_Q] \to R_{\overline{\sigma}_\ffrm, |N| \cup Q, \cD_Q}$ and the other from the homomorphism $\cO[\Delta_Q] \to \cO[\prod_{v \in Q} T(K_v)_\ell]$, are the same.

We need two key properties of the modules $H'_{Q, 1}$ (with fixed $j$  as in Proposition \ref{minimal_taylor_wiles_data}) : 
\begin{itemize}
\item The natural inclusion 
\[ C_{\text{cusp}}(U_0(Q), \cO) \subset C_{\text{cusp}}(U_1(Q), \cO) \]
induces an identification $H'_{Q, 0, \frn_{Q, 0}} = (H'_{Q, 1, \frn_{Q, 1}})^{\Delta_Q}$. 
\item By Lemma \ref{lem: freeness_over_rings_of_diamond_operators}, we then have that $(H'_{Q, 1})^{\ell^{j} \Delta_Q}$ is a free $\cO[\Delta_Q / \ell^{j} \Delta_Q]$-module. This property implies in turn that $(H'_{Q, 1, \frn_{Q, 1}})^{\ell^{j} \Delta_Q}$ is a free $\cO[\Delta_Q / \ell^{j} \Delta_Q]$-module. Observe that $\Delta_Q / \ell^{j} \Delta_Q \cong (\bbZ / \ell^{j} \bbZ)^{\oplus r \# Q}$, where $r = \rank \hG$.
\end{itemize}

For patching it is a bit more convenient to work with $H_Q := \Hom_\cO((H'_{Q, 1, \frn_{Q, 1}})^{\ell^{j} \Delta_Q}, \cO)$ and $H_0 := \Hom_\cO(H'_0, \cO)$. These are finite free $\cO$-modules, which we endow with their natural structures of $R_{\overline{\sigma}_\ffrm, |N| \cup Q, \cD_Q} \otimes_{\cO[\Delta_Q]} \cO[\Delta_Q / \ell^{j} \Delta_Q]$-module and $R_{\overline{\sigma}_\ffrm, |N|, \cD_Q}$-module, respectively, by dualization. We can summarize the preceding discussion as follows:
\begin{itemize}
\item The module $H_Q$ is a finite free $\cO[\Delta_Q / \ell^{j} \Delta_Q]$-module, where $\cO[\Delta_Q / \ell^{j} \Delta_Q]$ acts via the algebra homomorphism 
\[ \cO[\Delta_Q / \ell^{j} \Delta_Q] \to R_{\overline{\sigma}_\ffrm, |N| \cup Q, \cD_Q} \otimes_{\cO[\Delta_Q]}  \cO[\Delta_Q / \ell^{j} \Delta_Q]. \]
\item There is a natural surjective map $H_Q \surj H_0$, which factors through an isomorphism $(H_Q)_{\Delta_Q} \xrightarrow{\sim} H_0$, and is compatible with the isomorphism $R_{\overline{\sigma}_\ffrm, |N| \cup Q, \cD_Q} / \fra_Q \cong R_{\overline{\sigma}_\ffrm, |N|, \cD}$.
\end{itemize}

Let $h := h^1_{\cD}(\Gamma_{K, |N|}, \ad \bs)$. By Proposition \ref{minimal_taylor_wiles_data}, we can find for each $j \geq 1$ a Taylor--Wiles datum $(Q_j, \{ \varphi_v \}_{v \in Q_j})$ which satisfies the following conditions:
\begin{itemize}
\item $Q_j$ is disjoint from $|N|$. 
\item For each $v \in Q_N$, we have $q_v \equiv 1 \text{ mod }\ell^{j}$ and $\# Q_j= h$.
\item There exists a surjection $\cO[[x_1, \ldots, x_g]]\surj R_{\overline{\sigma}_\ffrm, |N| \cup  Q_j, \cD}$, where $g = hr$.
\end{itemize}
Define $R_\infty = \cO \llbracket X_1, \dots, X_g \rrbracket$. The situation is summarized in the diagram below
\[
\begin{tikzcd}[column sep = tiny]
R_{\infty}\ar[rrrr, twoheadrightarrow] & & & & R_{\overline{\sigma}_\ffrm, |N| \cup Q_j, \cD} \ar[d, twoheadrightarrow] & \acts & H_{Q_j}  \ar[d, twoheadrightarrow, "(-)_{\Delta_{Q_j}}"]\\
\cO[\Delta_Q/\ell^j \Delta_Q] \ar[urrrr] & & & & R_{\overline{\sigma}_\ffrm, |N|, \cD} & \acts  & H_0
\end{tikzcd}
\]

We now patch these objects together. This involves quotienting the objects in the diagram by open ideals to get diagrams of Artinian objects. Then because there are only finitely many isomorphisms classes of such diagrams, one can pass to the inverse limit. The details of this process are the same as in \cite[p.48-49]{BHKT}, so we skip to the conclusion. Define $\Delta_\infty := \bbZ_\ell^g$, $\Delta_j := \Delta_\infty / \ell^{j}\Delta_\infty$, $S_\infty := \cO\llbracket \Delta_\infty \rrbracket$,  $\frb_j := \ker (S_\infty \to \cO[[\Delta_j]])$, $\frb_0 := \ker( S_\infty \to \cO)$. By patching, we have the following objects. 

\begin{itemize}
\item $R^\infty$, a complete Noetherian local $\cO$-algebra with residue field $k$, which is equipped with structures of $S_\infty$-algebra and a surjective map $R_\infty \surj R^\infty$.
\item $H_\infty$, a finite $R^\infty$-module.
\item $\alpha_\infty$, an isomorphism $R^\infty / \frb_0 \cong R_{\overline{\sigma}_\ffrm, |N|, \cD}$.
\item $\beta_\infty$, an isomorphism $H_\infty / \frb_0 \cong H_0$. 
\end{itemize}
These objects have the following additional properties:
\begin{itemize}
\item $H_\infty$ is free as an $S_\infty$-module.
\item The isomorphisms $\alpha_\infty$, $\beta_\infty$ are compatible with the structure of $R_{\overline{\sigma}_\ffrm, |N|, \cD}$-module on $H_0$.
\end{itemize}

The situation is summarized in the diagram below
\[
\begin{tikzcd}[column sep = tiny]
R_{\infty}\ar[rrrr, twoheadrightarrow] & & & &  R^{\infty} \ar[d, twoheadrightarrow] &  \acts & H_\infty  \ar[d, twoheadrightarrow]\\
S_{\infty} \ar[u, dashed] \ar[urrrr, twoheadrightarrow] &  & & & R_{\overline{\sigma}_\ffrm, |N|, \cD}  & \acts &  H_0
\end{tikzcd}
\]
Since $S_{\infty}$ is topologically free, we may choose the dashed arrow so that the leftmost triangle commutes. 

We find that
\[ \dim R^\infty \geq \depth_{R^\infty} H_\infty \geq \depth_{S_\infty} H_\infty = \dim S_\infty = \dim R_\infty \geq \dim R^\infty, \]
and hence that these inequalities are equalities, $R_\infty \to R^\infty$ is an isomorphism, and (by the Auslander--Buchsbaum formula) $H_\infty$ is also a free $R^\infty$-module. It follows that $H_\infty / \frb_0 \cong H_0$ is a free $R^\infty / \frb_0 \cong R_{\overline{\sigma}_\ffrm, |N|, \cD}$-module, and that $R_{\overline{\sigma}_\ffrm, |N|, \cD}$ is an $\cO$-flat complete intersection. This in turn implies that $C_{\text{cusp}}(U, \cO)_\ffrm \cong \Hom_\cO(H_0, \cO)$ is a free $R_{\overline{\sigma}_\ffrm, |N|, \cD}$-module (since complete intersections are Gorenstein). This completes the proof of the theorem.
\end{proof}

\section{Formal smoothness of local deformations}\label{mwf}

This section is devoted to the proof of the following theorem, which guarantees that for all large enough $\ell$, the unrestricted local condition will be balanced at a place of ramification for any given cuspidal automorphic representation.  

\begin{theorem}\label{almostall}
Let $\Pi \subset C_{\cusp}(G,\overline{\bbQ})$ be a cuspidal automorphic representation of $G(\bbA_K)$ with coefficients in $\overline{\bbQ}$.  For any prime number
$\ell \neq p$ let
$\Pi_\ell = \Pi\otimes_{\ol{\Q}} \overline{\bbQ}_\ell$ denote its base change to $\overline{\bbQ}_\ell$. Suppose $\Pi$ is unramified outside the finite set $S$, and $\sigma_{\ell} \in \ul{\rho}_{\Pi, \ell}$ (notation as in Remark \ref{rem: weak correspondence}) has Zariski-dense image. (As explained in Remark \ref{rem: weak correspondence}, this implies a posteriori that $\sigma_{\Pi, \ell}$ exists, and then $\sigma_{\ell} = \sigma_{\Pi, \ell}$.) Then, letting $$\bs_{\ell}:  \G_{K,S} \ra \hG(\overline{\bbF}_\ell)$$
be a reduction of $\sigma_{\ell}$ mod $\ell$, there is an integer $c(\Pi)$, depending only on $\Pi$, such that, for all $v \in S$ and all $\ell > c(\Pi)$, we have
$$H^2(\Gv, \ad(\bs_{\ell})) = 0.$$ 
\end{theorem}

We concurrently found two very different proofs of Theorem \ref{almostall}, and include them both below, as they both establish general results along the way that may be of independent interest. The proof in \S \ref{ssec: automorphy approach} is based on a lifting theorem (Theorem \ref{thm_existence_of_liftings}) for global Galois representations, which is analogous to \cite[Theorem 4.3.1]{Bar14}; cf. also \cite[Corollary 4.7]{Kha09}. The lifting theorem implies that if $H^2(\Gv, \ad(\bs_{\ell})) \neq 0$, then there is another lift of $\bs_{\ell}$ to characteristic zero which is ``less ramified'' at $v$ than $\sigma_{\ell}$. By the Global Langlands correspondence for $\GL_n$ established by Lafforgue \cite{Laf02}, such a lift must come from another cuspidal automorphic representation -- of which only finitely many contribute at any given level structure --  congruent to $\Pi$ modulo $\ell$. The idea is then that there can only be finitely many such congruences, because otherwise strong Multiplicity One for $\GL_n$ would be violated. We remark that a similar strategy, for an analogous situation over number fields, was used in the proof of \cite[Corollary 7.11]{Gui20}. 
  
 The proof in \S \ref{ssec: geometric approach} is of a more geometric nature. We reduce Theorem \ref{almostall} to a torsion version of the Weight-Monodromy Conjecture, formulated by Ito \cite{Ito20}, for the compactifications constructed by Laurent Lafforgue to prove the global Langlands correspondence for $\GL_n$. Then we prove this torsion Weight-Monodromy Conjecture for the intersection cohomology of any proper variety over an equal characteristic local field (generalizing work of Ito in the smooth proper case \cite{Ito20}). A crucial tool is an integral version of the Decomposition Theorem recently proved by Cadoret-Zheng, which we use to reduce to the smooth and proper case established by Ito. This second argument actually yields a more general statement, Theorem \ref{thm: general H^2 vanishing}.

Although neither argument gives an effective estimate on the constant $c(\Pi)$, they give different starting points for obtaining such an estimate. (That is one reason why it seems valuable to include both arguments.) In the first approach, what must be controlled is the set of primes at which cusp form can have congruences to other cusp forms. In the second approach, what must be controlled has to do with the torsion in the integral intersection cohomology of L. Lafforgue's compactifications of moduli spaces of shtukas for $\GL_n$.

\subsection{Approach based on lifting theorems for Galois representations}\label{ssec: automorphy approach}
First we show that Theorem \ref{almostall} follows from the same statement when $G = \GL_n$. Pick a faithful irreducible representation of $G$ on a finite projective $\Z$-module, say of rank $n > 1$. This induces an injection $\mf{g} \inj \mf{gl}_n$, which is split over $\Z_{\ell}$ for all suffiicently large $\ell$ and therefore induces an injection $H^2(\Gv, \ad(\bs_{\ell})) \inj H^2(\Gv, \mf{gl}_n \circ \bs_{\ell})$. This reduces the claim to the case $\wh{G} = \GL_n$, so for the rest of this subsection we focus on the case $G = \GL_n$. 

\subsubsection{Some linear algebra}\label{subsec_linear_algebra}

We recall the following result of Deligne \cite[Proposition 1.6.1]{Del80}.
\begin{proposition}\label{prop_monodromy_filtration}
Let $k$ be a field and let $V$ be a finite-dimensional $k$-vector space, equipped with a nilpotent endomorphism $N : V \to V$. Then there exists a unique increasing filtration $M_\bullet$ of $V$ such that $N M_i \subset M_{i-2}$ and, for each $k \geq 0$, $N^k$ induces isomorphisms $\operatorname{gr}_k M_\bullet \to \operatorname{gr}_{-k} M_\bullet$.
\end{proposition}
If $N$ is a nilpotent endomorphism of an $n$-dimensional vector space, we will write $\operatorname{Jord}(N)$ for the partition of $n$ given by the sizes of the Jordan blocks of $N$. We recall that there is a partial order on partitions $n$ which can be viewed as  corresponding to the closure ordering of nilpotent orbits in the adjoint representation of $\GL_n$. 

Let $v$ be a place of $K$, and suppose given a prime $\ell > n$ not dividing $q$. Let $E / \bbQ_\ell$ be a finite extension with residue field $k_E = \cO_E / (\varpi_E)$, and suppose given a continuous representation $\rho : G_{K_v} \to \GL_n(\cO_E)$. Fix a choice of (geometric) Frobenius lift $\phi_v \in G_{K_v}$, and $t_\ell \in G_{K_v}$ a generator for the $\ell$-part of tame inertia. Let us say that $\rho$ is \emph{good} if the following conditions are satisfied:
\begin{itemize}
\item $\rho$ is unipotently ramified. Let $N = \log(\rho(t_\ell)) \in M_n(\cO_E)$; then $\phi_v N \phi_v^{-1} = q^{-1} N$. 
\item The eigenvalues of $\rho(\phi_v)$ all lie in $E$.
\item $\rho$ is pure of weight 0. This means that the eigenvalues of $\rho(\phi_v)$ on $V = E^n$ are $q_v$-Weil numbers and that, writing $V_i \subset V$ for the subspace where the eigenvalues are of weight $i$, for each $k \geq 0$ $N^k$ induces an isomorphism $V_k \to V_{-k}$. (Note that the filtration associated to $N$ by Proposition \ref{prop_monodromy_filtration} is then $M_k = \oplus_{i \leq k} V_i$.)
\item Let $f_i(X) = \det(X - \rho(\phi_v)|_{V_i}) \in \cO_E[X]$. Then for each $i \neq j$, the reduced polynomials $\overline{f}_i(X), \overline{f}_j(X) \in k_E[X]$ are coprime. 
\item For each $i, j$, the polynomials $\overline{f}_i(X)$, $\overline{f}_j(q_v^{-1} X) \in k_E[X]$ fail to be coprime only if $j = i + 2$. 
\end{itemize}
Let $L = \cO_E^n$, so that $V = L \otimes_{\cO_E} E$, and define $L_i = V_i \cap L$. Hensel's lemma implies that if $\rho$ is good then $L = \oplus_{i \in \bbZ} L_i$.
\begin{lemma}\label{lem_Jord_and_H2}
Suppose that $\rho$ is good, and let $\overline{N} \in M_n(k_E)$ denote the reduction of $N$ modulo $\varpi_E$. If $\operatorname{Jord}(\overline{N}) = \operatorname{Jord}(N)$, then $H^2(K_v, \ad \overline{\rho}) = 0$.
\end{lemma}
\begin{proof}
The equality $\operatorname{Jord}(\overline{N}) = \operatorname{Jord}(N)$ is equivalent to the equalities $\dim_E \ker N^k = \dim_{k_E} \ker \overline{N}^k$ for each $k \geq 0$. Let $\overline{L} = L \otimes_{\cO_E} k_E$ and $\overline{L}_i = L_i \otimes_{\cO_E} k_E$. We first claim that for each $k \geq 0$, $\overline{N}^k$ induces an isomorphism $\overline{L}_k \to \overline{L}_{-k}$. Equivalently, the intersection $(\ker \overline{N}^k) \cap \overline{L}_k$ is 0. 

The equality $\dim_E \ker N^k = \dim_{k_E} \ker \overline{N}^k$ implies that in fact $\ker \overline{N}^k = (\ker N^k \cap L) \otimes_{\cO_E} k_E$. Since $\ker N^k$ is invariant under $\Ad \rho(\phi_v)$, we find that in fact
\[ (\ker \overline{N}^k) \cap \overline{L}_k =  (\ker N^k \cap L_k) \otimes_{\cO_E} k_E = 0. \]
This establishes the claim. By linear algebra, we then have a decomposition
\[ \overline{L} = \oplus_{i \in \bbZ} \oplus_{j=0}^i \overline{N}^j \overline{L}(i), \]
where $\overline{L}(i) = \overline{L}_i \cap (\ker \overline{N}^{i+1})$.

By Tate duality, we need to show that $H^0(K_v, \ad \overline{\rho}(1)) = \Hom_{G_{K_v}}(\overline{L}, \overline{L}(1)) = 0$. An element of this Hom space determines a linear map $F : \overline{L} \to \overline{L}$ which commutes with $\overline{N}$ and satisfies $F(\overline{L}_i) \subset \overline{L}_{i+2}$ for each $i \in \bbZ$. To show the Hom space vanishes, it is enough to show that $F$ annihilates each $\overline{L}(i)$. However, if $x \in \overline{L}(i)$ then we find $F(x) \in \overline{L}_{i+2}$ and 
\[ \overline{N}^{i+1} F(x) = F( \overline{N}^{i+1} x ) = 0. \]
The restriction of $\overline{N}^{i+1}$ to $\overline{L}_{i+2}$ is injective, so this forces $F(x) = 0$, hence $F = 0$, as required.
\end{proof}

\subsubsection{Local lifting rings}

Let $\ell > n $ be a prime not dividing $q$ and let $E / \bbQ_\ell$ be a finite extension with residue field $k_E = \cO_E / (\varpi_E)$. Let $v$ be a place of $K$ and fix a continuous homomorphism $\overline{\rho}_v : G_{K_v} \to \GL_n(k_E)$ and a continuous character $\chi : G_K \to \cO_E^\times$ such that $\det \overline{\rho}_v = \chi|_{G_{K_v}}$. Let $\operatorname{CNL}_{\cO_E}$ denote the category of complete Noetherian local $\cO_E$-algebras with residue field $k_E$. Then the functor $L_v^\chi$ which associates to each $A \in \operatorname{CNL}_{\cO_E}$ the set of lifts $\rho_{v, A} : G_{K_v} \to \GL_n(A)$ of $\overline{\rho}_v$ of determinant $\chi|_{G_{K_v}}$ is represented by an object $R_v^{\square, \chi} \in \operatorname{CNL}_{\cO_E}$. This ring has a pleasant geometry:
\begin{proposition}\label{prop_dimension_of_components}
	Let assumptions be as above. Then:
	\begin{enumerate}
	\item 	The ring $R_v^{\square, \chi}$ is a reduced local complete intersection. 
	\item For each minimal prime $\mathfrak{q} \subset R_v^{\square, \chi}$, the quotient $R_v^{\square, \chi} / \mathfrak{q}$ is a domain of Krull dimension $n^2$. 
	\item Let $\mathfrak{q} \subset  R_v^{\square, \chi}$ be a minimal prime, and let $L_v^{\chi, \mathfrak{q}} \subset L_v^{\chi}$ be the corresponding subfunctor. Then $L_v^{\chi, \mathfrak{q}}$ is a local deformation problem, in the sense that for each $A \in \operatorname{CNL}_{\cO_E}$, the subset $L_v^{\chi, \mathfrak{q}}(A) \subset L_v^\chi(A)$ is invariant under the conjugation action of the group $\ker(\GL_n(A) \to \GL_n(k_E))$.
	\end{enumerate}
\end{proposition}
\begin{proof}
The first two points follow from \cite[Proposition 2.13]{Yao21}. The third is established in \cite[\S 1.3]{Bar14} in  the case where the determinant is not fixed and $K_v$ is replaced by a finite extension of $\bbQ_p$ ($p \neq \ell$); the same proof applies here. 
\end{proof}
We will use the following lemma.
\begin{lemma}\label{lem_minimal_lifting_ring}
	With assumptions as above, suppose further that $\ell > n$ and that $\rho_v$ is unipotently ramified. Then there exists a minimal prime $\mathfrak{q} \subset R_v^{\square, \chi}$ such that $R_v^{\square, \chi} / \mathfrak{q}$ is formally smooth over $\cO_E$ and such that for every homomorphism $R_v^{\square, \chi} / \mathfrak{q} \to \overline{\bbQ}_\ell$, corresponding to a lifting $\rho_v : G_{K_v} \to \GL_n(\overline{\bbZ}_\ell)$ of $\overline{\rho}_v$ of determinant $\chi|_{G_{K_v}}$, the following conditions are satisfied:
	\begin{enumerate}
		\item $\rho_v$ is unipotently ramified.
		\item If $t_{\ell}\in I_{K_v}(\ell)$ denotes a generator of the $\ell$-part of tame inertia, then $\operatorname{Jord}(\log(\rho_v(t_\ell))) = \operatorname{Jord}(\log(\overline{\rho}_v(t_\ell)))$.
	\end{enumerate}
\end{lemma}
\begin{proof}
	We can take $R_v^{\square, \chi} / \mathfrak{q}$ to be the functor parametrizing liftings which are minimally ramified, in the sense of \cite[Definition 2.4.14]{CHT08}, and of determinant $\chi|_{G_{K_v}}$. The claimed properties follow from \cite[Lemma 2.4.19]{CHT08}. 
\end{proof}

\subsubsection{Existence of lifts}

\begin{theorem}\label{thm_existence_of_liftings}
	Let $\ell > n $ be a prime such that $(\ell, 2q) = 1$, and let $\overline{\rho} : G_K \to \GL_n(\overline{\bbF}_\ell)$ be a continuous representation such that $\overline{\rho}|_{G_{K \cdot \overline{\bbF}_q}}$ is irreducible. Let $S$ be a finite set of places of $K$, containing the places at which $\overline{\rho}$ is ramified, and fix the following data:
	\begin{enumerate}
	\item A finite extension $E / \bbQ_\ell$ inside $\overline{\bbQ}_\ell$.
		\item A character $\chi : G_K \to \cO_E^\times$ of finite order lifting $\det \overline{\rho}$, unramified outside $S$.
		\item For each $v \in S$, a minimal prime $\mathfrak{q}$ of the universal lifting ring $R_v^{\square, \chi}$.
	\end{enumerate}
	Then we can find a continuous lifting $\rho : G_K \to \GL_n(\overline{\bbQ}_\ell)$ of $\overline{\rho}$ satisfying the following conditions:
	\begin{enumerate}
		\item $\det \rho = \chi$.
		\item For each place $v\not\in S$ of $K$, $\rho|_{G_{K_v}}$ is unramified.
		\item For each place $v \in S$ of $K$, $\rho|_{G_{K_v}}$ defines a homomorphism $R_v^{\square, \chi} / \mathfrak{q} \to \overline{\bbQ}_\ell$.
	\end{enumerate}
\end{theorem}
\begin{proof}
	Let $D : \operatorname{CNL}_{\cO_E} \to \operatorname{Sets}$ denote the functor of deformations $\rho_A : G_K \to \GL_n(A)$ of $\overline{\rho}$ satisfying the following conditions: 	
	\begin{enumerate}
		\item $\det \rho_A = \chi$.
		\item If $v \in S$, then $\rho_A|_{G_{K_v}}$ defines a map $R_v^{\square, \chi} / \mathfrak{q} \to A$.
		\item If $v \not\in S$, then $\rho_A|_{G_{K_v}}$ is unramified.
	\end{enumerate}
	Then $D$ is represented by an object $R \in \operatorname{CNL}_{\cO_E}$, and a standard computation in Galois cohomology and using Proposition \ref{prop_dimension_of_components} (see \cite[Corollary 5.4]{Yao21}) shows that $R$ has Krull dimension $\dim R \geq 1$. On the other hand, $R$ is a finite $\cO_E$-algebra (see \cite[Theorem 5.5]{Yao21}, where this is deduced from work of de Jong and Gaitsgory). It follows that $R$ is an $\cO_E$-flat complete intersection and in particular that there exists a homomorphism $R \to \overline{\bbQ}_\ell$, which implies the existence of a lift $\rho : G_K \to \GL_n(\overline{\bbZ}_\ell)$ with the desired properties.
\end{proof}

\subsubsection{Minimality}

\begin{theorem}
	Let $\pi$ be a cuspidal automorphic representation of $\GL_n(\A_K)$ over $\overline{\bbQ}$ of finite central character and such that the associated compatible system $(r_\lambda(\pi))_\lambda$ is connected, in the sense of \cite{BGP}. Let $S$ be a finite set of places of $K$, including the set of places at which $\pi$ is ramified. Then there exists an integer $N > \max(2, q, n)$ such that for each prime number $\ell > N$ and for each place $\lambda$ of $\overline{\bbQ}$ of residue characteristic $\ell$, the following conditions are satisfied:
	\begin{enumerate}
		\item $\overline{\rho}_\lambda|_{G_{K \cdot \overline{\bbF}_q}}$ is irreducible.
		\item For each $v \in S$, $H^2(K_v, \ad \overline{\rho}_\lambda) = 0$.
	\end{enumerate}
\end{theorem}
\begin{proof}
	Let $\chi$ denote the central character of $\pi$, which we identify with a character $\chi : G_K \to \overline{\bbQ}^\times$. After replacing $K$ by a finite Galois extension, we may assume that for each place $v \in S$, $\pi_v$ is unipotently ramified. We can associate to each $v \in S$ a partition $\operatorname{Jord}(\pi_v)$ of $n$, namely the one given by the Jordan decomposition of the nilpotent part of the Weil--Deligne representation $\operatorname{\rec}_{K_v}(\pi_v)$. 
	
	By \cite[Corollary 6.12]{BGP}, we can find an integer $N_0 > \max(2, q, n)$ such that for each prime number $\ell > N$ and for each place $\lambda$ of $\overline{\bbQ}$ of residue characteristic $\ell$, $\overline{\rho}_\lambda|_{G_{K \overline{\bbF}_q}}$ is irreducible. Since $S$ is finite, it suffices to fix a single $w \in S$ and show that the number of primes $\ell > N_0$ such that there is a place $\lambda | \ell$ of $\overline{\bbQ}$ such that $H^2(K_w, \ad \overline{\rho}_\lambda) \neq 0$ is finite. 
	
	Let $L \subset \overline{\bbQ}$ be a number field such that the compatible system $(r_\lambda(\pi))_\lambda$ is $L$-rational and such that the roots of $f(X) = \det(X - \rec_{K_w}(\pi_w)(\phi_w))$ lie in $L$. By the Ramanujan conjecture (established in \cite[Theorem VI.10]{Laf02}), $\rec_{K_w}(\pi_w)$ is pure of weight 0, and we may factor $f(X) = \prod_i f_i(X)$, where the roots of $f_i(X)$ are $q_w$-Weil numbers of weight $i$. The polynomials $f(X)$, $f_i(X)$ lie in $\cO_L[1/q, X]$. If $i \neq j$ then the polynomials $f_i(X), f_j(X)$ are coprime. If $i \neq j+2$ then the polynomials $f_i(X), f_{j+2}(q_w X)$ are coprime. After increasing $N_0$, we can therefore assume that for each prime $\ell > N_0$ and for each place $\lambda | \ell$ of $\overline{\bbQ}$, the representation $\rho_\lambda|_{G_{K_w}}$ (conjugated to take values in $\GL_n(\cO_E)$ for some $E / \bbQ_\ell$) is good, in the sense of \S \ref{subsec_linear_algebra}. 
	\begin{lemma}
	For each prime number $\ell > N_0$ and place $\lambda | \ell$ of $\overline{\bbQ}$ such that $H^2(K_w, \ad \overline{\rho}_\lambda) \neq 0$, there exists a cuspidal automorphic representation $\pi'$ of $\GL_n(\A_K)$ over $\overline{\bbQ}$ satisfying the following conditions:
		\begin{enumerate}
			\item $\pi'$ has central character $\chi$ and is unramified outside $S$.
			\item For each place $v \in S$, $\pi'_v$ is unipotently ramified. 
			\item $\operatorname{Jord}(\pi'_w) < \operatorname{Jord}(\pi_w)$.
			\item There is an isomorphism $\overline{r_\lambda(\pi)} \cong \overline{r_\lambda(\pi')}$.
		\end{enumerate}
	\end{lemma}
	\begin{proof}
We would like to apply Theorem \ref{thm_existence_of_liftings}. We need to specify a minimal prime of each lifting ring $R_v^{\square, \chi}$ ($v \in S$). If $v \neq w$ then we take any component containing the point corresponding to $r_\lambda(\pi)$. If $v = w$ then we take the minimal lifting ring described in Lemma \ref{lem_minimal_lifting_ring}. We then apply Theorem \ref{thm_existence_of_liftings} with these choices to obtain a lift $\rho$ of $\overline{r_\lambda(\pi)}$, which by the global Langlands correspondence for $\GL_n(\A_K)$ corresponds to a cuspidal automorphic representation $\pi'$ of $\GL_n(\A_K)$. The only property that remains to be justified is that $\operatorname{Jord}(\pi'_w) < \operatorname{Jord}(\pi_w)$. By local-global compatibility, we have $\operatorname{Jord}(\pi'_w) = \operatorname{Jord}(\rho|_{G_{K_w}})$ (where we define the Jordan partition of a unipotently ramified representation to be that of the image of a generator of tame inertia). The definition of the minimally ramified deformation problem shows that $\operatorname{Jord}(\rho|_{G_{K_w}}) = \operatorname{Jord}(\overline{r_\lambda(\pi)}|_{G_{K_w}})$. Finally, Lemma \ref{lem_Jord_and_H2} and our hypothesis that $H^2(K_w, \ad \overline{\rho}_\lambda) \neq 0$ show together that $\operatorname{Jord}(\overline{r_\lambda(\pi)}|_{G_{K_w}}) < \operatorname{Jord}(\pi_w)$.
	\end{proof}
	We now use the preceding Lemma to finish the proof of the Theorem. The set of cuspidal automorphic representations of $\GL_n(\A_K)$ over $\overline{\bbQ}$ of central character $\chi$, unramified outside $S$, and unipotently ramified inside $S$, is finite \cite[Corollary 1.2.3]{Har74}. Enlarge the number field $L$ so that each such automorphic representation is in fact defined over $L$. If the number of places $\lambda$ of $L$ for which $H^2(K_w, \ad \overline{\rho}_\lambda) \neq 0$ is infinite, then by the pigeonhole principle we can find an automorphic representation $\pi'$ in this set such that for infinitely many places $\lambda$ of $L$ there is an isomorphism $\overline{r_\lambda(\pi)} \cong \overline{r_\lambda(\pi')}$, and yet $\operatorname{Jord}(\pi'_w) < \operatorname{Jord}(\pi_w)$. In particular, $\pi \neq \pi'$. 
	
	If $v \not \in S$ is a place of $K$ then in particular we find $\det(X - r_\lambda(\pi)(\Frob_v)) \equiv \det(X - r_\lambda(\pi')(\Frob_v)) \text{ mod } \lambda$ for infinitely many places $\lambda$, and therefore $\pi_v \cong \pi'_v$. The strong Multiplicity One theorem then implies that $\pi = \pi'$, contradicting the inequality $\operatorname{Jord}(\pi'_w) < \operatorname{Jord}(\pi_w)$. This concludes the proof. 
\end{proof}

\subsection{Geometric approach based on the Torsion Weight-Monodromy Conjecture}\label{ssec: geometric approach}

We will actually establish the following more general statement, which handles arbitrary compatible systems of irreducible representations (as opposed to specifically the adjoint representation). 

\begin{theorem}\label{thm: general H^2 vanishing}
Let $\{\sigma_{\ell} \co \G_{K,S} \ra \GL_n(\overline{\Q}_{\ell})\}$ be any compatible system of irreducible $n$-dimensional representations. For each $v \in |X|$, Grothendieck's quasi-unipotence theorem implies that there is an open subgroup $J_v \subset I_{K_v}$ such that $\sigma_{\ell}(J_v)$ is unipotent for all $\ell \neq p$. Choose for each $v \in |X|$ a $t_v \in J_v$ projecting to a generator of the pro-$\ell$ part of tame inertia for all but finitely many $\ell \neq p$. Then all but finitely many primes $\ell$ have the properties that:
\begin{itemize}
\item $\bs_{\ell}$ is irreducible and 
\item for all $v$, $\sigma_{\ell}(t_v) - 1$ has torsion-free cokernel on some (equivalently any) lattice. 
\end{itemize}
\end{theorem}

\subsubsection{Some reductions} By assumption $G$ is semi-simple; we may reduce to the case where $G$ is simple since the statement of the Theorem is compatible with taking finite direct sums. 

The representation $\sigma_{\ell} $ is defined over a finite extension $E/\Q_{\ell}$. Let $\cO_{E}$ be the ring of integers of $E$, and $\wt{\sigma}_{\ell}$ be a lattice in $\sigma_{\ell}$, with mod-$\ell$ reduction $\bs_{\ell}$. 

\begin{lemma}
If the cokernel of $(\wt{\sigma}_{\ell}(t_v)-1)^i$ is torsion-free for every $i$, then $H^2(K_v, \bs_{\ell}) = 0$. 
\end{lemma}

\begin{proof}
By \cite[Lemma 3.7]{Ito20}, the torsion-freeness implies that $\Jord(\wt{\sigma}_{\ell}(t_v)-1) = \Jord(\bs_{\ell}(t_v)-1)$ in the notation of Lemma \ref{lem_Jord_and_H2}, and then applying the same proof as for Lemma \ref{lem_Jord_and_H2} gives the conclusion. 
\end{proof}

\begin{lemma}
Let $M$ be a finite flat $\Z_{\ell}$-module and $\gamma$ a $\Z_{\ell}$-linear endomorphism of $M$. If $(M \otimes_{\Z_{\ell}} M)_{\Delta\gamma}$ is torsion-free, then $M_{\gamma}$ is torsion-free. (Here $\Delta \gamma$ refers to the diagonal action of $\gamma$ on the tensor product.) 
\end{lemma}

\begin{proof}
Suppose $m \in M$ represents a non-zero $\ell^k$-torsion element of $M_{\gamma}$. Then $m \otimes m \in M \otimes_{\Z_{\ell}} M$ represents an $\ell^{2k}$-torsion element of $(M \otimes_{\Z_{\ell}} M)_{\Delta \gamma}$. Furthermore, its image under $(M \otimes_{\Z_{\ell}} M)_{\Delta \gamma} \surj M_{\gamma} \otimes_{\Z_{\ell}} M_{\Delta\gamma}$ is non-zero, so this shows that $(M \otimes_{\Z_{\ell}} M)_{\Delta \gamma}$ has a non-zero torsion element. 
\end{proof}

Hence it suffices to show that for all sufficiently large $\ell$, the cokernel of $(\wt{\sigma}_{\ell}(t_v)-1)^i \otimes (\wt{\sigma}_{\ell}(t_v)-1)^i$ is torsion-free for all $i$. The point of casting the problem this way is that this tensor product of Galois representations appears in the cohomology of shtukas, and through this we will be able to transform the question into a purely geometric one. 

Let $n = \dim \frg$. Laurent Lafforgue constructed a proper Deligne-Mumford stack $\ol{\Sht}_N \rightarrow (X-N) \times (X-N)$ (our $\ol{\Sht}_N$ refers to what is called in \cite{Laf02} $\ol{\mrm{Cht}^{n, \ol{P} \leq P}_N}'	 / a^{\Z}$ for sufficiently large Harder-Narasimhan polygon $P$, and appropriate $a$ corresponding to the central character of $\Pi$) such that (because $\sigma_{\ell}$ is irreducible) $\Xi \boxtimes (\sigma_{\ell} \otimes \sigma_{\ell})$ occurs, for some\footnote{Explicitly, $\Xi$ is the $U(N)$ invariants on the cuspidal automorphic representation corresponding to $\sigma_{\ell}$ via L. Lafforgue's Global Langlands Correspondence for $\GL_n$.} Hecke-module $\Xi$, as a Hecke-Galois stable summand of $\pc H^0 (\ol{\Sht}_N, \IC(E))$ \cite[Theorem VI.27]{Laf02}, where the latter refers to the $0$th perverse cohomology of the geometric generic fiber of $\ol{\Sht}_N \rightarrow (X-N) \times (X-N)$.

\begin{corollary}\label{cor: summand}
For any $v \in |X|$, for all but finitely many $\ell$, there exists $N$ disjoint from $v$ such that $\wt{\sigma}_{\ell} \otimes_{\cO_E} \wt{\sigma}_{\ell}$ occurs as a Hecke-stable\footnote{The Hecke algebra does not act on all of $\pc H^0 (\ol{\Sht}_N, \IC(E))$ or $\pc H^0 (\ol{\Sht}_N, \IC(\cO_E))$. There are Hecke correspondences on $\ol{\Sht}_N$ defined by normalization, which act on these cohomology groups, but not in a way satisfying the relations of the Hecke algebra. However, Lafforgue shows that on the ``essential'' summand $H^*_{\mrm{ess}}(\ol{\Q}_{\ell}) \subset \pc H^0 (\ol{\Sht}_N, \IC(\ol{\Q}_{\ell}))$, the Hecke operators satisfy the relations of the Hecke algebra. This is all that matters for us.} summand of $\pc H^0 (\ol{\Sht}_N, \IC(\cO_{E}))$ where $\ol{\Sht}_N \rightarrow (X-N) \times (X-N)$ is a proper morphism of schemes.
\end{corollary}

\begin{proof}
Fixing the Harder-Narasimhan polygon $p$, for all sufficiently large $N$ \cite[Proposition 2.16]{Var04} implies that the Deligne-Mumford stack $\ol{\Sht}_N$ is represented by a scheme. 

Since $\sigma_{\ell}$ is irreducible, the $\sigma_{\ell}$-isotypic part of $\pc H^0 (\ol{\Sht}_N, \IC(\ol{\Q}_{\ell}))$ occurs as a summand of the finite-dimensional subspace $H^*_{\mrm{ess}}(\ol{\Q}_{\ell}) \subset \pc H^0 (\ol{\Sht}_N, \IC(\ol{\Q}_{\ell}))$, where the ``essential part'' is defined by the condition that its Hecke eigenvalues are cuspidal. 

The Hecke operators act on $H^*_{\mrm{ess}}(\ol{\Q}_{\ell})$, and the subalgebra of $\End(H^*_{\mrm{ess}}(\ol{\Q}_{\ell}))$ generated by them forms a finite $\ol{\Q}_{\ell}$-algebra, so we may choose a finite number of Hecke operators $T_1, \ldots, T_k$ that forms a generating set. Each of these finitely many operators on $H^*_{N, \mrm{ess}}$ satisfies a characteristic polynomial which is defined over a number field, hence for all sufficiently large $\ell$ all their eigenvalues lie in $\ol{\Z}_{\ell}$. Since the projection to eigenspaces is given by universal polynomials in the eigenvalues, these projection operators are then defined over $\ol{\Z}_{\ell}$ for all sufficiently large $\ell$, giving a decomposition which then further descends to $\cO_E$. 
\end{proof}

For $\ol{\Sht}_N$ as in Corollary \ref{cor: summand}, by Grothendieck's quasi-unipotence theorem we may choose $t_v$ to act unipotently on all of $H^0 ((\ol{\Sht}_N)_{K_v^s}, \IC(\Z_{\ell}))$. It then suffices to show that for almost all $\ell$, the cokernel of $(t_v-1)^i$ is torsion-free on $H^0 (\ol{\Sht}_N, \IC(\Z_{\ell}))$ for all $i$. We will prove a much more general statement, for any smooth proper variety over $K_v$. 

\subsubsection{Monodromy weight filtrations for torsion intersection cohomology}
Let $Z$ be a proper variety over a characteristic $p$ local field $F$. Then we have an action of $\Gal(F^s/F)$ on its (geometric) intersection cohomology $H^*(Z_{F^s}, \IC_Z(\Z_{\ell}))$. 

There is an open subgroup $J \subset I_{F}$ such that its action is unipotent. Take $t \in J$ such that its projection to the $\ell$-part of tame inertia is non-trivial for all but finitely many $\ell \neq p$.

\begin{proposition}\label{prop: monodromy cokernel}
For almost all $\ell \neq p$, the cokernel of $(t-1)^i \co H^*(Z_{F^s}, \IC_Z(\Z_{\ell})) \rightarrow H^*(Z_{F^s}, \IC_Z(\Z_{\ell}))$ is torsion-free for all $i$. 
\end{proposition}

In the special case where $Z$ is smooth and proper, this is a result of Ito \cite{Ito20}, which ultimately relies on the Weight-Monodromy Conjecture. Our proof will reduce to the smooth and proper case using the following Proposition. 

\begin{proposition}\label{prop: IC direct summand integrally}
Let $f \co Y \rightarrow Z$ be an alteration of varieties over a field $L$, with $Y$ smooth over $L$. Then for all sufficiently large $\ell$, $Rf_* \Z_{\ell}$ has $\IC_Z(\Z_{\ell})$ as a direct summand. 
\end{proposition}

\begin{proof}[Proof of Proposition \ref{prop: monodromy cokernel} assuming Proposition \ref{prop: IC direct summand integrally}]
According to de Jong \cite{dJ96}, we may find an alteration $f \co Y \rightarrow Z$ with $Y$ being smooth and proper. The cokernels in question are then direct summands of the analogous cokernels for $Y$, which by  \cite[combination of Lemma 3.7, Proposition 3.9, and Theorem 3.6]{Ito20} are torsion-free for all sufficiently large $\ell$. 
\end{proof}

\begin{remark}
Our Proposition \ref{prop: IC direct summand integrally} also implies a torsion version of the Weight-Monodromy Conjecture for the intersection cohomology of a proper variety $Z$ with $\F_{\ell}$-coefficients, generalizing the case of smooth and proper $Z$ treated in \cite{Ito20}. 
\end{remark}

The proof of Proposition \ref{prop: IC direct summand integrally} will be given in the next subsection. If we were considering $\Q_{\ell}$-coefficients instead, we would be able to deduce the analogous statement from the Decomposition Theorem (although it is not an entirely trivial deduction when $F$ is not separably closed, since the splitting of the Decomposition Theorem over $F^s$ is not canonical, and does not descend to $F$ in general). For varieties over $\CC$, the statement with $\Z_{\ell}$-coefficients would follow for all sufficiently large $\ell$ from the $\Q_{\ell}$-coefficient version plus the existence of a $\Z$-structure on $Rf_* \Z_{\ell}$ and $\IC(\Z_{\ell})$. However, Proposition \ref{prop: IC direct summand integrally} in the stated generality seemed out of reach until we learned of the recent theorem of Cadoret-Zheng, which gives a version of the Decomposition Theorem with $\Z_{\ell}$-coefficients for all sufficiently large $\ell$. 

\begin{theorem}[{Cadoret-Zheng \cite{CZ}}]\label{Decomposition Theorem}
Let $f \co Y \rightarrow Z$ be a proper surjection of finite type schemes over a field $L$, with $Y$ is smooth over $L$. Then for all sufficiently large $\ell$, we have 
\[
Rf_* \Z_{\ell} \cong \bigoplus_i \pch^i (Rf_* \Z_{\ell})[-i],
\]
with furthermore each $\pch^i (Rf_* \Z_{\ell})$ being torsion-free as a perverse sheaf. 

If furthermore $L$ is separably closed, then for all sufficiently large $\ell$, each $\pch^i (Rf_* \Z_{\ell})$ is isomorphic to a direct sum of intersection complexes of semi-simple local systems.
\end{theorem}

\subsubsection{Perverse sheaves}

We would like to apply Theorem \ref{Decomposition Theorem} to the situation of Proposition \ref{prop: IC direct summand integrally}. A subtlety is that because the splitting provided by the Decomposition Theorem over $L^s$ is not canonical, it does not in general descend to $L$. We must therefore pay careful attention to rationality issues. We will use a trick that we learned from Hansen's blog\footnote{Available at \url{https://totallydisconnected.wordpress.com/2020/08/26/a-trick-and-the-decomposition-theorem/}} (and which is credited there to Bhargav Bhatt, although we were later informed by Cadoret-Zheng that the argument essentially appears already around \cite[\S 5.3.11]{BBD}) in order to deduce the necessary results. 

\begin{definition}\label{defn: generic part}
Let $\cF$ be a perverse sheaf on a variety $Z$. Let $j \co U \inj Z$ be the inclusion of the maximal (dense) open variety on which $\cF$ is a (shifted) local system. We define the \emph{generic part of $\cF$} to be the perverse sheaf $\cF^{\gen} := j_{!*} j^* \cF$. 
\end{definition}

We make no claim that there is a non-trivial map between $\cF$ and $\cF^{\gen}$ in general. However, we shall prove that in certain situations of interest, $\cF^{\gen}$ can be realized as a direct summand of $\cF$.

\begin{proposition}\label{prop: generic summand}
Let $\cF$ be a perverse sheaf on a variety $Z$ over $L$. Suppose that $\cF|_{Z_{L^s}}$ is isomorphic to a direct sum of IC sheaves. Then $\cF^{\gen}$ is a direct summand of $\cF$. 
\end{proposition}

\begin{proof}
Let $j \co U \rightarrow Z$ be as in Definition \ref{defn: generic part} and $B$ be the closed complement of $U$. The hypothesis implies that $\cF|_{Z_{L^s}}$ is isomorphic to $\cF^{\gen}_{Z_{L^s}} \oplus \cE$ where $\cE$ is a direct sum of IC sheaves supported on $B_{L^s}$.

Consider the maps $\alpha \co (\pc j_!) j^* \cF \rightarrow  \cF$ and $\beta \co \cF \rightarrow (\pc j_*) j^* \cF$. 

Let $\cG := \Ima(\alpha) \subseteq \cF$. Then we have a tautological map $\cG \surj \cF^{\gen}$. Since $(\pc j_! j^* \cF )_{L^s}$ has no proper quotients supported on $B_{L^s}$ \cite[Exercise 3.1.6]{Ach}, the composition $\cG_{L^s}  \rightarrow \cF|_{Z_{L^s}} \rightarrow \cE$ vanishes, so that $\cG_{L^s} \inj \cF^{\gen}|_{L^s}$. This implies that the given map $\cG \rightarrow \cF^{\gen}$ is an isomorphism, which gives an injection $\cF^{\gen} \cong \cG \inj \cF$. A dual argument shows that for $\cG' := \Ima(\beta)$, we have a surjection $\cF \surj \cG' \cong \cF^{\gen}$. We then check after base change to $L^s$ that the composition $\cF^{\gen} \inj \cF \surj \cF^{\gen}$ is the identity map. 

\end{proof}

\begin{corollary}\label{cor: direct summand}
Let $f \co Y \rightarrow Z$ be a proper map of $L$-varieties, with $Y$ smooth over $L$. Then for all sufficiently large $\ell$, $\pch^0 (Rf_* \IC_Y(\Z_{\ell}))^{\gen}$ is a direct summand of $Rf_* \IC_Y(\Z_{\ell})$. 
\end{corollary}

\begin{proof}
By Theorem \ref{Decomposition Theorem}, the hypothesis of Proposition \ref{prop: generic summand} holds, so that $\pch^0 (Rf_* \IC_Y(\Z_{\ell}))^{\gen}$ is a direct summand of $\pch^0 (Rf_* \IC_Y(\Z_{\ell}))$. Applying Theorem \ref{Decomposition Theorem} again, $\pch^0 (Rf_* \IC_Y(\Z_{\ell}))$ is a direct summand of  $Rf_* \IC_Y(\Z_{\ell})$. 
\end{proof}

\begin{proof}[Proof of Proposition \ref{prop: IC direct summand integrally}]
By the definition of alteration, there is an open dense subset $U \subset Z$ over which the map $f|_U$ is finite flat. By shrinking $U$ further if necessary, we may assume that $U$ is contained in the smooth locus of $Z$, so that $\IC_Z(\Z_{\ell})$ is the intermediate extension of $\IC_Z(\Z_{\ell})|_U \cong \Z_{\ell}[\dim U]$. Over $U$, the composition $\Z_{\ell} \rightarrow f_* \Z_{\ell} \xrightarrow{\Tr} \Z_{\ell}$ is multiplication by $\deg f|_U$, hence as long as $\ell > \deg f|_U$ it realizes the constant sheaf $(\Z_{\ell})_U$ as a direct summand of $ f_* ((\Z_{\ell})_{f^{-1}(U)} ) \cong (Rf_* \Z_{\ell})|_U$. Since restriction to open subsets is perverse t-exact, this says that $\IC_U(\Z_{\ell})$ is a summand of $\pch^0(Rf_* \IC_Y(\Z_{\ell}))|_U$, and then applying intermediate extension shows that $\IC_Z(\Z_{\ell})$ is a direct summand of $\pch^0(Rf_* \IC_Y(\Z_{\ell}))^{\gen}$. We then conclude using Corollary \ref{cor: direct summand}. 
\end{proof}

\section{Global cyclic base change}

\subsection{Automorphic Galois representations}\label{ssec: automorphic Galois reps}

Let $R = \ol{\F}_{\ell}$ or $\ol{E}_{\lambda}$. In \cite[\S 2]{F20}, the \emph{abstract excursion algebra} $\Exc(\Gamma_{K}, \wh{G}_R)$ is defined. This is an algebra over $R$ such that the set of characters $\chi \co \Exc(\Gamma_{K}, \wh{G}_R) \rightarrow R$ is in bijection with semi-simple (not necessarily continuous) representations $\sigma \in H^1( \Gamma_{K}, \wh{G}(R))$. In \cite{X20, Xa}, Cong Xue extends Lafforgue's methods to define an action of $\Exc(\Gamma_{K}, \wh{G}_R) $ on $C_c^{\infty}(G(K) \setminus G(\A_K), R)$, preserving the cuspidal subspace. The same statement holds replacing $\Gamma_{K}$ by $\Gamma_{K,N}$ and $G(K) \setminus G(\A_K)$ by $G(K) \setminus G(\A_K) / U(N)$. The algebra $\cB(U(N),R)$ considered earlier in \S \ref{sec: excursion algebra} is the image of $\Exc(\Gamma_{K,N}, \wh{G}_R)$ in $\End_R(C_{\cusp}(U(N), R))$. 

 Following \cite[Definition 5.5]{F20}, we say that a Galois representation $\sigma \co \Gamma_{K} \rightarrow \wh{G}(R)$ is \emph{automorphic} if the corresponding maximal ideal $\mf{m}_{\sigma}$ appears in the support of the $\Exc(\Gamma_{K}, \wh{G}_R)$-action on $C_c(U(N), R)$ for some $N$. A priori this does not imply that $\sigma$ is related to a \emph{cusp} form, but the following Lemma shows that this is necessarily the case if $\sigma$ is absolutely irreducible. 

\begin{lemma}\label{lem: irreducible implies cuspidal}
Let $\sigma \co \Gamma_{K, S} \rightarrow \wh{G}(R)$ be an absolutely irreducible Galois representation which is moreover automorphic. Let $\mf{m}_{\sigma}$ the corresponding maximal ideal of $\Exc(\Gamma_K, \wh{G}_R)$. If $\sigma$ is irreducible, then $\mf{m}_{\sigma}$ appears in the support of the $\Exc(\Gamma_K, \wh{G}_R)$-action on $C_{\cusp}(U(N),R)$ for some $N$, and for this $N$ there is an eigenvector $f \in C_{\cusp}(U(N),R)$ on which $\cB$ acts through the character $\cB/\mf{m}_{\sigma}$. 
\end{lemma}

\begin{proof}
By \cite[Proposition 7.6.2]{Xa} and \cite[\S 6.2]{X20}, for a parabolic subgroup $P \subset G$ with Levi $M$, the constant term map 
\[
\CT_G^P \co   C_c(G(K) \setminus G(\A_K), R) \rightarrow C_c(M(K) \setminus M(\A_K), R)
\]
intertwines the action of $\Exc(\Gamma_{K}, \wh{G}_R)$ on the source with the action of $\Exc(\Gamma_{K}, \wh{M}_R)$ on the target, with respect to $\wh{M}(R) \inj \wh{G}(R)$. Therefore, all Galois representations attached by Lafforgue's construction to the image of $\CT_G^P$ factor through $\wh{M}(R)$. Hence if $\sigma$ does not factor through such a parabolic, and $C_c(G,R)_{\mf{m}_{\sigma}} \neq 0$, then we must have $C_{\cusp}(G,R)_{\mf{m}_{\sigma}} \neq 0$.

For the last statement, we note that since $C_{\cusp}(G,R) = \varinjlim_N C_{\cusp}(U,R)$, we have $(C_{\cusp}(U,R))_{\mf{m}_{\sigma}} \neq 0$ for some $U$. Since $C_{\cusp}(U,R)$ is a finite-dimensional vector space over $R$, it is the direct sum of its generalized eigenspaces for the $\Exc(\Gamma_{K}, \wh{G}_R)$-action, and the assumption implies that the character $\Exc(\Gamma_{K}, \wh{G}_R) \rightarrow \Exc(\Gamma_{K}, \wh{G}_R) / \mf{m}_{\sigma}$ appears as among the systems of eigenvalues. 
\end{proof}

\subsection{Mod $\ell$ base change}
We suppose $\ell$ is odd and is a {\it good prime} for the 
reductive group $G$. Explicitly, this means that we require $p>2$ if $\wh{G}$ has simple factors of type $A,B,C$ or $D$; $p>3$ if $\wh{G}$ has simple factors of type $G_2, F_4, E_6, E_7$; and $p>5$ if $\wh{G}$ has simple factors of type $E_8$. Let $K'/K$ be a cyclic extension of function fields of degree $\ell$. 

\begin{theorem}\cite[Theorem 1.6]{F20}\label{smith}  With notation as above, let $\bs:  \Gamma_{K} \ra \hG(\overline{\bbF}_\ell)$ be an automorphic Galois representation.  Then the restriction 
$\bs|_{\Gamma_{K'}}$ is automorphic for $G_{K'}$.
\end{theorem}



\subsection{Existence of almost all cyclic base changes in characteristic 0}

\begin{theorem}\label{thm: global existence}  
Let $\Pi$ be a cuspidal automorphic representation of $G(\bbA_K)$ defined over a number field $E$\footnote{
We emphasize again that every cuspidal automorphic representation of $G(\bbA_K)$ admits a model over the number field, so this is not a restriction.}. Suppose the Mumford-Tate group of $\Pi$ is $\hG$, as in Definition
\ref{MTgroup}. For a place $\lambda$ of $E$, let $\sigma_{\Pi,\lambda}$ be Lafforgue's parameter attached to $\Pi$ (as noted in Remark \ref{rem: weak correspondence}, the hypothesis implies that $\ul{\rho}_{\Pi, \lambda}$ is a singleton, so that $\sigma_{\Pi, \lambda}$ exists). 

Then there is an integer $c(\Pi)$, depending only on $\Pi$, such that for any prime $\ell$ satisfying $\ell > c(\Pi)$, and any cyclic extension $K'/K$ of degree $\ell$, there is an automorphic representation
$\Pi'$ of $G(\bbA_{K'})$ attached to the $L$-parameter $\sigma_{\Pi, \lambda}|_{\Gal((K')^s/K')}$. 
\end{theorem}

\begin{proof}Let $\bs_{\Pi, \lambda}$ denote the reduction of $\sigma_{\Pi,\lambda}$ modulo
$\lambda$; let $k_\lambda$ denote the residue field of $\lambda$.  Our hypothesis on the Mumford-Tate group of $\Pi$, together with the main theorem of \cite{BGP}, implies that the image of
$\bs_{\Pi, \lambda}(\Gamma_{K(\zeta_\ell)})$ contains $\hG(k_\lambda)^+$ for all $\lambda$ of sufficiently large characteristic.  It then follows from Lemma \ref{abundantimage} that there is a constant $b_1(\Pi)$ such that, for  $\ell > b_1(\Pi)$ the image of $\bs_{\Pi, \lambda}(\Gamma_{K'(\zeta_\ell)})$ is abundant, and that $\bs_{\Pi, \lambda}$ is also absolutely irreducible.

Now, it follows from Theorem \ref{almostall} that there is a constant $b_2(\Pi)$ such that for all $\ell > b_2(\Pi)$, the restriction of $\bs_{\Pi, \lambda}$ to $\Gamma_{K'}$ satisfies
condition (v) of Theorem \ref{thm_R_equals_B}. Taking $c (\Pi) := \max(b_1(\Pi), b_2(\Pi))$, we may then invoke Theorem \ref{smith} and Lemma \ref{lem: irreducible implies cuspidal} to apply Theorem \ref{thm_R_equals_B}. The freeneess guarantees the existence of an eigenform supported over a single point on the generic fiber of $R_{\overline{\sigma}_\ffrm, |N|, \cD}$, and we may take $\Pi'$ to be an irreducible cuspidal subquotient of the automorphic representation generated by it. 
\end{proof}

\section{Local cyclic base change}\label{sec: local base change}
Now let $F$ be a local function field. In this section we will prove a \emph{local} analogue of Theorem \ref{thm: global existence}, establishing existence of cyclic base change for irreducible smooth representations of $G(F)$, along any $\Z/\ell \Z$-extension of $F$ for almost all $\ell$. 

\subsection{The Genestier-Lafforgue correspondence}\label{ssec: GL}
Let $F$ be a local function field with ring of integers $\cO_F$ and residue characteristic $\ell \neq p$. Let $W_F$ be the Weil group of $F$. Let $G$ be a reductive group over $F$.  In \cite[Th\'{e}or\`{e}m 0.1]{Gen17}, Genestier-Lafforgue construct a map
\begin{align*}
\left\{ \begin{array}{@{}c@{}}  \text{irreducible admissible representations} \\  \text{$\pi$ of $G(F)$ over $\ol{\Q}_{\ell}$}\end{array} \right\}/\sim  & \longrightarrow \left\{ \begin{array}{@{}c@{}}  \text{semi-simple $L$-parameters} \\ 
\sigma_{\pi} \co W_F \rightarrow \wh{G}(\ol{\Q}_{\ell})  \end{array} \right\}/\sim.
\end{align*}
We shall use the two properties of the correspondence $\pi \mapsto \sigma_{\pi}$ recalled below. 

\subsubsection{Local-global compatibility}\label{sssec: LGC} For any automorphic representation $\Pi \cong \bigotimes_{x \in |X|} \Pi_x$ of $G(\A_K)$ that is associated to the $L$-parameter $\sigma_{\Pi}$ by V. Lafforgue's global Langlands parametrization (\S \ref{sec: excursion algebra}), $\sigma_{\Pi_x}$ is conjugate to the semi-simplification of $\sigma|_{W_{K_x}}$ for any $\sigma \in \ul{\rho}_{\Pi}$ (notation as in Remark \ref{rem: weak correspondence}).

\subsubsection{Compatibility with parabolic induction}\label{sssec: PI} If $P$ is a parabolic subgroup of $G$ with Levi quotient $M$, and $\tau$ is an irreducible admissible representation of $M(F)$, and $\pi$ is an irreducible subquotient of $\Ind_{P(F)}^{G(F)} \tau$ (with parabolic induction formed using the unitary normalization), then $\sigma_{\pi}$ is conjugate to the composition $W_F \xrightarrow{\sigma_{\tau}} \wh{M}(\ol{\Q}_{\ell}) \inj \wh{G}(\ol{\Q}_{\ell}) $.

\subsection{Existence of almost all cyclic base changes in characteristic 0}

\begin{definition}\label{def: base change}
Let $\pi$ be an irreducible admissible representation of $G(F)$ over $\ol{\Q}_{\ell}$. For a separable field extension $F'/F$, we say that an irreducible admissible representation $\pi'$ of $G(F')$ over $\ol{\Q}_{\ell}$ is a \emph{base change lifting of $\pi$ to $G(F')$} if $\sigma_{\pi'} \cong \sigma_{\pi}|_{W_{F'}}$. 
\end{definition}

This definition is an approximation to the notion of base change for $L$-packets. An $L$-packet for $G(F')$ should be said to be a base change lifting of an $L$-packet for $G(F)$ if the corresponding $L$-parameters are related by restriction on Weil groups. Since we lack a definition of $L$-packets for general groups and representations, we use the fibers of the Genestier-Lafforgue correspondence as a substitute for $L$-packets.

\begin{theorem}\label{thm: local existence}
Assume that $p$ is good for $G$ if $G$ is not simply laced. Let $\pi$ be an irreducible admissible representation of $G(F)$ over $\ol{\Q}$. There exists a constant $c(\pi)$ such that for all primes $\ell > c(\pi)$, for any $\Z/\ell \Z$-extension $F'/F$ there exists a base change lifting of $\pi$ to $G(F')$.
\end{theorem}

\begin{remark}In \cite[Theorem 1.1]{F20}, a version of this result was established for mod $\ell$ representations, when the extension is cyclic of degree equal to the same prime $\ell$, by completely different methods. 
\end{remark}


\begin{proposition}\label{prop: globalization}
Assume that $p$ is good for $G$ if $G$ is not simply laced. Let $\pi$ be an irreducible admissible supercuspidal representation of $G(F)$ over $\ol{\Q}$.  Then there exists a global field $K$ and a place $v \in |K|$ with $K_v \cong F$, and a cuspidal automorphic representation $\Pi$ of $G(\A_K)$ with Mumford-Tate group $\hG$ such that $\Pi_v \cong \pi$.  
\end{proposition}

\begin{proof}[Proof of Theorem \ref{thm: local existence} assuming Proposition \ref{prop: globalization}]
First we reduce to the case where $\pi$ is supercuspidal. Indeed, any $\pi$ can be realized as an irreducible subquotient of a parabolic induction of supercuspidal representations, of the form $\Ind_{P(F)}^{G(F)} \tau$. By \S \ref{sssec: PI}, the parabolic induction from $M(F')$ to $G(F')$ of a base change of $\tau$ to $M(F'	)$ will have an $L$-parameter of the desired form. So it suffices to treat the supercuspidal case. 

Hence we may and do assume for the rest of the argument that $\pi$ is supercuspidal. Then we may apply Proposition \ref{prop: globalization} to embed $\pi$ as the local component at $v$ of cuspidal automorphic representation $\Pi$ over a global field $K$ with $K_v \cong F$. Using Remark \ref{rem: weak correspondence}, we can replace $\Pi$ by an isomorphic $G(\A_K)$-representation (realized differently in the space of cuspidal functions on $G(K) \backslash G(\A_K) $) to assume that $\Pi$ is attached to an $L$-parameter $\sigma_{\Pi}$. Then Theorem \ref{thm: global existence} applies, so let $c(\pi) := c(\Pi)$ be as in Theorem \ref{thm: global existence}. For any $\Z/\ell\Z$-extension $F'/F$, we can find a $\Z/\ell \Z$-extension $K'/K$ with a place $v'$ lying over $v$ such that $K_{v'} \cong F'$. By Theorem \ref{thm: global existence}, there exists a cuspidal automorphic representation $\Pi'$ of $G(\A_{K'})$ with $L$-parameter $\sigma_{\Pi'} \cong \sigma_{\Pi}|_{\Gamma_{K'}}$. If $\Pi'_{v'}$ is the local component of $\Pi'$ at $v'$, then the local-global compatibility of \S \ref{sssec: LGC} ensures that $\sigma_{\Pi'_{v'}} \cong \sigma_{\Pi'}|_{W_{K'_{v'}}}$, hence  
\[
\sigma_{\Pi'_{v'}}  \cong \sigma_{\Pi'}|_{W_{K'_{v'}}}  \cong \sigma_{\Pi}|_{W_{F'}} \cong \sigma_{\pi}|_{W_{F'}}, 
 \]
 so $\Pi'_{v'}$ is the desired local base change. 
\end{proof}

\begin{remark}
Since all supercuspidal representations of $G(F)$ admits a model over $\ol{\Q}$, at least after twisting by a central character (which is unnecessary for us because our $G$ is semi-simple), Theorem \ref{thm: local existence} applies to all supercuspidal representations. In fact, we may present any irreducible admissible representation in terms of parabolic inductions and twists of representations defined over $\ol{\Q}$, and then the compatibility of the Genestier-Lafforgue correspondence with parabolic induction and twisting would allow to formulate a meaningful extension of Theorem \ref{thm: local existence} to all irreducible admissible representations. 
\end{remark}

\subsection{Globalization of supercuspidal representations} This subsection is devoted to the proof of Proposition \ref{prop: globalization}. We will use an argument due to Beuzart-Plessis, based on the Deligne-Kazhdan simple trace formula, to construct a globalization $\Pi$ of $\pi$ with specified local components at several auxiliary places that automatically force the image of Lafforgue's corresponding parameter to be Zariski dense.

\subsubsection{Genestier-Lafforgue parameters of simple supercuspidals} The notion of ``simple supercuspidal representations'' was singled out by Gross-Reeder \cite{Gro10}. Let us recall the definition.  

\begin{definition}
Let $G$ be a split semisimple group over a non-archimedean local field $F$ of residue characteristic $p \neq \ell$. A \emph{simple supercuspidal representation} is a representation of $G(F)$ that arises in the following way. Let $B \subset G$ be a Borel subgroup, with unipotent radical $U$. Let $I \subset G(\cO_F)$ be the corresponding Iwahori subgroup, and $I(1) \subset I$ its pro-unipotent radical. For an affine generic character $\phi \co I(1) \rightarrow \Q_{\ell}(\mu_p)$, we define $V_{\phi} := \cInd_{I(1) \times Z(G)}^{G(F)} (\phi \otimes \mathbbm{1})$. 
\end{definition}

Gross-Reeder showed that simple supercuspidal representations are irreducible and supercuspidal. They anticipated the shape of the associated $L$-parameters. To explain this, let $\cJ \subset W_{F}$ be the inertia subgroup. Let us call \emph{simple supercuspidal parameter} a discrete\footnote{Recall that a parameter $\sigma \co \cW \times \wh{G}(\ol{\Q}_{\ell})$ is called \emph{discrete} if its image is not contained in a proper parabolic subgroup of $\wh{G}$.} Langlands parameter $\sigma \co W_F \rightarrow \wh{G}(\ol{\Q}_{\ell})$ such that the adjoint representation $\Ad \circ \sigma$ has $\wh{\mf{g}}^{\cJ}  =0$ and Swan conductor equal to the rank of $\wh{G}$. We emphasize that a discrete Langlands parameter, a fortiori a simple supercuspidal parameter, is already semi-simple; this will be used when invoking the local-global compatibility of \S \ref{sssec: LGC}. 



\begin{proposition}\label{prop: simple wild parameter}
Assume that $p$ is good for $G$ if $G$ is not simply laced. Then the Genestier-Lafforgue parameter of $V_{\phi}$ is a simple supercuspidal parameter.
\end{proposition}

\begin{proof}
The idea is to find a globalization of $V_{\phi}$ to a cuspidal automorphic representation on $\mathbf{P}^1$, whose associated global $L$-parameter can be computed explicitly. A particularly convenient such globalization was studied by Heinloth-Ng\^{o}-Yun \cite{HNY} and we only need to collect the relevant consequences of their work. 

Let $K = \F_q(\PP^1)=\F_q(t)$. Let $I_0$ be the Iwahori subgroup of $G(\F_q((t)))$, and $I(1)_{\infty}$ be the pro-unipotent radical of the Iwahori subgroup of $G(\F_q((1/t)))$. As a special case of \cite[Proposition 2.7]{Yun16}, there is a unique automorphic representation $\Pi \cong \otimes'_{x \in |\mathbf{P}^1|} \Pi_x$ of $G(\A_K)$ such that:
\begin{itemize}
\item $\Pi_x$ is unramified if $x \neq 0, \infty$. 
\item $\Pi_0$ has a vector invariant under $I_0$. 
\item $\Pi_{\infty}$ has an eigenvector under $I(1)_{\infty}$, on which the action is given by $\phi$. 
\end{itemize}
Moreover, this $\Pi$ is cuspidal and appears with multiplicity one in the automorphic spectrum, and 
\[
\dim C_c^{\infty} \left( G(K) \backslash G(\A_K) / \prod_{x \neq 0, \infty} G(\cO_x), \ol{\Q} \right)^{I_0 \times (I(1)_{\infty}, \phi)} = 1.
\]
Since the excursion algebra acts on this 1-dimensional $\ol{\Q}_{\ell}$-vector space, its action is automatically through a character, which determines the $L$-parameter $\sigma_{\Pi}$ associated to $\Pi$ by Lafforgue's correspondence. 

In this case, Heinloth-Ng\^{o}-Yun construct a Hecke eigensheaf $A_{\phi}$ a moduli stack $\Bun_{G(0,2)}$ of $G$-bundles on $\PP^1$ with level structure at $0$ and $\infty$, whose associated Frobenius trace function is a non-zero $f_{\phi} \in C_c^{\infty} \left( G(K) \backslash G(\A_K) / \prod_{x \neq 0, \infty} G(\cO_x), \ol{\Q} \right)^{I_0 \times (I(1)_{\infty}, \phi)}$. They further prove that $A_{\phi}$ is a Hecke eigensheaf, with corresponding local system the (generalized) Kloosterman local system $\mrm{Kl}_{\wh{G}}(\phi) \co \pi_1(\mathbb{P}^1 - \{0, \infty\}, \wh{G}(\ol{\Q}_{\ell}))$. By \cite[Theorem 2, Corollary 2.15]{HNY}, letting $\cJ_{\infty} \subset W_{K_{\infty}}$ be the inertia group at $\infty$, the local monodromy representation $ \mrm{Kl}_{\wh{G}}(\phi)|_{\cJ_{\infty}} \co \cJ_{\infty} \rightarrow \wh{G}(\ol{\Q}_{\ell})$ at $\infty$ is irreducible and it is the restriction to $\cJ_\infty$ of a simple supercuspidal parameter.

By \cite[Proposition 6.4]{BHKT}, the parameter associated to $f_{\phi}$ by Lafforgue's correspondence must coincide with $\mrm{Kl}_{\wh{G}}(\phi)$. By local-global compatibility, the Genestier-Lafforgue parameter of $V_{\phi}$ therefore agrees upon restriction to $\cJ_{\infty}$ with $\mrm{Kl}_{\wh{G}}(\phi)|_{\cJ_{\infty}}$. 
\end{proof}

\subsubsection{Globalization}

The following Lemma is proved by Beuzart-Plessis in the appendix to \cite{GHS}, using the Deligne-Kazhdan simple trace formula:

\begin{lemma}\label{lem: globalize with local prescription}
Let $\pi$ be a supercuspidal representation of $G(F)$ over $\ol{\Q}$. There exists a global curve $X$, with function field $K = \F_q(X)$ and places $v,w,w', w'' \in |X|$, and a cuspidal representation $\Pi \cong \bigotimes'_{x \in |X|} \Pi_x$ of $G(\A_K)$ such that:
\begin{enumerate}
\item $K_v \cong F$ and $\Pi_v \cong \pi$ as $G(F)$-representations,
\item $\Pi_w$ is a simple supercuspidal representation of $G(K_w)$,
\item $\Pi_{w'}$ has non-zero trace against the pseudo-coefficient for the Steinberg representation given by the Euler-Poincar\'{e} function \cite[Theorem 8.2.1]{Lau96}, and 
\item $\Pi_{w''}$ is an unramified representation parametrized by a regular element of $\wh{T}//W$.  Moreover, the parameter can be chosen to avoid any finite
union of proper subtori of $\wh{T}$. 
\end{enumerate}
\end{lemma}

We will need a characterization of irreducible subgroups of $\wh{G}$ containing a principal unipotent element.  The following result is probably known but
we were unable to find a reference.

\begin{lemma}\label{Dynkin}  Let $\frg$ be a semisimple Lie algebra and $\frh$ a proper semisimple subalgebra containing a regular nilpotent element of $\frg$.
Then the rank of $\frh$ is strictly less than that of $\frg$.
\end{lemma}

\begin{proof}  Let $E$ be the regular nilpotent of $\frg$ which lies in $\frh$. We first claim that $E$ is also regular in $\frh$. We can complete $E$ to an $\mathfrak{sl}_2$-triple 
$(E, H, F)$ in $\frh$. The element $H$ is regular semisimple in $\frg$, hence in $\frh$, which implies by \cite[VIII. 11.4, Proposition 7]{Bou05} that $E$ is regular nilpotent in $\frh$.

Let $r(\frg)$ and $\ell(\frg)$ denote respectively the ranks of $\frg$ and its length as a representation of the $3$-dimensional subalgebra spanned by $(E,H,F)$; define $r(\frh)$ and 
$\ell(\frh)$ analogously.  The  Proposition just cited from \cite{Bou05} implies that
$$r(\frg) = \ell(\frg); ~~ r(\frh) = \ell(\frh).$$
Thus $r(\frg) = r(\frh)$ if and only if $\frg$ and its subspace $\frh$ have the same length as $(E,H,F)$-modules; but this is only possible if $\frg = \frh$.
\end{proof}

\begin{corollary}\label{finitesl2}  Let $\cH$ denote the set of conjugacy classes of proper reductive subgroups of $\wh{G}$ containing a given principal unipotent element, and for each $H \in \cH$,
let $\wh{T}_H \subset \wh{T}//W$ denote the conjugacy class of its maximal torus.  Then the complement of the union of the $\wh{T}_H$ is Zariski dense in $\wh{T}//W$.
\end{corollary}

\begin{proof} A classification of $\cH$ in all types appears in Saxl-Seitz \cite[Theorem A, Theorem B]{SS97} (in characteristic 0, the result has been credited to earlier work of Dynkin). It is finite, it follows from Lemma \ref{Dynkin} that each $\wh{T}_H$ is of codimension at least $1$ in $\wh{T}//W$ (this is also clear from inspection of $\cH$). Therefore, $\bigcup_{H \in \cH}$ is a finite union of positive-codimension subvarieties in $\wt{T}//W$. 
\end{proof}


\subsubsection{Determination of global monodromy}
Let $\Pi$ be as in Lemma \ref{lem: globalize with local prescription}. Let $\ul{\rho}_{\Pi}$ be the set of Remark \ref{rem: weak correspondence} and $\sigma \in \ul{\rho}_{\Pi}$. We claim that $\sigma$ has Zariski-dense image. This will complete the proof of Proposition \ref{prop: globalization}.

So it only remains to prove the claim. By local-global compatibility, the semi-simplification of $\sigma|_{K_x}$ corresponds to $\Pi_x$ under the Genestier-Lafforgue correspondence for all $x \in |X|$. By Proposition \ref{prop: simple wild parameter} and local-global compatibility, the semi-simplification of $\sigma_{\Pi}|_{K_w}$ is already absolutely irreducible, so $\sigma_{\Pi}$ is absolutely irreducible. By \cite[Lemma 11.4]{ST21}, for any representation $V$ of $\wh{G}$ the corresponding local system is $V \circ \sigma_{\Pi}$ is pure of weight zero. By condition (iii) and compatibility of the Genestier-Lafforgue correspondence with parabolic induction, $\sigma|_{W_{K_{w'}}}$ has the same semi-simplification as the Steinberg representation.  It follows as in \cite[\S 4.3]{HNY} from purity of the weight-monodromy filtration \cite[Theorem I.8.4]{Del80} that the image of $\sigma|_{K_{w'}}$ contains a principal unipotent element.

By \cite[Corollaire I.3.9]{Del80} plus the property that $V \circ \sigma_{\Pi}$ is pure of weight zero for every $V$, the neutral component of the Zariski closure of the image of $\sigma$ in $\wh{G}$ is semi-simple. Furthermore, we have just seen that it does not lie any proper parabolic subgroup, and also that it contains a principal unipotent element. By Lemma \ref{Dynkin}, 
proper subgroups of $\wh{G}$ have smaller rank. 
But this is ruled out by (iv), where the proper subtori are those in the statement of Corollary \ref{finitesl2}. \qed 

\section{Cyclic base change for toral supercuspidal representations}\label{sec: toral supercuspidal}

In this section we will investigate cyclic base change more explicitly for a certain class of representations singled out in \cite[\S 6]{Kal19}, called ``toral supercuspidal representations''. The strategy is to first explicate cyclic base change for the mod $\ell$ reductions of these representations, using a Conjecture of Treumann-Venkatesh (established in \cite{F20} for the Genestier-Lafforgue correspondence) that ``base change functoriality is realized by Tate cohomology''. Then, using the deformation theory of Galois representations, we will lift the result to characteristic zero.  

The main new work is in calculating the Tate cohomology of toral supercuspidal representations, and what facilitates this calculation is a geometric model for these representations established by Chan-Oi \cite{CO21}, as compact inductions from parahoric subgroups of ``generalized Deligne-Lusztig inductions'' studied by Chan-Ivanov. The generalized Deligne-Lusztig representations are produced from the cohomology of ``Deligne-Lusztig type varieties'' built out of group schemes coming from the Moy-Prasad filtration, analogously to the way in which Deligne-Lusztig varieties are built from reductive groups over finite fields. Hence our computation naturally breaks into two steps: (1) studying of Tate cohomology of Deligne-Lusztig type varieties, which we do in \S \ref{ssec: tate of CI} and (2) studying Tate cohomology of compact inductions, which we do in \S \ref{ssec: compact inductions}. 

Kaletha constructs an explicit Local Langlands parametrization of the toral supercuspidal representations in \cite{Kal19}. On the other hand, we identify base change relations among the Genestier-Lafforgue parameters of toral supercuspidal representations. This gives some evidence for consistency between the Genestier-Lafforgue correspondence \cite{Gen17} and Kaletha's construction of $L$-packets of regular supercuspidal representations. 

We emphasize that in this section, $F$ is a non-archimedean local field having residue characteristic $p$, but we allow $F$ to have characteristic zero in all results up through Theorem \ref{thm: tate coh of Chan-Oi}.

\subsection{Work of Chan-Ivanov} We briefly recall the generalized Deligne-Lusztig representations appearing in \cite{CI19}. 

Let $G$ be a reductive group over a non-archimedean local field $F$. Let $T \inj G$ be an unramified maximal torus and $x \in \cB(G/F)$ be a point of the Bruhat-Tits building of $G$ that lies in the apartment of $T$. If $F'/F$ is a tamely ramified extension, then we may also regard $x$ as a point of $\cB(G_{F'}/F')$ using the identification $\cB(G/F) = \cB(G_{F'}/F')^{\Gal(F'/F)}$.  Corresponding to $x$ we have by Bruhat-Tits theory a parahoric group scheme $\cG/\cO_F$, whose generic fiber is $G/F$. 

Let $\F_q$ be the residue field of $F$. By assumption, $T$ splits over $\breve{F} := W(\ol{\F}_q)$. Let $U$ be the unipotent radical of an $\breve{F}$-rational Borel subgroup of $G_{\breve{F}}$ containing $T_{\breve{F}}$. 

For $r \in \Z_{\geq 0}$, we have group schemes $\GG_r, \TT_r, \UU_r$ over $\F_q$ as in \cite[\S 2.5, 2.6]{CI19} corresponding to subquotients of the Moy-Prasad filtration at $x$, such that 
\[
\GG_r(\F_q) = G_{x, 0:r+} := G_{x,0}/G_{x,r+}, \quad \TT_r(\F_q) = T_{0:r+} := T_{x, 0}/T_{x,r+}
\]
and $\UU_r \subset (\GG_r)_{\ol{\F}_q}$. 

\subsubsection{Deep level Deligne-Lusztig varieties}
We recall certain schemes constructed in \cite[\S 4]{CI19}, generalizing Deligne-Lusztig varieties. Let 
\[
S_{\TT_r,\UU_r} := \{ x \in \GG_r \co x^{-1} \Fr_q(x) \in \UU_r\}.
\]
(The variety $S_{\TT_r,\UU_r}$ is called $X_r$ in \cite{CO21}.) It is a separated, smooth, finite type scheme over $\ol{\F}_q$, with an action of $G_{x,0:r+} \times T_{x,0:r+}$ by multiplication on the left and right, and also a free action of $\UU_r \cap \Fr_q^{-1} (\UU_r)$ by right translation. 

It is actually more convenient for us to work with 
\[
Y_{\TT_r,\UU_r} := S_{\TT_r,\UU_r}/\UU_r \cap \Fr_q^{-1} (\UU_r).
\]
Since the natural map $S_{\TT_r,\UU_r} \rightarrow Y_{\TT_r,\UU_r}$ is a bundle in affine spaces, the compactly supported (geometric) cohomology groups of source and target are identified, equivariantly for the $G_{x,0:r+} \times T_{x,0:r+}$-action, up to Tate twist and an (even) shift of cohomological degrees, so they will lead to the same (virtual) representations. The $Y_{\TT_r, \UU_r}$ are called ``deep level Deligne-Lusztig varieties''. 

\begin{example}When $r=0$, the definition of $Y_{\TT_r,\UU_r}$ specializes to that of a classical Deligne-Lusztig variety. 
\end{example}

\begin{definition}[Generalized Deligne-Lusztig induction]\label{def: DL induction} Let $\Lambda \in \{\ol{\F}_{\ell}, \ol{\Q}_{\ell}, k, \cO\}$ be an $\ell$-adic coefficient ring and $\theta \co \TT_r(\F_q) \rightarrow \Lambda^{\times}$ be a character. We denote 
\[
H_c^*(Y_{\TT_r,\UU_r}; \Lambda)_{\theta} := H_c^*(Y_{\TT_r,\UU_r}; \Lambda)_{\theta}  \otimes_{\Lambda[T_{x, 0:r^+}]} \theta,
\]
which is a graded representation of $G_{x, 0:r^+}$.  We define the virtual representation of $G_{x, 0:r^+}$,
\[
R^{\GG_r}_{\TT_r,\UU_r }(\theta) := \sum_i (-1)^i [H_c^i(Y_{\TT_r,\UU_r}; \Lambda)_{\theta}] \in K_0(G_{x,0:r+}; \Lambda). 
\]
The version of this definition with $\Lambda = \ol{\Q}_{\ell}$ is considered in \cite[Definition 4.4]{CI19}, while we will also be interested in $\Lambda = \ol{\F}_{\ell}$.
\end{definition}

\begin{remark}
By inflation, we may view $\theta$ as a character of $T_{x,0}$ that is trivial on $T_{x,r+}$, and $R^{\GG_r}_{\TT_r,\UU_r}(\theta) \in K_0(G_{x, 0}; \Lambda)$. In practice, $\theta$ will come by restriction from a character of $T(F)$. 
\end{remark}


\begin{definition}
Following \cite[\S 2.10]{CI19} in the case of $\ol{\Q}_{\ell}$-coefficients, we say that $\theta \co \TT_r(\F_q)  \rightarrow k^{\times}$ is \emph{regular} if for any absolute root $\alpha$ and any $d \geq 1$ such that $\Fr_q^d(\alpha) = \alpha$, the restriction of 
\[
\TT_r(\F_{q^d}) \xrightarrow{\Nm} \TT_r(\F_q) \xrightarrow{\theta}  \ol{\F}_{\ell}^{\times}
\]
to the subgroup (defined in \cite[\S 2.6]{CI19}) $\TT_r^{\alpha, r} \subset \TT_r(\F_{q^d})$ is non-trivial. 
\end{definition}

For a finite group $\Gamma$ and a virtual representation $V \in K_0(\Gamma; \Lambda)$ which is a representation up to sign, we define $|V|$ to be the underlying representation and $(-1)^V$ to be the sign of $V$, so that $V = (-1)^V |V|$. 

\begin{lemma}\label{lem: independence} 
Let $\Lambda \in \{ \ol{\Q}_{\ell}, \ol{\F}_{\ell}\}$. 

(i) If $\theta$  is regular, then $R^{\GG_r}_{\TT_r,\UU_r}(\theta) \in K_0(\GG_r(\F_q); \Lambda)$ is independent of $\UU_r$. 

(ii) If $\theta$ is regular and the stabilizer of $\theta$ in the Weyl group of the special fiber of $\cG$ is trivial, then $\pm R_{\TT_r,\UU_r}^{\theta}$ is a (non-virtual) representation of $G_{x,0:r+}$ and $|R_{\TT_r,\UU_r}^{\theta}|$ is irreducible. 
\end{lemma}

\begin{proof}
Both (i) and (ii) are established in \cite{CI19} when $\Lambda = \ol{\Q}_{\ell}$, so we shall simply explain the reduction from this case to $\Lambda = \ol{\F}_{\ell}$. 

(i) Any $\theta \co T_{x,0:r+} \rightarrow  \ol{\F}_{\ell}^{\times}$ lifts canonically to $\wt{\theta} \co T_{x,0:r+}  \rightarrow \ol{\Z}_{\ell}^{\times} \subset \ol{\Q}_{\ell}^{\times}$ via the Teichm\"{u}ller map. Moreover, $R^{\GG_r}_{\TT_r,\UU_r}(\theta)$ is the image of $R^{\GG_r}_{\TT_r,\UU_r}(\wt{\theta})$, the characteristic $0$ Deligne-Lusztig induction studied in \cite[\S 4.1]{CI19}, under the reduction map $K_0(G_{x,0:r+} ; \ol{\Q}_{\ell}) \rightarrow K_0(G_{x,0:r+};  \ol{\F}_{\ell})$. If $\theta$ is regular then $\wt{\theta}$ is regular in the sense of \cite[\S 2.10]{CI19}, so by \cite[Corollary 4.7]{CI19} $R^{\GG_r}_{\TT_r,\UU_r}(\wt{\theta})$ is independent of $U$.  

(ii) This follows from lifting to characteristic zero and a similar argument as in (i), using \cite[Corollary 4.7]{CI19} for the analogous statement in characteristic zero. We note that $|R_{\TT_r,\UU_r}^{\theta}|$ is non-zero by \cite[Corollary 4.6]{CI19}.
\end{proof}


\subsection{Tate cohomology of some generalized Deligne-Lusztig inductions}\label{ssec: tate of CI}
\subsubsection{Recollections on Tate cohomology}

Let $\sigma$ be an order $\ell$ endomorphism of an abelian group $V$. Write $N := 1 + \sigma + \ldots + \sigma^{\ell-1} \in \Z[\langle \sigma \rangle]$. The \emph{Tate cohomology} groups of $V$ (with respect to the $\sigma$-action) are defined as:
\begin{align*}
\rT^0 (\sigma, V) = \rT^0(V) := \frac{\ker( 1- \sigma \co V \rightarrow V)}{N \cdot V},\\
\rT^1 (\sigma, V)= \rT^1(V) := \frac{\ker(N \co V \rightarrow V) }{(1-\sigma) \cdot V}.
\end{align*}
It is sometimes convenient to extend the definition of $T^i(V)$ to all $i \in \Z$, so that $T^i (-)= T^{i+2}(-)$. 

Given a short exact sequence 
\[
0 \rightarrow V' \rightarrow V \rightarrow V'' \rightarrow 0,
\]
there is a (periodic) long exact sequence on Tate cohomology 
\begin{equation}\label{eq: LES}
\ldots \rightarrow \rT^0(V') \rightarrow \rT^0(V) \rightarrow \rT^0(V'') \rightarrow \rT^1(V') \rightarrow \rT^1(V) \rightarrow \rT^1(V'') \rightarrow \rT^2(V') \rightarrow \ldots 
\end{equation}

In \cite[\S 3.4]{F20}, we defined the notion of Tate cohomology for a scheme $Y$ with an admissible action of $\Z/\ell \Z \cong \langle \sigma \rangle$. Admissibility automatically holds if $Y$ is a quasiprojective variety over a field, and we will always be in this situation when invoking this theory, so let us assume $Y$ is such. Another useful description of the (compactly supported) \emph{Tate cohomology group $\rT^i(\sigma, Y;  \Lambda) = \rT^i(Y;  \Lambda)$ with coefficients in $\Lambda$} is the $i$th cohomology of the totalization of 
\begin{equation}\label{eq: tate double complex}
\begin{tikzcd}
\ldots  & \deg -1 & \deg 0 & \deg 1  & \ldots \\
\ldots \ar[r, "1-\sigma"] & R\Gamma_c(Y;  \Lambda) \ar[r, "N" ]  & R\Gamma_c(Y;  \Lambda) \ar[r, "1-\sigma" ] & R\Gamma_c(Y;  \Lambda) \ar[r, "N"] & \ldots 
\end{tikzcd}
\end{equation}
It is immediate from this definition that $\rT^i(-) \cong \rT^{i+2}(-)$. We shall be interested in coefficients such as $\Lambda \in \{\ol{\F}_{\ell}, \ol{\Z}_{\ell},  k, \cO,  \cO/\ell^n\}$.

The double complex \eqref{eq: tate double complex} leads to two spectral sequences abutting to $\rT^*(Y; \ol{\F}_{\ell})$. 
\begin{itemize}
\item The ``vertical then horizontal'' spectral sequence has first page 
\[
E_1^{ij} = H^j_c(Y; \Lambda).
\]
Its $j$th row is the complex computing $\rT^i(H^j_c(Y; \Lambda))$, namely
\[
\ldots \xrightarrow{N} H^j_c(Y; \Lambda) \xrightarrow{1-\sigma} H^j_c(Y; \Lambda) \xrightarrow{N} H^j_c(Y; \Lambda) \xrightarrow{1-\sigma} \ldots 
\]
Therefore, the second page is 
\[
E_2^{ij} = \rT^i(\sigma, H^j_c(Y;  \Lambda)) .
\]
\item The ``horizontal then vertical'' spectral sequence has
\[
E_2^{ij} = H^j_c(Y^{\sigma};  \Lambda)
\]
and is moreover degenerate starting from $E_2$ \cite[Theorem 4.4]{TV16}. 
\end{itemize}

\subsubsection{Tate cohomology of representations}\label{ssec: Tate cohomology}

Let $V$ be a representation over $k$ of a finite group $\Gamma$. Suppose $\Z/\ell \Z \cong \langle \sigma \rangle$ acts on $\Gamma$. Then $V \mapsto V \circ \sigma$ defines an action of $\sigma$ on isomorphism classes of representations of $\Gamma$ over $\ol{\F}_{\ell}$.

\begin{lemma}
If $V$ is irreducible and $V \cong V \circ \sigma$ as $\Gamma$-representations, then there is a unique extension of the $\Gamma$-action on $V$ to an action of $\Gamma \rtimes \langle \sigma \rangle$. 
\end{lemma}

\begin{proof}
If $A \co V \xrightarrow{\sim} V \circ \sigma$ as $\Gamma$-representations, then Schur's Lemma implies that the composition $A^\ell \co  V \xrightarrow{\sim} V \circ \sigma^\ell = V$ is multiplication by a scalar, say $\lambda$. Then defining $\sigma$ to act as $\lambda^{-1/\ell}  A$ gives an extension of the desired form. Schur's Lemma also implies that $A$ is unique up to scalar, and scaling $A$ evidently results in the same extension. 
\end{proof}

For a $\Z/\ell\Z \cong \langle \sigma \rangle $-module $V$, we have defined Tate cohomology groups $\rT^0(V)$ and $\rT^1(V)$. If the $\sigma$-action on $V$ extends to an action of $\Gamma \rtimes \langle \sigma \rangle$, then $\rT^0(V)$ and $\rT^1(V)$ inherit an action of the subgroup of $\Gamma$ fixed by $\sigma$, which we denote $Z_\Gamma(\sigma)$.

\subsubsection{Torsion in integral cohomology}\label{sssec: torsion}

We write 
\[
H_c^i(Y_{\TT_r, \UU_r}; \ol{\Z}_{\ell}) \cong 
H_c^i(Y_{\TT_r, \UU_r}; \ol{\Z}_{\ell})_{\mrm{tf}}  \oplus 
H_c^i(Y_{\TT_r, \UU_r}; \ol{\Z}_{\ell})_{\tors}
\]
for the decomposition into the torsion-free and torsion summands, respectively. For a character $\wt{\theta} \co \TT_r(\F_q) \rightarrow \ol{\Z}_{\ell}^{\times}$, we write 
\[
H_c^i(Y_{\TT_r, \UU_r}; \ol{\Z}_{\ell})_{\mrm{tf}, \wt{\theta}} := H_c^i(Y_{\TT_r, \UU_r}; \ol{\Z}_{\ell})_{\mrm{tf}} \otimes_{\ol{\Z}_{\ell}[\TT_r(\F_q)]}  \wt{\theta}.
\]
Thus $H_c^i(Y_{\TT_r, \UU_r}; \ol{\Z}_{\ell})_{\mrm{tf}, \wt{\theta}} $ is a lattice in $H_c^i(Y_{\TT_r, \UU_r}; \ol{\Q}_{\ell})_{\wt{\theta}}$. 

Write $\theta := \wt{\theta} \otimes_{\ol{\Z}_{\ell}} \ol{\F}_{\ell} \co \TT_r(\F_q) \rightarrow \ol{\F}_{\ell}^{\times}$. We consider the hypothesis, 
\begin{equation}\label{eq: one degree}
\text{$H_c^*(Y_{\TT_r, \UU_r}; \ol{\F}_{\ell})_{\theta}$ is non-zero in only one degree.}
\end{equation}
We will then be interested in studying $\rT^i(\sigma, H_c^*(Y_{\TT_r, \UU_r}; \ol{\F}_{\ell})_{\theta})$ as a representation of $G_{x, 0:r+}$. For regular $\theta$, it is expected that $H_c^*(Y_{\TT_r, \UU_r}; \ol{\Q}_{\ell})_{\theta}$ is non-zero in only one degree. Therefore, if $\theta$ is regular and $q$ is sufficiently large, then one would expect \eqref{eq: one degree} to hold for all large enough $\ell$, as a case of the more general expectation that for any (finite type) variety over a separably closed field, the \'{e}tale cohomology with coefficients in $\Z_{\ell}$ should be torsion-free for all sufficiently large $\ell$; however, this does not appear to be known in general. When $r=0$, in which case the $Y_{\TT_r, \UU_r}$ are usual Deligne-Lusztig varieties, this concentration of degree and torsion-freeness for $\theta$ in general position is proved in \cite[Lemma 3.5]{Br90}, which already furnishes many interesting examples where \eqref{eq: one degree} is known. Based on this, it seems natural to guess: 

\begin{conjecture}\label{conj: 0-toral}
If $\theta$ is a 0-toral character \cite[Definition 3.7]{CO21}, then \eqref{eq: one degree} holds. 
\end{conjecture}

Later we will prove theorems concerning 0-toral $\theta$ under the assumption that \eqref{eq: one degree} holds. It would be interesting to make these unconditional by proving Conjecture \ref{conj: 0-toral}. As discussed, we already have many unconditional depth zero examples. 

\begin{lemma}\label{lem: lattice}
Suppose $\theta \co \TT_r(\F_q) \rightarrow \ol{\F}_{\ell}^{\times}$ satisfies \eqref{eq: one degree} and let $\wt{\theta} \co \TT_r(\F_q) \rightarrow \ol{\Z}_{\ell}^{\times}$ be the composition of $\theta$ with the Teichm\"{u}ller lift. Let $i$ be the unique degree in which $H_c^i(Y_{\TT_r, \UU_r}; \ol{\F}_{\ell})_{\theta}$ is non-zero. Then $H_c^i(Y_{\TT_r, \UU_r}; \ol{\Z}_{\ell})_{\wt{\theta}}$ is a lattice in $H_c^i(Y_{\TT_r, \UU_r}; \ol{\Q}_{\ell})_{\wt{\theta}}$.
\end{lemma}

\begin{proof}
We must show that $H_c^i(Y_{\TT_r, \UU_r}; \ol{\Z}_{\ell})_{\wt{\theta}}$ is torsion-free. The universal coefficient theorem implies that $\left(H_c^i(Y_{\TT_r, \UU_r}; \ol{\Z}_{\ell})_{\wt{\theta}}\right)_{\tors} \otimes_{\ol{\Z}_{\ell}} \ol{\F}_{\ell}$ injects into $H_c^{i+1}(Y_{\TT_r, \UU_r}; \ol{\F}_{\ell})_{\theta} = 0$, so $\left(H_c^i(Y_{\TT_r, \UU_r}; \ol{\Z}_{\ell})_{\wt{\theta}}\right)_{\tors} = 0$ by Nakayama's Lemma. 
\end{proof}

\subsubsection{Base change}\label{sssec: base change}

Recall that $F_{\ell}$ is the degree $\ell$ unramified extension of $F$. Let $T' = \Res_{F_{\ell}/F} (T_{F_\ell})$, $G' = \Res_{F_{\ell}/F}(G_{F_\ell})$. We may view $x$ as a point in the apartment of $T'$, in the building of $G'/F$. 

Given $U \subset G_{\breve{F}}$, we have several possible choices for $U' \subset G'_{\breve{F}}$. We will make a particular choice that is advantageous for the upcoming computation. A choice of generator $\sigma \in \Gal(F_{\ell}/F)$ induces an isomorphism 
\[
F_{\ell} \otimes_F \breve{F}  \cong \underbrace{\breve{F} \times \ldots \times \breve{F} }_{\ell \text{ times}} =: \breve{F}^{\ell}
\]
sending $x \otimes y  \mapsto (xy, \sigma(x)y, \ldots, \sigma^{\ell-1}(x) y)$. This in turn induces 
\begin{equation}\label{eq: weil restriction}
G'_{\breve{F}} \cong G_{\breve{F}} \times \ldots \times G_{\breve{F}}.
\end{equation}
The action of $\sigma$ on the LHS transports to the cyclic permutation on the RHS. We define $U' \subset G'_{\breve{F}}$ to be the image of $U \times \ldots \times U \subset G_{\breve{F}}^{\ell}$ under the isomorphism \eqref{eq: weil restriction}. Then $U'$ is stable under the $\sigma$-action on $G'_{\breve{F}}$. The resulting deep level Deligne-Lusztig variety $Y_{\TT_r', \UU_r'}$ therefore carries an action of $\sigma$, such that for $g \in G'_{x, 0:r+}  = \GG_r'(\F_q)$ and $y \in Y_{\TT_r', \UU_r'}(R)$, 
\[
\sigma \cdot (g \cdot y) = \sigma(g) \cdot \sigma(y)
\]
and similarly for $t \in T'_{x,0:r+} = \TT_r'(\F_q)$. This induces a $\sigma$-action on $H_c^i(Y_{\TT_r', \UU_r'}; \ol{\F}_{\ell})$ that is compatible in the same manner with the $\GG'_r(\F_q) \times \TT_r'(\F_q)$-action. We will consider Tate cohomology with respect to this action.

Recall that for a representation $V / \ol{\F}_{\ell}$ of a finite group $\Gamma$, the \emph{Frobenius twist} of $V$ is the representation
\[
V^{(\ell)} := V \otimes_{\ol{\F}_{\ell}, \Frob_{\ell}} \ol{\F}_{\ell}.
\]
Now note that $H_c^i(Y_{\TT_r, \UU_r}; \ol{\F}_{\ell})$ has an $\F_{\ell}$-structure induced by cohomology with coefficients in $\F_{\ell}$, and therefore has an action of $\Aut(\ol{\F}_{\ell})$. Let $\theta \co T_{x, 0:r+} = \TT_r(\F_q) \rightarrow \ol{\F}_{\ell}^{\times}$ be a character and $\theta' \co \TT_r'(\F_q) \rightarrow \ol{\F}_{\ell}^{\times}$ be the composition of $\theta$ with $\Nm \co \TT_r'(\F_q)  \rightarrow \TT_r(\F_q) $. Applying $\Frob_{\ell} \in \Aut(\ol{\F}_{\ell})$ induces an isomorphism of $\GG_r(\F_q)$-representations
\begin{equation}\label{eq: Frob twist of DL}
(H_c^i(Y_{\TT_r, \UU_r}; \ol{\F}_{\ell}) \otimes_{\ol{\F}_{\ell}[\TT_r(\F_q)]} \theta)^{(\ell)} \xrightarrow{\sim} H_c^i(Y_{\TT_r, \UU_r}; \ol{\F}_{\ell}) \otimes_{\ol{\F}_{\ell}[\TT_r(\F_q)]} \theta^{\otimes \ell}.
 \end{equation}

\begin{theorem}\label{thm: DL tate cohomology}
Let $\theta \co \TT_r(\F_q) \rightarrow \ol{\F}_{\ell}^{\times}$. Assume that $\ell \nmid \# \TT_0(\F_q)$. 

(i) If $H^*_c(Y_{\TT_r, \UU_r};\ol{\F}_{\ell})_{\theta} \neq 0$, then $\rT^i(\sigma, H_c^{i'}(Y_{\TT_r', \UU_r'};\ol{\F}_{\ell})_{\theta'}) \neq 0$ for some $i \in \{0,1\}$ and some $i'$.

(ii) Suppose that $\theta$ and $\theta'$ both satisfy \eqref{eq: one degree}, and are non-zero in degrees $j,j'$ respectively. Then for $i \in \{0,1\}$ we have 
\[
\rT^i (\sigma, H_c^{j'}(Y_{\TT_r', \UU_r'}; \ol{\F}_{\ell})_{\theta'}  \otimes_{\ol{\Z}_{\ell}} \ol{\F}_{\ell} ) \cong  (H^j_c(Y_{\TT_r, \UU_r, F};\ol{\F}_{\ell})_{\theta})^{(\ell)}
\]
as representations of $\GG_r(\F_q)$. 
\end{theorem}

\begin{proof} We begin with some analysis of the relevant varieties. For a scheme $Y$ over $\F_q$ we will write $\ol{Y}$ for its base change to $\ol{\F}_q$. 

Parallel to \eqref{eq: weil restriction}, the choice of generator $\sigma \in \Gal(\F_{q^\ell}/\F_q)$ identifies
\begin{equation}\label{eq: weil restriction base changed}
\ol{\GG}'_r \xrightarrow{\sim} \ol{\GG}_r \times \ldots \times \ol{\GG}_r =: \ol{\GG}_r^\ell.
\end{equation}
The isomorphism \eqref{eq: weil restriction base changed} transports:
\begin{itemize}
\item The automorphism $\sigma$ on the LHS to the cyclic rotation action on the RHS.
\item The action of $\Fr_q$ on the LHS (coming from the $\F_q$-rational structure $\GG'_r$ for $\ol{\GG}'_r$) to the composition of the cyclic rotation with $(\Fr_q, \ldots , \Fr_q)$ on the RHS, where each factor of $\Fr_q$ comes from the $\F_q$-rational structure $\GG_r$ on $\ol{\GG}_r$. Let us denote this endomorphism of $\ol{\GG}_r^{\ell}$ by $F := \sigma \circ (\Fr_q, \ldots , \Fr_q)$. 
\end{itemize}
Similar remarks apply to $\TT'_r \cong \Res_{\F_{q^\ell}/\F_q} (\TT_{r, \F_{q^{\ell}}})$.

Under \eqref{eq: weil restriction base changed}, the variety $Y_{\TT_r', \UU_r'} $ has the presentation
\[
Y_{\TT_r', \UU_r'}  \cong \{ (g_1, \ldots, g_{\ell}) \in \ol{\GG}_r^{\ell}  \co (g_1, \ldots, g_{\ell})^{-1} F( g_1 \ldots, g_{\ell}) \in  \UU_r^{\ell} \} / \UU_r^{\ell} \cap F^{-1} \UU_r^{\ell}.
\]
The action of $\sigma$ on $Y_{\TT_r', \UU_r'} $ transports to the cyclic rotation on factors in the presentation on the RHS. 	

\begin{lemma}\label{lem: fixed points of CI}
The diagonal map $Y_{\TT_r, \UU_r} \rightarrow Y_{\TT_r', \UU_r'}$ identifies with the inclusion of the $\sigma$-fixed points on $Y_{\TT_r', \UU_r'}$, equivariantly for the action of $\GG_r(\F_q) \times \TT_r(\F_q)$ by left and right translation. 
\end{lemma}

\begin{proof}
Immediate upon writing down the definitions.
\end{proof}

\begin{corollary}\label{cor: equiv localization}
For every $i$, restriction induces an isomorphism
\[
\rT^i(Y_{\TT_r' , \UU_r'}; \ol{\F}_{\ell}) \cong \rT^i(Y_{\TT_r, \UU_r}; \ol{\F}_{\ell})
\]
which is equivariant with respect to the action of $Z_{\GG_r'(\F_q) \times \TT_r'(\F_q)}(\sigma)  = \GG_r(\F_q) \times \TT_r(\F_q)$, induced by left and right translation on $Y_{\TT_r' , \UU_r'}$ and $Y_{\TT_r, \UU_r}$. 
\end{corollary}

\begin{proof} This follows from Lemma \ref{lem: fixed points of CI} and the equivariant localization theorem of \cite[\S 3.4.5]{F20}.
\end{proof}

With these preparations, we are now ready for the proofs of parts (i) and (ii) of the Theorem. 

\emph{Proof of (i).} Restricting $\theta'$ along $\TT_r(\F_q) \inj \TT_r'(\F_q)$ yields $\theta'|_{\TT_r(\F_q)} = \theta^{\otimes \ell}$, so Corollary \ref{cor: equiv localization} induces an isomorphism of $\GG_r(\F_q)$-representations
\begin{equation}
\rT^*(Y_{\TT_r', \UU_r'}; \ol{\F}_{\ell}) \otimes_{\ol{\F}_{\ell}[\TT_r(\F_q)]} (\theta'|_{\TT_r(\F_q)})  \cong  \rT^*(Y_{\TT_r, \UU_r}; \ol{\F}_{\ell}) \otimes_{\ol{\F}_{\ell}[\TT_r(\F_q)]} \theta^{\otimes \ell}.
\end{equation}
Similarly to \eqref{eq: Frob twist of DL}, we have 
\[
(\rT^*(Y_{\TT_r, \UU_r}; \ol{\F}_{\ell}) \otimes_{\ol{\F}_{\ell}[\TT_r(\F_q)]} \theta)^{(\ell)} \xrightarrow{\sim} \rT^*(Y_{\TT_r, \UU_r}; \ol{\F}_{\ell}) \otimes_{\ol{\F}_{\ell}[\TT_r(\F_q)]} \theta^{\otimes \ell}.
 \]
 Note that Frobenius twisting shows that $H^*_c(Y_{\TT_r, \UU_r};\ol{\F}_{\ell})_{\theta} \neq 0 \iff H^*_c(Y_{\TT_r, \UU_r};\ol{\F}_{\ell})_{\theta^{\otimes \ell}} \neq 0$. 

As explained in \S \ref{ssec: Tate cohomology}, there are two spectral sequences abutting to $\rT^i(Y_{\TT'_r, \UU_r'};  \ol{\F}_{\ell})$. One degenerates at $E_2$ and has $E_2^{ij} = H^j_c((Y_{\TT'_r, \UU_r'})^{\sigma}; \ol{\F}_{\ell})$, which by Lemma \ref{lem: fixed points of CI} is $\GG_r(\F_q) \times \TT_r(\F_q)$-equivariantly isomorphic to $H^j_c(Y_{\TT_r, \UU_r};  \ol{\F}_{\ell})$, and the other has $E_2^{ij} = \rT^i (\sigma, H^j_c(Y_{\TT'_r, \UU_r'}; \ol{\F}_{\ell}))$.

Since $\ker (\TT_r(\F_q) \rightarrow \TT_0(\F_q))$ is pro-$p$, the assumption that $\ell \nmid \TT_0(\F_q)$ implies that $\ell \nmid \TT_r(\F_q)$. Therefore, $- \otimes_{ \ol{\F}_{\ell}[\TT_r(\F_q)]} \theta$ implements projection to a summand. By the preceding paragraph, the Jordan-H\"{o}lder factors (as $\GG_r(\F_q)$-representations) of $H^j_c(Y_{\TT_r, \UU_r}; \ol{\F}_{\ell})_{\theta^{\otimes \ell}}$ are also subquotients of some $\rT^i (H^{j'}_c(Y_{\TT'_r, \UU_r'}; \ol{\F}_{\ell})) \otimes_{ \ol{\F}_{\ell}[\TT_r(\F_q)]} \theta^{\otimes \ell}$. Therefore the proof of (i) will be concluded by the following Lemma, which allows us to  ``commute'' the formation of Tate cohomology with the projection to an isotypic component. 

\begin{lemma}\label{lem: exchange Tate with char}
We have that 
\[
\rT^i (H^{j'}_c(Y_{\TT'_r, \UU_r'}; \ol{\F}_{\ell})) \otimes_{ \ol{\F}_{\ell}[\TT_r(\F_q)]} \theta^{\otimes \ell} \cong \rT^i (H^{j'}_c(Y_{\TT'_r, \UU_r'}; \ol{\F}_{\ell})_{\theta'})
\]
as $\GG_r(\F_q)$-representations.
\end{lemma}

\begin{proof}
Since $\# \TT_r(\F_q)$ is coprime to $\ell$, we have that 
\[
\rT^i (H^{j'}_c(Y_{\TT'_r, \UU_r'}; \ol{\F}_{\ell})) \otimes_{ \ol{\F}_{\ell}[\TT_r(\F_q)]} \theta^{\otimes \ell} \cong \rT^i (H^{j'}_c(Y_{\TT'_r, \UU_r'}; \ol{\F}_{\ell})\otimes_{ \ol{\F}_{\ell}[\TT_r(\F_q)]} \theta^{\otimes \ell} ).
\]
as $\GG_r(\F_q)$-representations. Here on the RHS we have projected to a $\TT_r(\F_q)$-isotypic component before forming Tate cohomology, and we need to show that the same answer is computed if we instead project to a particular $\TT_r'(\F_q)$-isotypic component. 

Let $\{ \vartheta_l \}$ be the set of characters of $\TT_r'(\F_q)$ that extend $\theta^{\otimes \ell}$ on $\TT_r(\F_q) \subset \TT_r'(\F_q)$, indexed so that $\vartheta_1 = \theta'$. Then 
\begin{equation}\label{eq: sublemma 1}
\rT^i (H^{j'}_c(Y_{\TT'_r, \UU_r'}; \ol{\F}_{\ell})\otimes_{ \ol{\F}_{\ell}[\TT_r(\F_q)]} \theta^{\otimes \ell} ) \cong \rT^i \left( \bigoplus_l H^{j'}_c(Y_{\TT'_r, \UU_r'}; \ol{\F}_{\ell})\otimes_{ \ol{\F}_{\ell}[\TT_r'(\F_q)]} \vartheta_l 	\right). 
\end{equation}

We claim that $\theta'$ is the only $\sigma$-equivariant extension of $\theta^{\otimes \ell}$ to $\TT_r'(\F_q)$. Indeed, any $\sigma$-equivariant character on $\TT_r'(\F_q)$ factors through the norm map $\TT_r'(\F_q) \xrightarrow{\Nm} \TT_r(\F_q)$, and the composition $\TT_r(\F_q) \inj \TT_r'(\F_q) \xrightarrow{\Nm} \TT_r(\F_q)$ is multiplication by $\ell$. Since $\ol{\F}_{\ell}^{\times}$ is $\ell$-torsion-free, there is only one character $\TT_r(\F_q)  \rightarrow \ol{\F}_{\ell}^{\times}$ that inflates to $\theta^{\otimes \ell}$ under multiplication by $\ell$, namely $\theta$. 

Therefore, $\sigma$ permutes the set $\{\vartheta_l \co  l>1\}$ without any fixed points, necessarily grouping them into free orbits, and it therefore freely permutes the summands of
\begin{equation}\label{eq: free part}
\bigoplus_{l>1} H^{j'}_c(Y_{\TT'_r, \UU_r'}; \ol{\F}_{\ell})\otimes_{ \ol{\F}_{\ell}[\TT_r'(\F_q)]} \vartheta_l.
\end{equation}
Since Tate cohomology of a free $\sigma$-module vanishes, 
\[
\rT^i (
\bigoplus_{l>1} H^{j'}_c(Y_{\TT'_r, \UU_r'}; \ol{\F}_{\ell})\otimes_{ \ol{\F}_{\ell}[\TT_r'(\F_q)]} \vartheta_l) = 0.
\]
Therefore, \eqref{eq: sublemma 1} is $\GG_r(\F_q)$-equivariantly isomorphic to $ \rT^i (H^{j'}_c(Y_{\TT'_r, \UU_r'}; \ol{\F}_{\ell})_{\theta'})$. 
\end{proof}

\emph{Proof of (ii).} We consider the Tate spectral sequence, as in part (i). The additional assumptions imply that
\begin{itemize}
\item $R\Gamma_c(Y_{\TT_r', \UU_r'}; \ol{\F}_{\ell})_{\theta'}$ is concentrated in a single degree $j'$, and so  is quasi-isomorphic to $H_c^{j'}(Y_{\TT_r', \UU_r'}; \ol{\F}_{\ell})_{\theta'}$. 
\item $R\Gamma_c(Y_{\TT_r, \UU_r}; \ol{\F}_{\ell})_{\theta}$ is concentrated in a single degree $j$, and so is quasi-isomorphic to $H_c^{j}(Y_{\TT_r, \UU_r}; \ol{\F}_{\ell})_{\theta}$.
\end{itemize}
Corollary \ref{cor: equiv localization} then implies that 
\begin{equation}\label{eq: thmii 1}
\rT^i(\sigma, H_c^{j'}(Y_{\TT_r', \UU_r'}; \ol{\F}_{\ell})) \cong  H_c^{j}(Y_{\TT_r, \UU_r}; \ol{\F}_{\ell})
\end{equation}
as $\GG_r(\F_q) \times \TT_r(\F_q)$-representations. Using Lemma \ref{lem: exchange Tate with char} we find that 
\begin{equation}\label{eq: thmii 2}
\rT^i (\sigma, H^{j'}_c(Y_{\TT_r', \UU_r'}; \ol{\F}_{\ell})_{\theta '} ) \cong \rT^{i+j'}(Y_{\TT_r', \UU'_r}; \ol{\F}_{\ell}) \otimes_{\ol{\F}_{\ell}[\TT_r(\F_q)]} \theta^{\otimes \ell}  
\end{equation}
as $\GG_r(\F_q)$-representations. Projecting \eqref{eq: thmii 1} to the $\theta^{\otimes \ell}$-isotypic component and  then using \eqref{eq: Frob twist of DL} to relate the Frobenius twist of the $\theta$-isotypic component with the $\theta^{\otimes \ell}$-isotypic component, the proof is concluded. 
\end{proof}

\subsection{Tate cohomology of compact inductions}\label{ssec: compact inductions} We study the relationship between compact induction and Tate cohomology. In this section, we let $\msf{G}'$ be any $p$-adic group with an action of $\Z/\ell \Z \cong \langle \sigma \rangle$, $\msf{H}' \subset \msf{G}'$ a $\sigma$-invariant closed subgroup, $\msf{G} = (\msf{G}')^{\sigma}$ and $\msf{H} = (\msf{H}')^{\sigma}$. 

\begin{proposition}\label{prop: c-ind}
Let $\pi$ be a finite-dimensional representation of $\msf{H}$. If $\msf{G}/\msf{H} \xrightarrow{\sim} (\msf{G}'/\msf{H}')^{\sigma}$, then
\[
\rT^i(\cInd_{\msf{H}'}^{\msf{G}'} \pi) \cong \cInd_{\msf{H}}^{\msf{G}}(\rT^i \pi)
\]
as $\msf{G}'$-representations. 
\end{proposition}

\begin{proof}[Proof of Proposition \ref{prop: c-ind}] 
A special case appears in \cite[Proposition 14]{R16}, which already contains the main ideas of the proof. We use Bernstein-Zelevinsky's perspective of $l$-sheaves on $l$-spaces. There is an equivalence of categories between $\msf{G}'$-equivariant sheaves on $\msf{G}'/\msf{H}'$ and representations of $\msf{H}'$, which we denote $\cF_{\pi} \leftrightarrow \pi$. Furthermore, the Tate cohomology of $\cF_{\pi}$ as a sheaf transports to the Tate cohomology of $\pi$ as a representation.

Under this equivalence and the analogous equivalence between $\msf{G}'$-equivariant sheaves on a point and $\msf{G}'$-representations, the $\msf{G}'$-representation $\cInd_{\msf{H}'}^{\msf{G}'}(\pi)$ corresponds to the $\msf{G}'$-equivariant sheaf $\pr_! \cF_{\pi}$ on $\pt$, where $\pr \co \msf{G}' / \msf{H}' \rightarrow \pt$, under the functor of taking global sections for $\msf{G}'$-equivariant sheaves on $\pt$. By \cite[\S 3.3]{TV16}, the restriction map on sections induces
\[
\rT^i(\cInd_{\msf{H}'}^{\msf{G}'} \pi) \xrightarrow{\sim} \cInd_{\msf{H}}^{\msf{G}}(\rT^i \pi),
\]
which completes the proof. 
\end{proof}

Next we work out some situations where the hypothesis of Proposition \ref{prop: c-ind} is satisfied. We have the long exact sequence 
\[
0 \rightarrow \msf{H} \rightarrow \msf{G} \rightarrow \msf{G}/\msf{H} \rightarrow H^1(\langle \sigma \rangle, \msf{H}') \rightarrow \ldots 
\]
so $\msf{G}/\msf{H} \xrightarrow{\sim} (\msf{G}'/\msf{H}')^{\sigma}$ if $H^1(\langle \sigma \rangle, \msf{H}') = 0$.

\begin{lemma}\label{lem: tate coh connected reductive}
Let $F_{\ell}/F$ be the unramified $\Z/\ell\Z$-extension of local fields of characteristic $\ell \neq p$. Let $H$ be a connected algebraic group over $\cO_F$ and $\msf{H}' = H(\cO_{F_{\ell}})$, with the action of $\Gal(F_{\ell}/F) \cong \langle \sigma \rangle$ by Galois conjugation on points. Then $H^1(\sigma, \msf{H}')= 0$. 
\end{lemma}

\begin{proof}
Let $\msf{H}^+$ be the kernel of the reduction map $\msf{H}' = H(\cO_{F_{\ell}}) \rightarrow H(\F_{q^\ell})$. From the LES of cohomology, we have the exact sequence
\[
\ldots \rightarrow H^1(\sigma, \msf{H}'_+) \rightarrow H^1(\sigma, \msf{H}') \rightarrow H^1(\sigma, H(\F_{q^\ell})) \rightarrow \ldots.
\] 
Since $\msf{H}'_+$ is pro-$p$, we have $H^1(\sigma, \msf{H}'_+) = 0$. By Lang's Theorem, $H^1(\sigma, H(\F_{q^\ell})) = 0$. Therefore, $H^1(\sigma, \msf{H}')  = 0$. 
\end{proof}

We shall be particularly interested in the case where $\msf{G}' = G(F_{\ell})$ and $\msf{H'} = Z(\msf{G}') G(F_{\ell})_{x,0}$ for $x \in \cB(G/F_{\ell})$ fixed by $\Gal(F_{\ell}/F)$, which guarantees that $\msf{H'}$ is stable under $\Gal(F_{\ell}/F)$. By unramified descent in Bruhat-Tits theory, $G(F_{\ell})_{x,0}^{\Gal(F_{\ell}/F)} = G(F)_{x,0}$, so we have $\msf{G} := (\msf{G}')^{\Gal(F_{\ell}/F)} = G(F)$ and $\msf{H} := (\msf{H}')^{\Gal(F_{\ell}/F)}  = Z(\msf{G}) G(F)_{x,0}$. 

Let $Z \subset G$ be the maximal central torus. Then $Z_{F_{\ell}}$ is the maximal central torus of $G_{F_{\ell}}$.

\begin{lemma}\label{lem: fixed points}
Assume that 
\begin{enumerate}
\item[(i)] The action of $\Gal(F_{\ell}/F)$ on the cocharacter group $X^*(Z_{F_{\ell}})$ is trivial. 
\item[(ii)] The component group $Z(\msf{G}')/Z$ has order coprime to $\ell$. 
\end{enumerate}
Let $H' = Z(\msf{G}')G(F_\ell)_{x,0}$. Then $H^1(\sigma, \msf{H}') = 0$. 
\end{lemma}

\begin{example}
The assumptions (i) is satisfied, for example, if $G$ is split reductive, or whenever $G$ is (not necessarily split and) semi-simple. For any given $G$, assumption (ii) is satisfied for all large enough $\ell$. 
\end{example}

\begin{proof}
Since $G$ is assumed to be split over an unramified extension of $F$, $Z$ is an unramified torus, hence has a canonical integral model $\cZ/\cO_F$ (the ``connected N\'{e}ron model''), which has the property that $\cZ(\cO_F)$ is the maximal bounded subgroup of $Z(F)$. We have a split short exact sequence
\[
0 \rightarrow \cZ(\cO_{F_{\ell}}) \rightarrow Z(F_{\ell}) \rightarrow Z(F_{\ell}) / \cZ(\cO_{F_{\ell}}) \rightarrow 0
\]
where $ Z(F_\ell) / \cZ(\cO_{F_{\ell}}) \cong X_*(Z_{F_{\ell}})$ is torsion-free, and $\cZ(\cO_{F_{\ell}})  = Z(F_{\ell}) \cap G(F_{\ell})_{x,0} \subset G(F_{\ell})_{x,0}$. Therefore, $Z(F_{\ell}) G(F_{\ell})_{x,0} \cong X_*(Z_{F_{\ell}}) \times  G(F_{\ell})_{x,0}$ so 
\[
H^1(\sigma, Z(F_{\ell}) G(F_{\ell})_{x,0})  \cong H^1(\sigma, X_*(Z_{F_{\ell}})) \times H^1(\sigma, G(F_{\ell})_{x,0}).
\]
Now, $H^1(\sigma, X_*(Z_{F_{\ell}}))  = 0$ because the assumed condition (i) implies $H^1(\sigma, X_*(Z_{F_{\ell}})) \cong \Hom(\langle \sigma \rangle, X_*(Z_{F_{\ell}})) = 0$. Since $G(F_{\ell})_{x,0}$ is the group of $\cO_{F_{\ell}}$-points of a connected Bruhat-Tits group scheme, the group $H^1(\sigma, G(F_{\ell})_{x,0})$ vanishes by Lemma \ref{lem: tate coh connected reductive}. Therefore, $H^1(\sigma, Z(F_{\ell}) G(F_{\ell})_{x,0}) = 0$. 

Finally, assumption (ii) implies that the index of $Z(F_{\ell}) G(F_{\ell})_{x,0}$ is normal subgroup of $Z(\msf{G}')G(F_{\ell})_{x,0}$ with finite index coprime to $\ell$, so the long exact sequence implies that $H^1(\sigma, Z(\msf{G}')G(F_{\ell})_{x,0})  = 0$. 
\end{proof}

\subsection{Toral supercuspidal representations}\label{ssec: toral supercuspidal base change} 
We will prove base change results for a class of representations studied by Chan-Oi in \cite{CO21}, which they call \emph{toral supercuspidal} representations.

\subsubsection{Assumptions}\label{sssec: assumptions} We impose the same assumptions as in \cite[\S 7]{CO21}. In particular, $G$ is a reductive group over $F$, $p$ is odd and not bad for $G$, $p \nmid \pi_1(G_{\mrm{der}})$ and $p \nmid \pi_1(\wh{G}_{\mrm{der}})$. The maximal torus $T \subset G$ is unramified elliptic, and $\wt{\theta} \co T(F) \rightarrow \ol{\Z}_{\ell}^{\times}$ is ``$0$-toral'' (what was also previously called ``regular'') and of depth $r >0$, i.e. trivial on $T_{x,r+}$ for some $x \in \cB(G/F)$. We write $\theta \co T \rightarrow k^{\times}$ for the reduction of $\wt{\theta}$ modulo $\ell$. We assume that $R_{\TT_r, \UU_r}^{\GG_r}(\wt{\theta})$ is non-zero, which is automatic under a mild regularity hypothesis (see proof of Lemma \ref{lem: independence}(ii)). 

\subsubsection{Work of Chan-Oi} The assumptions imply that $R_{\TT_r, \UU_r}^{\GG_r}(\wt{\theta}) \otimes_{\ol{\Z}_{\ell}} \ol{\Q}_{\ell}$ is irreducible and independent of the choice of $\UU_r$, hence we simply abbreviate it as $R_{\TT_r}^{\GG_r}(\wt{\theta})_{\ol{\Q}_{\ell}}$. Furthermore, by \cite[Theorem 7.2]{CO21} the $\ol{\Q}_{\ell}$-representation $\cInd_{T(F)G(F)_{x,0}}^G( | R_{\TT_r}^{\GG_r}(\wt{\theta})_{\ol{\Q}_{\ell}}|)$ is irreducible supercuspidal. 

Given $\wt{\theta}$, there is defined in \cite[Definition 4.8]{CO21} a certain sign character $\epsilon^{\mrm{ram}}[\wt{\theta}] \co T(F) \rightarrow \ol{\Q}_{\ell}^{\times}$. We define 
\[
\pi_{T, \wt{\theta} \cdot \epsilon^{\mrm{ram}}[\wt{\theta}] } := \cInd_{T(F)G(F)_{x,0}}^{G(F)} | R_{\TT_r, \UU_r}(\wt{\theta} )_{\ol{\Q}_{\ell}}|.
\]
The twisting by $\epsilon^{\mrm{ram}}[\wt{\theta}]$ is for consistency with Kaletha's indexing of regular supercuspidal representations \cite{Kal19}; 
$\pi_{T, \wt{\theta} \cdot \epsilon^{\mrm{ram}}[\wt{\theta}] } $ is a regular supercuspidal representation by \cite[Theorem 7.2]{CO21}. In particular, it is irreducible.

Since $\epsilon^{\mrm{ram}}[\wt{\theta}]$ takes values in $\{\pm 1\}$, we may regard it also a character valued in $\ol{\F}_{\ell}^{\times}$ by reduction modulo $\ell$. We may define
\[
\pi_{T, \theta \cdot \epsilon^{\mrm{ram}}[\theta]} :=  \cInd_{T(F)G(F)_{x,0}}^{G(F)}  R_{\TT_r, \UU_r}(\theta ) \in K_0(G(F); \ol{\F}_{\ell}).
\]
A priori this is only a class in $K_0(G(F); \ol{\F}_{\ell})$, but there are two circumstances in which we can lift it to an honest representation, which will be denoted $|\pi_{T, \theta \cdot \epsilon^{\mrm{ram}}[\theta]}|$. 
\begin{itemize}
\item If $\ell$ is banal, in which case we may lift $\pm R_{\TT_r, \UU_r}(\theta )$ to a representation $|R_{\TT_r, \UU_r}(\theta )|$. 
\item If $H_c^{*}(Y_{\TT_r, \UU_r}; E)_{\wt{\theta}}$ concentrates in a single degree, in which case we define 
\[
|\pi_{T, \theta \cdot \epsilon^{\mrm{ram}}[\theta]}| := \cInd_{T(F)G(F)_{x,0}}^{G(F)} ( H_c^{*}(Y_{\TT_r, \UU_r}; \cO)_{\mrm{tf}, \wt{\theta}} \otimes_{\cO} k ).
\]
In this case $|\pi_{T, \theta \cdot \epsilon^{\mrm{ram}}[\theta]}|$ is the mod $\ell$ reduction of a $\ol{\Z}_{\ell}$-lattice in $\pi_{T, \wt{\theta} \cdot \epsilon^{\mrm{ram}}[\wt{\theta}] } $, namely the one induced by the lattice $H_c^{*}(Y_{\TT_r, \UU_r}; \cO)_{\mrm{tf}, \wt{\theta}} \subset  |R_{\TT_r, \UU_r}(\wt{\theta})|$. Hence, if $\ell$ is banal, then this definition coincides with the in the previous bullet point, so that $|\pi_{T, \theta \cdot \epsilon^{\mrm{ram}}[\theta]}|$ is unambiguously defined. Furthermore, if \eqref{eq: one degree} is satisfied, then by Lemma \ref{lem: lattice} we have 
\[
|\pi_{T, \theta \cdot \epsilon^{\mrm{ram}}[\theta]}| \cong \cInd_{T(F)G(F)_{x,0}}^{G(F)} ( H_c^{*}(Y_{\TT_r, \UU_r}; k)_{\theta}).
\]
\end{itemize}

\begin{proposition}\label{modellirrcusp} If $\ell$ is banal for $G$, then (under our running hypotheses) the representation $\pi_{T, \theta \cdot \epsilon^{\mrm{ram}}[\theta]}$ is irreducible and cuspidal.
\end{proposition}

\begin{proof}
Indeed, \cite[Theorem 7.2]{CO21} identifies 
$\pi_{T, \wt{\theta} \cdot \epsilon^{\mrm{ram}}[\wt{\theta}] } $ with the compact induction of Yu's representation, denoted $^\circ \tau_d$ in {\it loc. cit.}.  Since $\ell$ is banal
the mod $\ell$ reduction of $^\circ \tau_d$ is again irreducible and obtained from Yu's construction. Its compact induction is identified with  $\pi_{T, \theta \cdot \epsilon^{\mrm{ram}}[\theta]}$. The result then follows from \cite[Theorem 6.1]{Fi}.
\end{proof}

\begin{remark}  Marie-France Vign\'eras has indicated another proof of Proposition \ref{modellirrcusp}.  Let 
$$\tau = | R_{\TT_r}^{\GG_r}(\wt{\theta})_{\ol{\Q}_{\ell}}|$$
and let $[\tau]$ denote the reduction mod $\ell$ of any $\ol{\Z}_{\ell}$-lattice in $\tau$; this is independent of the choice because $\ell$ is banal.
Write $U = T(F)G(F)_{x,0}$.
Lemma 3.2 of her article \cite{Vig01} proves a simple
criterion for irreducibility of the compact induction of $[\tau]$, which we adapt to our present notation:
\begin{itemize}
\item[(a)] $\End_{\ol{\F}_{\ell}[G]}(\cInd_{U}^G[\tau]) = \ol{\F}_{\ell}$;
\item[(b)]  Let $\pi$ be any irreducible $\ol{\F}_{\ell}$-representation of $G$.  If $[\tau]$ is contained in the restriction of $\pi$ to $U$
then $[\tau]$ is also a quotient of $\pi|_U$.  
\end{itemize}
Now (b) is automatic because $\ell$ is banal and therefore prime to the pro-order of $U$.    On the other hand,
it follows from the supercuspidality of $\cInd_U^G \tau$ that for any $g \in G \setminus U$ the restriction of $\tau$ and $g(\tau)$ to the
intersection $U \cap g Ug^{-1}$ are disjoint.  Again, since $\ell$ is banal, the same holds for the restriction of
$[\tau]$ and $g([\tau])$.   Point (a) then follows, and this implies irreducibility.  Cuspidality is then a consequence of \cite[Theorem II.2.7]{Vig96},
because the matrix coefficients of the compactly induced representation are compact modulo center.
\end{remark}


\subsubsection{Base change for toral supercuspidal representations}

Let $G$ be a reductive group over $F$ and $F'/F$ a $\Z/\ell\Z$-extension. Let $\pi'$ be an irreducible (admissible) representation of $G(F')$ over $\ol{\F}_{\ell}$. Choose a generator $\sigma$ of $\Gal(F'/F) \cong \Z/\ell \Z$. Then $\sigma$ acts on $G(F')$ through its Galois action on $E$. We say that $\pi'$ is \emph{$\sigma$-fixed} if $\pi' \cong \pi' \circ \sigma$ as $G(F')$-representations. 

\begin{lemma}[{\cite[Proposition 6.1]{TV16}}]\label{lem: sigma-fixed}
If $\pi'$ is $\sigma$-fixed, then the action of $G(F')$ extends uniquely to a $G(F') \rtimes \langle \sigma \rangle$-action.
\end{lemma}

The \emph{Tate cohomology} groups of $\pi'$, with respect to the $\sigma$-action, are 
\[
\rT^0(\pi') := \frac{\ker( 1- \sigma \co \pi' \rightarrow \pi')}{(1+ \sigma + \ldots + \sigma^{\ell-1}) \cdot \pi'} , \quad 
\rT^1(\pi') := \frac{\ker(1+ \sigma + \ldots + \sigma^{\ell-1} \co \pi' \rightarrow \pi') }{(1-\sigma) \cdot \pi'}.
\]
The $G(F')$-action on $\pi'$ induces an action of $G(F)$ on $\rT^i(\pi')$.

Let $\pi$ be an irreducible admissible representation of $G(F)$ over $\ol{\F}_{\ell}$, and $\pi'$ be an irreducible admissible representation of $G(F')$ over $\ol{\F}_{\ell}$. Recall that in Definition \ref{def: base change} we defined what it means for $\pi'$ to be a \emph{base change} of $\pi$. In this situation we say that $\pi$ is a \emph{base change descent} of $\pi'$. Addressing \cite[Conjecture 6.5]{TV16}, it was proved in \cite[Theorem 1.3]{F20} that if $F$ has characteristic $p\neq \ell$, and $\ell$ is odd and good for $\wh{G}$, then any irreducible $G(F)$-subquotient $\pi$ of $\rT^i(\pi')$ base changes to $(\pi')^{(\ell)}$, the Frobenius twist of $\pi'$. The Theorem below computes Tate cohomology of mod $\ell$ toral supercuspidal representations for $F' = F_{\ell}$ the unramified extension of order $\ell$.

\begin{theorem}\label{thm: tate coh of Chan-Oi}
Let $G,T, \theta$ be as in \S \ref{sssec: assumptions}. Let $G', T', \theta'$ be as in \S \ref{sssec: base change}. Assume that $T(F_{\ell})$ is elliptic. 

(i) Assume that $\ell \nmid \# \TT_0(\F_q)$ and $G$ satisfies the assumptions of Lemma \ref{lem: fixed points}. If $H_c^*(Y_{\TT_r, \UU_r}, \ol{\F}_{\ell})_{\theta} \neq 0$, then there exists $i \in \{0,1\}$, $i'$ such that 
\[
\rT^i \left( \cInd_{T(F_{\ell})G(F_{\ell})_{x,0}}^{G(F_{\ell})}  H_c^{i'}(Y_{\TT_r' ,\UU_r'}; \ol{\F}_{\ell} )_{\theta'}  \right)\neq 0.
\]

(ii) In addition to the assumptions from (i), suppose further that $\theta$ and $\theta'$ both satisfy \eqref{eq: one degree}. Then we have an isomorphism of $G(F)$-representations for $i \in \{0,1\}$,
\[
\rT^i(\pi_{T', \theta' \cdot \epsilon^{\mrm{ram}}[\theta' ]}) \cong \pi_{T, \theta \cdot \epsilon^{\mrm{ram}}[\theta]}^{(\ell)}.
\]
\end{theorem}

\begin{proof}(i) Since $T(F_{\ell})$ is elliptic, we have $T(F_{\ell})G(F_{\ell})_{x,0} = Z(G(F_{\ell})) G(F_{\ell})_{x,0}$. By Lemma \ref{lem: fixed points}, we may then apply Proposition \ref{prop: c-ind} to deduce that
\[
\rT^i \left( \cInd_{T(F_{\ell})G(F_{\ell})_{x,0}}^{G(F_{\ell})}  H_c^{i'}(Y_{\TT_r' ,\UU_r'}; \ol{\F}_{\ell} )_{\theta'}  \right) =    \cInd_{T(F_{\ell})G(F_{\ell})_{x,0}}^{G(F_{\ell})}  \rT^i (H_c^{i'}(Y_{\TT_r' ,\UU_r'}; \ol{\F}_{\ell} )_{\theta'}  ). 
\]
By Theorem \ref{thm: DL tate cohomology}, there exist $i$ and $i'$ for which $\rT^i (H_c^{i'}(Y_{\TT_r' ,\UU_r'}; \ol{\F}_{\ell} )_{\theta'} \neq 0$.

(ii) Let $j$ and $j'$ be the non-vanishing degrees of $H_c^{*}(Y_{\TT_r, \UU_r}; \ol{\F}_{\ell})_{\theta}$ and $H_c^{*}(Y_{\TT_r', \UU_r'}; \ol{\F}_{\ell})_{\theta'}$, respectively. By Lemma \ref{lem: lattice}, the assumptions imply that 
\begin{itemize}
\item $H_c^{*}(Y_{\TT_r', \UU_r'}; \ol{\Z}_{\ell})_{\mrm{tf}, \wt{\theta}'}  \otimes_{\ol{\Z}_{\ell}} \ol{\F}_{\ell} \cong H_c^{*}(Y_{\TT_r', \UU_r'}; \ol{\F}_{\ell})_{\theta'}$ as representations of $T(F_{\ell})G(F_{\ell})_{x,0}$ where $T(F_{\ell})$ acts through $\theta'$ and $G(F_{\ell})_{x,0}$ acts through inflation from $G(F_{\ell})_{x,0:r+}$.
\item $H_c^{*}(Y_{\TT_r, \UU_r}; \ol{\Z}_{\ell})_{\mrm{tf}, \wt{\theta}}  \otimes_{\ol{\Z}_{\ell}} \ol{\F}_{\ell} \cong H_c^{*}(Y_{\TT_r, \UU_r}; \ol{\F}_{\ell})_{\theta}$ as representations of $T(F)G(F)_{x,0}$ where $T(F)$ acts through $\theta$ and $G(F)_{x,0}$ acts through inflation from $G(F)_{x,0:r+}$.
\end{itemize}


By Lemma \ref{lem: fixed points}, we may apply Proposition \ref{prop: c-ind} to deduce that 
\begin{equation}\label{eq: tate coh 1}
\rT^i(\cInd_{T(F_{\ell})G(F_{\ell})_{x,0}}^{G(F_{\ell})} H_c^{j'}(Y_{\TT_r', \UU_r'}; \ol{\F}_{\ell})_{\theta'} )  \cong  \cInd_{T(F)G(F)_{x,0}}^{G(F)}   \rT^i(H_c^{j'}(Y_{\TT_r', \UU_r'}; \ol{\F}_{\ell})_{\theta'}) 
\end{equation}
as $G(F)$-representations. Next, the assumptions that $\theta$ and $\theta'$ satisfy \eqref{eq: one degree} allow us to apply Theorem \ref{thm: DL tate cohomology}(ii) to deduce that 
\[
\rT^i (H_c^{j'}(Y_{\TT_r', \UU_r'}; \ol{\F}_{\ell})_{\theta'})   \cong (H_c^j(Y_{\TT_r, \UU_r}; \ol{\F}_{\ell})_{\theta})^{(\ell)}
\]
as $T(F)G(F)_{x,0}$-representations. We then conclude by noting that Frobenius twist commutes with compact induction.
\end{proof}

\begin{remark}
For a $\sigma$-fixed irreducible representation $\pi'$ of $G(F_{\ell})$, it is typically not obvious that $\rT^i(\pi') \neq 0$. (Thanks to \cite[Theorem 1.3]{F20}, this would already imply that $\pi$ has a base change descent to $G(F)$.) We emphasize that our traction on the toral supercuspidal representations $\pi_{T, \theta \cdot \epsilon^{\mrm{ram}}[\theta]}$ comes from the \emph{geometric} description of these representations developed in \cite{CI19}. Another advantage of the geometric description, which was observed in \cite{CO21}, is that it naturally incorporates the twisting character $\epsilon^{\mrm{ram}}[\theta]$. 
\end{remark}


\begin{corollary}\label{cor: BC for CO} Let assumptions be as in Theorem \ref{thm: tate coh of Chan-Oi}(i),(ii) and suppose $\ell$ is odd and banal for $G$. Then $\pi_{T, \theta \cdot \epsilon^{\mrm{ram}}[\theta]}$ base changes to $\pi_{T', \theta' \cdot \epsilon^{\mrm{ram}}[\theta']}$. 
\end{corollary}

\begin{proof}
Under our assumptions, Proposition \ref{modellirrcusp} implies that  $\pi_{T', \theta' \cdot \epsilon^{\mrm{ram}}[\theta']}$ and $\pi_{T, \theta \cdot \epsilon^{\mrm{ram}}[\theta]}$ are irreducible and cuspidal (noting that $\ell \neq p$ is banal for $G(F)$ if and only if $\ell$ is banal for $G(F_{\ell})$, since the residue field cardinalities of $F_{\ell}$ and of $F$ are congruent modulo $\ell$). The claim then follows from \cite[Theorem 1.3]{F20}, using Theorem \ref{thm: tate coh of Chan-Oi}(ii) to calculate the Tate cohomology of $\pi_{T', \theta' \cdot \epsilon^{\mrm{ram}}[\theta']}$. 
\end{proof}

\begin{remark}
One would expect compatibility between the Genestier-Lafforgue correspondence and Kaletha's correspondence for regular supercuspidal representations \cite{Kal19}. No results towards such compatibility are known at present for general groups. In fact, Kaletha's work has not yet been extended to function fields or to mod $\ell$ representations, but Corollary \ref{cor: BC for CO} appears to be in accordance with what one would expect from such an extension. Namely, if we instead let $F$ be a local field of characteristic \emph{zero} and residue characteristic $p \neq \ell$ sufficiently large relative to $G$, and let $\theta$ be a character $\wt{\theta} \co T(F) \rightarrow \ol{\Q}_{\ell}^{\times}$, then it is computed in \cite[\S 8]{CO21} that the $L$-parameter of $\pi_{T,  \wt{\theta} \cdot \epsilon^{\mrm{ram}}[\wt{\theta}]} $ (according to Kaletha's correspondence for regular supercuspidal representations) is 
\begin{equation}\label{eq: Kaletha L-param}
W_F \xrightarrow{\varphi_{\wt{\theta}}} \ld T(\ol{\Q}_{\ell}) \xrightarrow{\ld j} \ld G(\ol{\Q}_{\ell})
\end{equation}
where $\varphi_{\wt{\theta}}$ corresponds to $\wt{\theta}$ under local class field theory, and $\ld j$ is determined by $T \inj G$. In particular, this implies that $\pi_{T,  \wt{\theta}  \cdot \epsilon^{\mrm{ram}}[\wt{\theta}]}$ base changes to $\pi_{T', \wt{\theta}' \cdot \epsilon^{\mrm{ram}}[\wt{\theta}']}$ under Kaletha's correspondence for regular supercuspidal representations. One would then expect the same of the mod $\ell$ reductions and local function fields $F$, which suggests the statement of Corollary \ref{cor: BC for CO}.
\end{remark}

\end{document}